\documentclass[journal]{IEEEtran}
\usepackage{amsmath,amssymb,amsfonts}
\usepackage{amsthm}
\usepackage{xcolor}
\usepackage{graphicx}
\allowdisplaybreaks

\usepackage{nomencl}
\makenomenclature

\newtheorem{theorem}{Theorem}
\newtheorem{lemma}{Lemma}
\newtheorem{remark}{Remark}
\newtheorem{corollary}{Corollary}
\newtheorem{assumption}{Assumption}
\usepackage{hyperref}
\usepackage[capitalise]{cleveref}
\usepackage{cite}
\usepackage{subfigure}
\ifCLASSINFOpdf
\else
\fi
\usepackage{soul}

\usepackage{caption}
\newcommand{\revise}{}
\newcommand{\ilr}{}
\newcommand{\ilrr}{}
\newcommand{\icr}{}

\newcommand{\icl}{}
\newcommand{\ic}{}
\newcommand{\il}{}
\newcommand{\ml}{}
\newcommand{\mm}{}
\newcommand{\mml}{}
\newcommand{\lmm}{}
\newcommand{\ill}{}

\usepackage[prependcaption,colorinlistoftodos]{todonotes}

\begin{document}
\title{Distributed Optimal Secondary Frequency Control \ic{in Power Networks}
%\icl{Networks: Primal-Dual algorithms}
%and
with \\ \ic{Delay} Independent Stability}

\author{Mengmou Li,~Jeremy Watson,~Ioannis Lestas% <-this % stops a space
%\thanks{Corresponding author: Ioannis Lestas.}
\thanks{%\revise{This work was carried out when all the authors were}
%\ilrr{
The authors are with the Department of Engineering, University of Cambridge, U.K. (e-mail: mmli.research@gmail.com; jeremy.watson@canterbury.ac.nz; icl20@cam.ac.uk)}}% <-this % stops a space


\markboth{Journal of \LaTeX\ Class Files,~Vol.~6, No.~1, January~2007}%
{Shell \MakeLowercase{\textit{et al.}}: Bare Demo of IEEEtran.cls for Journals}



\maketitle

\begin{abstract}
%\boldmath
\icl{Distributed secondary frequency control for power systems, is a problem that has been extensively studied in the literature, and one of its key features is that an additional communication network is required to achieve optimal power allocation.}
%It is thus important to consider the effect of delays as
\icl{%Communication channels are though subject to delays,
Therefore\ml{,} being able to provide stability guarantees in the presence of communication delays is an important requirement.
%The most important distributed controllers for secondary frequency control are the
Primal-dual and distributed averaging proportional-integral (DAPI) protocols, respectively, \ill{are %\ml{the}
two main} control schemes that have been proposed in the literature. Each has its own relative merits, with the former allowing to incorporate general \mml{cost functions} and additional operational constraints, and the latter being more straightforward in its implementation. Although delays have been addressed in DAPI schemes, there are currently no theoretical guarantees for the stability of primal-dual schemes for frequency control, when these are \mm{subject} to communication delays. In fact, simulations illustrate that even small delays can destabilize such schemes. In this paper\ml{,} we show how a novel formulation of primal-dual schemes allows to construct a distributed algorithm with \mm{delay independent} stability guarantees. We also show that this algorithm can incorporate many key features of these schemes such as tie-line power flow requirements, generation constraints\ml{,} and the relaxation of demand measurements with an observer layer. }
%This paper addresses
%therefore proposes controllers with delay independent stability in the primal-dual scheme.
%We show that even small delays can cause instability for the primal-dual scheme.
%To address this challenge, an equivalent passive reformulation is derived from the original primal-dual controller, which eliminates the virtual edge dynamics. The reformulated algorithm is then incorporated with a novel form of scattering transformation to robustify the communication channels against delays.
%Using the proposed algorithm, we also show extensions to tie-line power flow and generation boundedness constraints in the presence of delays, while the requirement for demand measurements can be relaxed by a suitable observer layer.
Finally, we illustrate our results through simulations on a 5-bus example and on the IEEE-39 test system.
\end{abstract}


% \fi
\IEEEpeerreviewmaketitle

\section{Introduction}
\IEEEPARstart{F}{requency} control is an essential task in the operation of power systems. \ic{The development of efficient control policies} for frequency control is becoming increasingly important due to the increasing penetration of %unpredictable
renewable generation \ic{which inevitably exhibits fluctuations in the supply.} %and the reduced inertia of \icl{power systems}.
Frequency control in power systems is normally categorized into three hierarchical layers, namely, primary, secondary, and tertiary control \cite{kundur2007power}.
In this paper, we focus on secondary frequency control in which controllers are designed to restore \ic{the frequency to its} nominal value while maintaining the net area power balance.

\textbf{Literature review:}
Recently, various schemes for distributed secondary frequency control in power systems have been proposed,
%the most important being
\ic{which include} the primal-dual control scheme \cite{li2015connecting,mallada2017optimal,kasis2019stability}, \ic{the} distributed averaging proportional integral (DAPI) control scheme \cite{andreasson2012distributed,zhao2015distributed,trip2016internal,kasis2020distributed}, or a combination of both \cite{yang2020distributed}.
The primal-dual \ic{scheme is} %controllers are
derived from \ic{saddle point formulations by means of Lagrange multipliers} %\mm{the} \icl{Lagrange multiplier methods}
\cite{li2015connecting,mallada2017optimal} and the DAPI \ic{one is} %ones are
derived from distributed proportional-integral controllers in multi-agent systems \cite{andreasson2012distributed,zhao2015distributed}.
%\revise{The \ic{two schemes have various relative trade-offs which include the fact that} the primal-dual scheme allows the consideration of general strictly convex cost functions along with a wider range of constraints but requires generation and demand measurements, while the DAPI scheme requires only local frequency measurements but cannot easily % hardly
%accommodate line and power flow constraints and is restricted to only quadratic cost functions.}
\revise{The two schemes exhibit distinct trade-offs. The primal-dual scheme can handle general strictly convex cost functions along with a wider range of \icr{constraints.
%However, it requires generation and demand measurements.
On the other hand,} the DAPI scheme \icr{has a simpler implementation requiring only local} frequency measurements, but it is limited to quadratic cost functions and cannot \ilrr{incorporate additional operational constraints.}}
%line and power flow constraints while maintaining stability guarantees}.}
Recently, both types of controllers have been extended to consider more complex models of the physical system \cite{trip2016internal,kasis2019stability,kasis2020distributed}.
Despite their differences, the two schemes both include %inherit
a layer of dynamic average consensus dynamics \cite{kia2019tutorial}, which will be \ic{discussed} %revealed
in this paper.

%The physical power system and its controllers constitute a closed-loop cyber-physical system.  In distributed secondary frequency %control,
\icl{\ic{A key feature in distributed secondary control, which affects the structure of the control policies, is the fact that} the frequency is recovered to its nominal value. \ic{This implies that the frequency deviation} can no longer be used as a synchronizing variable through which optimal power sharing is achieved. Therefore, an additional communication network among buses is needed.
%, \icl{because the frequency deviation can no longer be used as a synchronizing variable in %order to achieve an optimal power allocation in a distributed way.}
%in the cyber layer to balance power supply and demand in the physical layer.
%This implies, however that
Communication delays are, however, inevitably present} due to the spatial distribution of buses.
%, which, however, is not considered by the aforementioned literature.
The problem of communication delays in distributed secondary frequency regulations for power systems has been considered in \cite{zhang2014redesigning,schiffer2017robustness,alghamdi2018conditions,koerts2021secondary} for the DAPI scheme.
% \todoiny{I think we need to make our claim more clearly and say in the first sentence that delays have been considered only for DAPI schemes.}
 The work \cite{zhang2014redesigning} shows the robustness of the DAPI control algorithms against arbitrary and bounded constant delays, \cite{koerts2021secondary} considers arbitrary constant delays in the Kron-reduced microgrid, \cite{schiffer2017robustness} derives sufficient conditions for designing parameters for the DAPI algorithms under heterogeneous time-varying delays and \cite{alghamdi2018conditions} extends the results to power systems with second-order turbine governor dynamics.
%However,
\ic{The stability conditions in
\cite{schiffer2017robustness,alghamdi2018conditions} %require explicit knowledge of delays, and verifying stability conditions
involve linear matrix inequalities formulated with global network information.}
\ic{Furthermore,} \il{a key feature of the existing literature} %\cite{zhang2014redesigning,schiffer2017robustness,alghamdi2018conditions,koerts2021secondary}
\il{is
the fact that %all address
delays have been addressed in the DAPI scheme only,} and there \ic{\il{are currently no %existing
results} in the literature relating} to delays in the primal-dual control scheme \icl{\ill{for secondary} frequency control}.
%This is possibly due to technical difficulties in addressing delays in the communication channel involving the incidence matrix in the `virtual swing equations', and proper characterization of delays is hampered by the existence of virtual edge dynamics.

\ic{In fact,} the primal-dual control algorithms \ic{in their conventional formulations} are very sensitive to communication delays. %this paper, we \icl{demonstrate}%reveal
%via simulations that
\icl{Even} \mm{small delays} could destabilize \il{power systems implementing such schemes} or cause the controllers to fail to restore the frequency \ic{to its nominal value}, \icl{which is %something that is
\ic{also} demonstrated in this paper \ic{in various simulation examples}.}
As \il{previously mentioned}, the primal-dual controllers are preferable to the DAPI ones when more \ic{advanced operational specifications} are considered such as general convex cost functions, generation boundedness constraints, and tie-line power flow constraints \cite{mallada2017optimal,kasis2019stability}.
Therefore, addressing delay issues in the primal-dual scheme is \icl{a significant problem \ic{of practical relevance}.}
%of great importance.

\il{Providing delay robustness guarantees in primal-dual schemes for optimal secondary frequency control is in general an involved problem. This is, for example, reflected in the lack of such results in the existing literature, and the sensitivity of existing such schemes to delays.
%As it will be discussed in the paper,
%On the one hand, existing primal-dual based schemes can lead to problematic
%behaviour even when small delays are introduced.
Furthermore,
the scattering transform \cite{hokayem2006bilateral,li2020smooth}, \cite{chopra2006passivity,hatanaka2018passivity}, which is a protocol that can lead to  delay independent stability, and has been applied to DAPI schemes \cite{koerts2021secondary}, is not directly applicable to conventional implementations of primal-dual schemes for secondary frequency control. This is due to the presence of "virtual edge dynamics", as it will be discussed in more detail within the paper. Therefore\ml{,} novel formulations of primal-dual control policies are needed in order to achieve robustness to communication delays, which is one of the contributions of this paper.}

\il{In particular, we show that an appropriate novel reformulation of  primal-dual schemes for secondary control allows to design distributed control policies with delay independent stability guarantees. Furthermore, we show that this reformulation allows to incorporate the operational constraints associated with generation and power flows that primal-dual schemes can handle.}

%\il{Another complication associated with the study of delays in power system models, arises form the fact that power system stability is often analyzed by means of energy like Lyapuonv functions and the use of invariance principles due to the higher order dynamics present. Nevertheless, the infinite dimensional setting relevant for the analysis of delays complicates the use of invariance principles since boundedness of trajectories does not directly follow from level sets of energy like functionals.}
%
%The scattering \ic{transformation, initially used in bilateral teleoperation, is a useful tool to} guarantee delay independent stability \cite{hokayem2006bilateral,li2020smooth}. It has been successfully incorporated into passive network systems to achieve synchronization \cite{chopra2006passivity,hatanaka2018passivity}.
%%In this paper, we aim to construct controllers independent of constant delays using %\icl{the} scattering transformation.
%However, unlike its application to the DAPI scheme \cite{koerts2021secondary},
%%this technique is not readily applicable to the primal-dual controlled power networks.
%%\icl{In particular, a main complication arises from
%\icl{the presence of `edge dynamics' in `virtual swing equations' hinders a direct application to primal-dual schemes for frequency control.} %, \ic{which is also reflected in the fact that such results do not exist in the literature.}}

\il{The analysis in the paper is also of independent interest making use of appropriately constructed Lyapunov functionals, and an invariance principle. %.} % \ic{that }
%Invariance principles are often used in power system stability analysis in the undelayed case, due to %the higher order dynamics involved, however
Despite the significance of invariance principles in power system stability analysis,
their application in an infinite dimensional setting is more involved.
%The latter is more involved to apply due to the infinite dimensional and differential-algebraic %character of the problem.
In particular, the level sets of energy like Lyapunov functionals do not necessarily lead to boundedness of trajectories, thus requiring a further exploitation of the structure of power system dynamics.}

%\ic{One of the main contributions of this paper is to show that an appropriate novel reformulation of  primal-dual schemes allows to design distributed control policies with delay independent stability. Furthermore, we show that this reformulation allows to incorporate the operational constraints associated with generation and power flows that primal-dual schemes can handle.
%%Therefore conventional primal-dual schemes need to be appropriately reformulated \mm{with a novel form of scattering transformation} in order to achieve delay independent stability, which is one of the main contributions of this paper.}
%\ic{The analysis in the paper is also of independent interest making use of appropriately constructed Lyapunov functionals, and an invariance principle. %.} % \ic{that }
%The latter is more involved to apply due to the infinite dimensional and differential-algebraic character of the problem, and requires an exploitation of the structure of the dynamics considered.}
%%. In particular convergence is deduced via an invariance principle,}
%%Due to the infinite dimensional and differential algebraic aspect of the problem considered as the analysis in the paper is also of independent interest as the boundedness of trajectories, a necessary condition to be able to apply an invariance principle, does not follow from conventional Lyapunov functionals, and the structure of the power system dynamics needs to be exploited.
\ic{The main contributions are outlined below.}
%\lmm{The novel scattering transformation also introduces new variables whose boundedness cannot be shown by the Lyapunov %functional, rendering a non-trivial application of the invaraince principle to the delay differential algebraic equation.}
%On \icl{the} one hand, the existence of the edge dynamics in the `virtual swing equations' hinders the proper characterization of time delays.
%On the other hand, the scattering transformation needs to be carefully designed to guarantee synchronization for the specific distributed control protocol. \icl{In particular, one of the main contributions of the paper is to show that by reformulating the primal-dual dynamics to a multivariable distributed communication protocol where additional states get communicated, allows to design distributed algorithms for secondary control with delay independent stability guarantees.}

\textbf{Contributions:}
%We address heterogeneous communication delays for the primal-dual algorithms by integrating passivity-based techniques with the power system model and controllers that inherit passivity properties.
\begin{enumerate}
	\item \icl{\ic{We propose, for the first time,}} distributed \icl{control policies} with delay independent stability \icl{guarantees} for \icl{primal-dual \ic{based}} secondary frequency control.
	\item We show that, similar to the conventional primal-dual controllers, the proposed \icl{control scheme allows to incorporate various additional constraints such as %extra
tie-line power flow constraints and constraints on \il{generation. It also} allows to relax the requirement for demand measurements with an observer}.
\end{enumerate}

\textbf{Paper Organization:} Basic notation and preliminaries are given in \Cref{sec:Preliminaries}. The power system model and optimization problem to be considered are introduced in \Cref{sec:Power System Model}. \il{The primal-dual control algorithm is then} introduced and \icl{the proposed control policy} with delay independent stability is \icl{presented} in \Cref{sec:Primal-dual Scheme}. \il{Extensions} of the primal-dual controller are given in \Cref{sec:Extensions} to account for operational constraints such as tie-line power flow and generation bounds, and to relax the requirement of explicit knowledge of the demand.
In \Cref{sec:Simulations}, we demonstrate our results through simulations.
Finally, conclusions are drawn in \Cref{sec:Conclusion}. \icl{The proofs} of the main results are given in the Appendix.

\section{Preliminaries}\label{sec:Preliminaries}
The \icl{notation} used in this paper \icl{is} summarized in \Cref{table:1}.
Given a group of vectors $x_1,\ldots, x_N$, \icl{we use} $x$ without subscripts \icl{to denote} their aggregates \il{unless} specified otherwise, i.e., $x = (x_1^T, \ldots, x_N^T)^T$. \icl{Also for $x\in\mathbb{R}^n$, $\|x\|$ denotes its Euclidean norm.}
\begin{table}[h!]
\centering
\begin{tabular}{l  l}
\hline
Indices\\
\hline
$\mathbf{0}$ & all-zero vector or matrix of proper size\\
$\mathbf{1}_n$ & $n$-dimensional all-one vector\\
$I_n$ &  $n \times n$ identity matrix\\
$\mathcal{G}$ & graph index\\
$D$ & Incidence matrix of the physical transmission network\\
$\tilde{D}$ & Incidence matrix of the communication graph\\
$\tilde{L}$ & Laplacian matrix of the communication graph\\
$\tilde{L}_{{K}}$ & Laplacian of communication graph without inter-area lines\\
$x^*$ & equilibrium point of variable $x$\\
\hline
Sets\\
\hline
$\mathbb{R}$ & set of real numbers\\
$\mathbb{R}^{n}$ & set of $n$-dimensional real vectors\\
$\mathbb{R}^{n \times m}$ & set of $n \times m$ real matrices\\
$G$ & set of generation buses\\
$L$ & set of load buses\\
$N$ & set of buses satisfying $N = G \cup L$\\
$E$ & set of physical transmission lines\\
$\tilde{E}$ & set of communication lines\\
$\tilde{N}_j$ & set of buses that have direct communication with bus $j$\\
${K}$ & set of communication areas\\
${C}_{k}$ & set of buses in the control area $k$\\
${B}_{k}$ & set of physical lines that connect area $k$ to other areas\\
$\tilde{B}$ & set of communication lines connecting different areas\\
\hline
Variables\\
\hline
$\eta_{ij}$ & power angle difference between bus $i$ and bus $j$\\
$\omega_j$ & frequency deviation from the nominal frequency at bus $j$\\
$Y_{ij}$ & line susceptance between bus $i$ and bus $j$\\
$M_j$ & generator inertia at bus $j$\\
$\Lambda_j$ & frequency damping coefficient at bus $j$\\
$p_j^L$ & \revise{uncontrollable demand} at bus $j$\\
$p_j^M$ & mechanical power injection at bus $j$\\
$p_{ij}$ & power transfer from bus $i$ to bus $j$\\
$p_j^c$ & power command signal at bus $j$\\
$r_{ij}^{x}$ & variable that bus $j$ receives from bus $i$ about variable $x$\\
$u_j$ & control input to generation bus $j$\\
$\hat{P}_{k}$ & scheduled tie-line power flow for area $k$\\
\hline
\end{tabular}
\caption{Main symbols used in this paper.}
\label{table:1}
\end{table}

The power network model is described by an \il{undirected, connected} graph $\mathcal{G}(N, E)$ where $N = \{ 1, \ldots, |N| \}$ is the set of buses and $E \subseteq N \times N$ the set of edges representing the transmission lines connecting the buses.
Since generators have inertia, it is reasonable to assume that only buses with inertia have non-trivial
generation dynamics. We define $G = \{1,\ldots,|G|\}$ and $L = \{1,\ldots, |L|\}$ as the sets of buses with and without inertia, respectively, such that $|G| + |L| = |N|$.
The edge $(i, j)$ denotes the link connecting buses $i$ and $j$.
For each $j \in N$, we use $i : i \rightarrow j$ and $k : j \rightarrow k$ to denote the sets of buses that precede and succeed bus $j$ respectively.
%The graph
%$\mathcal{G}(N, E)$ is assumed to be of an arbitrary direction.
%, so that if $(i, j) \in E$ then $(j, i) \notin E$.
\ic{We define} the directed incidence matrix $D \in \mathbb{R}^{|N| \times |E|}$ such that the element $D_{ij} = -1$ if the edge $j$ leaves node $i$, $D_{ij} = 1$ if the edge $j$ enters node $i$
and $0$ otherwise.
It should be noted that the form of the power system dynamics is not affected by
the ordering of nodes, and our results are independent of the choice of direction.
In addition to the \il{power network,} %\il{(referred to as the physical layer)},
we define \il{a} communication network described by a connected graph $\mathcal{G}(N, \tilde{E})$,  %\il{(referred to as the cyber layer)},
\il{and its} incidence matrix $\tilde{D}$ is defined similarly.
Moreover, $\tilde{N}_j$ denotes the neighboring set for bus $j$ in the communication graph such that $i \in \tilde{N}_j$ if either $(i,j) \in \tilde{E}$ or $(j,i) \in \tilde{E}$. We also define the Laplacian matrix for the communication graph as $\tilde{L} = \tilde{D} W \tilde{D}^T$ where $W \in \mathbb{R}^{|\tilde{E}| \times |\tilde{E}|}$ is a positive diagonal matrix representing edge weights. Then, we have $\tilde{L} \geq 0$.

Let ${K} = \{ 1, \ldots, |{K}|\}$ be the set of all control areas in the network. Let ${C}_{k}$ \il{denote} the set of buses in the $k^{th}$ control area, \il{which} satisfies $N = {C}_{1} \cup \ldots \cup {C}_{|{K}|}$, \mm{$C_{i} \cap C_{j} = \emptyset$, for all $ i \neq j$}. Define ${B} \subseteq {E}$, $\tilde{B} \subseteq \tilde{E}$ as the sets of physical lines and communication lines that connect different control areas, respectively. Let ${B}_k \subset {B}$ be the set of boundary lines for area $k$.
Let $\mathcal{G}(N,\tilde{E} / \tilde{B})$ be the subgraph of $\mathcal{G}(N,\tilde{E})$ by deleting all boundary lines connecting different areas.
Let $\tilde{L}_{{K}}$ be the Laplacian matrix for $\mathcal{G}(N,\tilde{E} / \tilde{B})$. It also holds that $\tilde{L}_{{K}} \geq 0$.
Moreover, we assume that the subgraph $\mathcal{G}(N,\tilde{E} / \tilde{B})$ has $|K|$ connected components.

%We will investigate the solutions of systems of differential equations with arbitrary bounded time delays. See \cite{hale2013introduction} for a detailed study of such systems.
%\todoiny{I think the sentence above should be a footnote or omitted.}
\icl{Delay differential equations.} Let $\mathcal{C} \left( [ -r, 0], \icl{\mathbb{R}^n} \right)$ denote the Banach space of continuous functions mapping $\left[-r, 0\right] \subseteq \mathbb{R}$ into $\icl{\mathbb{R}^n}$, with the norm of an element $\varphi$ in $\mathcal{C}$ given by $\|\varphi\| = \sup_{ \revise{- \il{r}} \leq x \leq 0} \| \varphi(x) \|$.
Let $x_t$ denote the function in $\mathcal{C} \left( [ -r, 0], \icl{\mathbb{R}^n} \right)$ given by $x_{t}(\theta) = x(t + \theta)$ for $\theta \in [-r, 0]$ and $t \in [0, \infty )$.
A general delay differential equation can be written as
$\dot{x} (t) = f (t,  x_t )$,
where $f: \mathbb{R}_{+} \times \mathcal{C} \left( [ -r, 0], \icl{\mathbb{R}^n}  \right)  \rightarrow \mathbb{R}^{n}$ is continuous in its first argument and locally Lipschitz, uniformly in $t$, in its second argument, which guarantee the existence and uniqueness of solutions and their continuous dependence on the initial condition \cite{hale2013introduction}.
A delay differential system
\begin{align}\label{eq:del}
\dot{x} =  f(\mm{t}, x_t, u), \quad y = h(x, u)
\end{align}
%\todoiny{\st{give the definition of function $h$ and the dimensions of $u, y$. Also do we need here $h$ to be a functional of $x_t$ or is being a function of $x$, as stated above, sufficient for the analysis in the paper, when passivity is invoked?}}
where $f: \mm{\mathbb{R}_{+}} \times \mathcal{X} \times \mathcal{U} \rightarrow \mathbb{R}^{n}$, $h: \mathbb{R}^{n} \times \mathcal{U} \rightarrow \mathbb{R}^{m}$, $\icl{\mathcal{X}} \subseteq \mathcal{C} \left( [ -r, 0], \icl{\mathbb{R}^n}\right)$, $\icl{\mathcal{U} \subseteq \mathbb{R}^m}$,
with equilibrium point $(\icl{x_t = \mathbf{0}},~ u = \mathbf{0})$
is said to be \mm{locally passive} if there exist an open neighborhood of $\icl{(x_t = \mathbf{0} ,~u = \mathbf{0} )} \in \icl{\mathcal{X} \times \mathcal{U}}$, and
a continuously differentiable positive semidefinite functional \icl{$V(x_t)$} such that
\begin{align*}
\dot{V} \icl{(x_t)} \leq u^T y, \quad \text{\icl{for all} }u \in \mathcal{U},~ \icl{x_t \in \mathcal{X}}.
\end{align*}
It is passive if the above inequality holds globally, \icl{i.e. for all $u \in \mathbb{R}^m,~ \icl{x_t \in \mathcal{C}}$.}
\icl{Note that the definition of passivity used}  for \icl{the time-delayed system \eqref{eq:del}} is \icl{analogous} %similarly
to that \icl{for undelayed} systems \icl{(e.g. \cite{khalil1996nonlinear})}, but a `storage functional' %storage functional
is used instead of a storage function.

\section{Power System Model}\label{sec:Power System Model}
We make the following assumptions for the network:
\begin{enumerate}
	\item Bus voltage magnitudes are $|V_j| = 1$ p.u. for all $j \in N$.
	\item Lines $(i, j) \in E$ are lossless and characterized by their susceptances $Y_{ij} = Y_{ji} > 0$.
	\item Reactive power flows do not affect bus voltage phase angles and frequencies.
\end{enumerate}
These assumptions are generally valid at medium to high voltages and are standard in the analysis of secondary frequency control \cite{arthur2000power}.

The power system model is described by the swing \icr{equation} at generation buses \revise{\eqref{eq:power system model eta}, \eqref{eq:power system model generator power balance}} and power balance at load buses \revise{\eqref{eq:power system model power balance}, \icr{and is given as follows,}}
\begin{subequations}\label{eq:power system model}
\begin{align}
& \hspace{-2mm} \dot{\eta}_{ij} = \omega_i - \omega_j, ~ (i, j) \in E,\label{eq:power system model eta}\\
& \hspace{-3mm} M_j \dot{\omega}_j = - p_j^L + p_j^M  \hspace{-1mm} - \hspace{-1mm}  \Lambda_j \omega_j - \hspace{-2mm} \sum_{k: j \rightarrow k} p_{jk} + \hspace{-2mm} \sum_{i: i \rightarrow j} p _{ij}, ~ j \in G,\label{eq:power system model generator power balance}\\
& \hspace{-2mm} 0 = -p_j^L - \Lambda_j \omega_j - \hspace{-2mm} \sum_{k: j \rightarrow k} p_{jk} + \hspace{-1.5mm} \sum_{i: i \rightarrow j} p _{ij}, ~j \in L,\label{eq:power system model power balance}\\
&  \hspace{-2mm} p_{ij} = Y_{ij} \sin \eta_{ij}, ~ (i,j) \in E \label{eq:power system model edge dynamics}
\end{align}
\end{subequations}
where \revise{the variables are defined in \Cref{table:1}.}
In particular, the positive constants $M_j$ and $\Lambda_j$ represent the generator inertia at generation bus $j$ and the frequency damping coefficient at any bus $j$, respectively, $p_j^L$ denotes the frequency-independent uncontrollable load at bus $j$, \revise{which \icr{could include a step change} in the demand.}

For simplicity, we consider first-order generation dynamics given by
\begin{align}\label{eq:generation dynamics}
\begin{array}{rl}
	\tau_j \dot{p}_j^M = - p_j^M + k_{g,j} u_j, \quad j \in G
	\end{array}
\end{align}
for some constants $\tau_j > 0$ and $k_{g,j} > 0$. The generation input $u_j$ is
%will be generated by
an appropriate control policy to be \icl{designed.}
%System \eqref{eq:generation dynamics} implies that the %re always exists an asymptotically %stable
%generation output $\bar{p}_j^M$ \icl{is always asymptotically stable} given any constant %input $u_j = \bar{u}_j$.
\begin{remark}
\revise{\ilrr{The first-order system considered here facilitates the passivity analysis in the paper. The dissipativity condition in  \cite{kasis2019stability} can be used when higher order turbine-governor dynamics are present, however this extension is
%The first-order system is considered here to effectively conduct passivity analysis. It is possible %to consider higher-order turbine governor dynamics which becomes passive when combined with the %damping term \cite{kasis2016primary,kasis2019stability}},
%\il{though these are
omitted} for simplicity in the presentation}.
%\todoiny{Can the dissipativity property in \cite{kasis2019stability} be handled or is it just the passivity property of generation dynamics?}
 \il{Also for brevity} in the presentation, we do not consider controllable loads on the demand side, which may be included in $p_j^M$.
	\end{remark}
%\begin{equation}
%	\begin{cases}
%		\dot{x}_j^M = A_j x_j^M + B_j u_j,\\
%		p_j^M = C_j x_j^M + D_j u_j,
%	\end{cases}
%	, ~ j \in G
%\end{equation}

It is desired that the generation is adjusted to match the uncontrollable demand with minimal cost. This goal can be represented by an optimization problem, which is termed the optimal generation regulation (OGR) problem:
\begin{align}\label{optimal generation regulation problem}
	\begin{array}{rl}
	\textbf{OGR:} \quad & \underset{p^M}{\min} \displaystyle \sum_{j \in G} Q_j (p_j^M),\\
	& \text{subject to }  \displaystyle \sum_{j \in G} p_j^M = \displaystyle \sum_{j \in N} p_j^L,
	\end{array}
\end{align}
where $Q_j(\cdot)$ is the \ic{cost function associated with generation} at bus $j$ and the constraint represents power \ic{balance.} % between supply and demand.
For \ic{the feasibility of this problem}, the following assumption is \ic{made.} %considered.

\begin{assumption}\label{assumption convex cost function}
The cost functions $Q_j$, $j \in G$ are continuously differentiable and strictly convex, and optimization problems considered in this work are feasible, i.e., there exists an equilibrium point of the system \eqref{eq:power system model}, \eqref{eq:generation dynamics} \ilr{with $\omega=0$} such that the corresponding $p^M$ is a solution to the considered optimization problem.
\end{assumption}
\revise{This assumption implies that \icr{power balance} %between the generation $p^M$ and load $p^L$
\ilr{can be satisfied at each bus after a disturbance, for power allocations that are solutions to the optimization problem in \eqref{optimal generation regulation problem}.}}
%, which is a necessary requirement} for secondary frequency control.}

\begin{figure}[htbp]
	\centering
	\includegraphics[width = 0.6\linewidth, height = 0.6\linewidth]{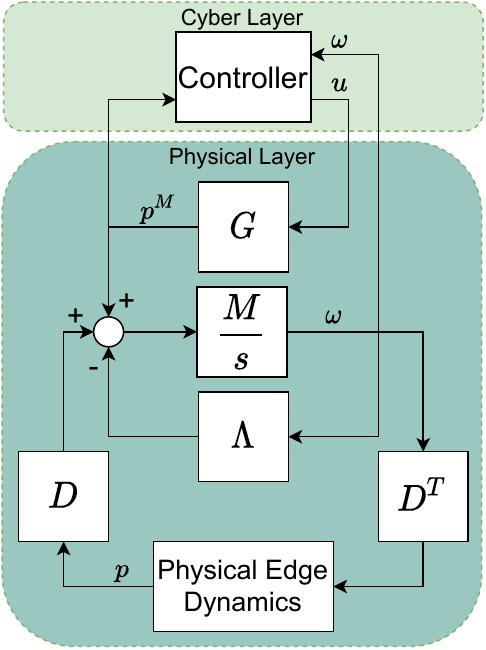}
	\caption{Schematic overview of the power system model. The physical edge dynamics are described by \revise{\eqref{eq:power system model eta} and} \eqref{eq:power system model edge dynamics}. The component $G$ represents the generation dynamics \eqref{eq:generation dynamics} and the controller is to be designed.  \icl{Though not explicitly illustrated in the figure, it should be noted that the cyber layer is required to be based on a distributed %\mm{a}
	communication protocol so as to implement a distributed secondary frequency control policy}.}
	\label{fig:Schematic overview}
\end{figure}

The physical layer \icl{(power system dynamics)} and cyber layer \icl{(control policy)} of the power network \il{are} depicted in \cref{fig:Schematic overview}, where $G$ denotes the generation dynamics and the physical edge dynamics are represented by \revise{\eqref{eq:power system model eta} and }\eqref{eq:power system model edge dynamics}; \icl{$u$ is} the generation input to be designed.
The goal in \il{distributed}
secondary frequency regulation is to design a \icl{controller} %in the cyber layer
such that the \textup{OGR} problem \eqref{optimal generation regulation problem} is solved in a distributed manner and the frequency is restored to its nominal value. It should be noted that the cyber layer should \icl{be based on a distributed %\mm{a}
communication protocol} among buses \icl{so as to implement a distributed secondary frequency control policy}. The equilibrium of \eqref{eq:power system model}, \eqref{eq:generation dynamics} \icl{satisfies}
\begin{subequations}
	\begin{align}
		& 0 = \omega_i^* - \omega_j^*, ~ (i, j) \in E, \label{eq:equilibrium of omega}\\
		& 0 = - p_j^L + p_j^{M,*} - \Lambda_j \omega_j^* - \hspace{-2mm} \sum_{k: j \rightarrow k} p_{jk}^* + \hspace{-1.5mm} \sum_{i: i \rightarrow j} p _{ij}^*, ~ j \in G, \label{eq:equilibrium of swing equations}\\
		&0 = -p_j^L - \Lambda_j \omega_j^* - \hspace{-2mm} \sum_{k: j \rightarrow k} p_{jk}^* + \hspace{-1.5mm} \sum_{i: i \rightarrow j} p _{ij}^*, ~j \in L, \label{eq:equilibrium of load balance}\\
		& p_{ij}^* = Y_{ij} \sin \eta_{ij}^*, ~ (i,j) \in E,\\
		& p_j^{M,*} = k_{g,j} u_j^*, ~ j \in G \label{eq:equilibrium of generation dynamics}
	\end{align}
\end{subequations}
with $p_j^{M,*}$ being the optimal solution of the \textup{OGR} problem \eqref{optimal generation regulation problem}.

We adopt the following assumption that is widely used in the power systems literature.
\begin{assumption}\label{assumption angle}
	$|\eta_{ij}^*| < \frac{\pi}{2}$ for all $(i,j) \in E$.
\end{assumption}
This assumption can be regarded as a security constraint that generally holds under normal operating conditions.

\section{Primal-dual Scheme With Communication Delays}\label{sec:Primal-dual Scheme}
In this section, we first review the primal-dual secondary frequency control algorithm for solving the \textup{OGR} problem \eqref{optimal generation regulation problem}. Then, we show that \ic{the} classical scheme is incapable of %addressing
\icl{efficiently incorporating} communication delays.  Next, we give an equivalent passive reformation of the primal-dual scheme \ic{that involves communicating an additional state. This is a key feature that allows to combine the control policy with the scattering transformation so as to achieve} delay independent stability.
\subsection{Controllers}
\icl{Equations \eqref{eq:generation input primal-dual}, \eqref{eq:controller without delays primal-dual} below describe the primal-dual scheme for optimal secondary frequency control that has been proposed in the literature\revise{\cite{mallada2017optimal,kasis2019stability}}. In particular,} for the generation dynamics \eqref{eq:generation dynamics}, the generation input is given by
\begin{align}\label{eq:generation input primal-dual}
\begin{array}{rl}
	u_j = k_{c,j} \left( p_j^c - \omega_j \right) +  \frac{p_j^M}{k_{g,j}} - k_{c,j} Q_j'\left(p_j^M\right), ~ j \in G
\end{array}
\end{align}
where the parameter $k_{g,j}$ is given in \eqref{eq:generation dynamics}, $k_{c,j} > 0$,  \revise{$p_j^c$ is a control variable \icr{through which %to schedule
power sharing is achieved, and is referred to} as the power command signal}, and $Q_j'$ represents the gradient of $Q_j$ for bus $j$.
The dynamics for the power command signal $p_j^c$ \revise{can be found in the literature \cite[Eqn (18)]{mallada2017optimal} \cite[Eqn (6)]{kasis2019stability}},
\begin{subequations}\label{eq:controller without delays primal-dual}
	\begin{align}
		\hspace{-2mm} \gamma_{ij} \dot{\psi}_{ij} = & p_i^c - p_j^c, ~ (i,j) \in \tilde{E} \label{eq:controller without delays edge dynamics primal-dual} \\
		\gamma_{j} \dot{p}_j^c = & - \left( p_j^M - p_j^L \right) - \sum_{k: j \rightarrow k} \psi_{jk} + \sum_{i: i\rightarrow j} \psi_{ij}, ~ j \in N \label{eq:controller without delays bus dynamics primal-dual}
	\end{align}
\end{subequations}
where $p_j^M := 0$, for $j \notin G$, $\gamma_j$ and $\gamma_{ij}$ are positive constants, $\psi_{ij}$ is a state of the controller that integrates the power command difference of communicating buses $i$ and $j$.
Controller \eqref{eq:controller without delays primal-dual} is referred to as \il{`virtual swing equation'} \cite{mallada2017optimal} since it has a similar structure \ic{to that of} the system model \eqref{eq:power system model}. \ic{Equation} \eqref{eq:controller without delays edge dynamics primal-dual} \ic{can similarly be seen as representing 'virtual edge dynamics'.} The set of communication lines $\tilde{E}$ here can be either the same or different to the set of physical transmission lines $E$.

The %optimality and
\il{convergence to an optimal equilibrium point using the control policy in} \eqref{eq:generation input primal-dual}, \eqref{eq:controller without delays primal-dual} can be found in the literature but we include it here \icl{to facilitate the derivation of} subsequent results and for completeness.

\begin{lemma}[Optimality]\label{lem:optimality without delays primal-dual}
Let Assumption \ref{assumption convex cost function} hold.
Any equilibrium of system \eqref{eq:power system model}, \eqref{eq:generation dynamics}, \eqref{eq:generation input primal-dual}, \eqref{eq:controller without delays primal-dual} is an optimal solution to the \textup{OGR} problem \eqref{optimal generation regulation problem} with $\omega^* = \mathbf{0}_{|N|}$ and $p^{c.*} = \text{Im}(\mathbf{1}_{|N|})$.
\end{lemma}
\mm{The proof can be found in \cite{kasis2019stability}.}
%The proof is given in Appendix\ref{appendix proof of lemma optimality without delays}.

\begin{lemma}[Convergence]\label{lem:convergence without delays primal-dual}
Consider an equilibrium of \eqref{eq:power system model}, \eqref{eq:generation dynamics}, \eqref{eq:generation input primal-dual}, \eqref{eq:controller without delays primal-dual} in which Assumption \ref{assumption angle} holds. Then, there exists an open neighborhood about the equilibrium such that solutions of \eqref{eq:power system model}, \eqref{eq:generation dynamics}, \eqref{eq:generation input primal-dual}, \eqref{eq:controller without delays primal-dual} asymptotically converge to a set of equilibria that solve the \textup{OGR} problem \eqref{optimal generation regulation problem} with $\omega^* = \mathbf{0}_{|N|}$.
\end{lemma}
\mm{\ic{A sketch} of the proof} is given in Appendix\ref{appendix proof of lem convergence without delays primal-dual}.

\subsection{Equivalent Reformulation of the Primal-dual Control}\label{subsection:Equivalent Reformulation of the Primal-dual Control}
Notice that the \il{`virtual swing equation' \eqref{eq:controller without delays primal-dual} contains} both the virtual bus dynamics \eqref{eq:controller without delays bus dynamics primal-dual} and the virtual edge dynamics \eqref{eq:controller without delays edge dynamics primal-dual}.
In practice, the virtual edge dynamics \eqref{eq:controller without delays edge dynamics primal-dual} are \ic{implemented at} each \icl{bus; that is,} each bus possesses and updates the dynamics of its corresponding edges. As a result, there is redundant information, e.g., the bus $j$ possesses the edge information of $\psi_{ij}$, denoted by $\psi_{ij}^{j}$, while its neighboring bus $i$ possesses the same edge information, denoted by $\psi_{ij}^{i}$. This scheme works satisfactorily in the undelayed case where $\psi_{ij}^{j} = \psi_{ij}^{i}$.
However, when there are heterogeneous delays in the communication channel, the updates of the same edge dynamics in the two buses become
\begin{equation}\label{eq:redundant update under delays}
\begin{aligned}
	\gamma_{ij} \dot{\psi}_{ij}^{j} = p_i^c(t - T_{ij}) - p_j^c, ~ (i,j) \in \tilde{E}\\
	\gamma_{ij} \dot{\psi}_{ij}^{i} = p_i^c - p_j^c ( t - T_{ji}), ~ (j, i) \in \tilde{E}
\end{aligned}
\end{equation}
where $T_{ij}$ and $T_{ji}$ represents the communication delays in the channel $i \rightarrow j$ and $j \rightarrow i$, respectively. It is apparent that $\psi_{ij}^{j} \neq \psi_{ij}^{i}$, and the goal of the secondary frequency control is not guaranteed, as \ic{also illustrated in the} simulations in \Cref{sec:Simulations}.
Even if we assume in addition that $T_{ij} = T_{ji}$, it still requires the additional knowledge of the delays to update the state with self-induced delays
\begin{equation}
	\gamma_{ij} \dot{\psi}_{ij}^{j} = p_i^c (t - T_{ij}) - p_j^c ( t - T_{ij}), ~ (i,j) \in \tilde{E}
\end{equation}
such that the two variables ${\psi}_{ij}^{i}$, ${\psi}_{ij}^{j}$ are rendered equal.
%These requirements are in fact restrictive.
\ic{Moreover, the} system becomes unstable when the homogeneous delay is large under this protocol \cite{schiffer2017robustness}.
%One may try to introduce extra computing units for the virtual edge dynamics, which increases the communication burden and is also not economical in a large-scale power network with numerous communication lines.

To ease these restrictions and deal with unknown and heterogeneous delays, we reformulate the primal-dual control algorithm \eqref{eq:controller without delays primal-dual} into an equivalent form in which it becomes \ic{possible} %more convenient
to address delays using the %with the tool of
scattering transformation.
Let $\gamma_{ij} = \gamma_{j} = 1$ in \eqref{eq:controller without delays primal-dual} and the \ml{communication} Laplacian weight \ml{in this work be $W = I$ for the ease of presentation, i.e., $\tilde{L} = \tilde{D} \tilde{D}^T $}, \revise{where $\tilde{D}$ is the incidence matrix.}
% \todoing{which Laplacian are you referring to here?}.
We can write \eqref{eq:controller without delays primal-dual} into the compact form
\begin{align}\label{eq:controller without delays compact primal-dual}
		\dot{\psi} =  - \tilde{D}^T p^c, \quad \dot{p}^c = - \left( p^M - p^L \right) + \tilde{D} \psi
\end{align}
where $\psi \in \mathbb{R}^{|\tilde{E}|}$ is the aggregated vector of $\psi_{ij}$, $\forall (i,j) \in \tilde{E}$.
Denote $\xi = \tilde{D} \psi$; \ic{since} $\tilde{D} \tilde{D}^T = \tilde{L}$ for undirected communication graphs, \eqref{eq:controller without delays compact primal-dual} becomes
\begin{align}\label{eq:controller without delays reformulation-1 primal-dual}
	\dot{\xi} =  - \tilde{L} p^c,
	\quad	\dot{p}^c =  - \left( p^M - p^L \right) + \xi
\end{align}
with initial condition satisfying $ \mathbf{1}_{|N|}^T \xi (0) = 0$ due to $\mathbf{1}_{|N|}^T \tilde{D} \psi \equiv 0$.
Since $\xi \in \mathbb{R}^{|N|}$, the virtual edge dynamics \ic{are} thus eliminated and each bus only needs to possess and update the dynamics of $\xi_j$ by communicating the variable $p_j^c$.
This reformulation also implies that the primal-dual control \eqref{eq:controller without delays reformulation-1 primal-dual} in fact inherits a layer of dynamic average consensus dynamics \cite{kia2019tutorial}.
However, simulations in  \Cref{sec:Simulations} show that \eqref{eq:controller without delays reformulation-1 primal-dual} is \mm{problematic} under \il{delays} because the delayed system will converge to undesirable equilibrium points.
\revise{This is because \ilrr{in order to have $\mathbf{1}_{|N|}^T (p^M - p^L)=0$} at the equilibrium point \ilrr{(which \icr{implies} $\omega=0$ from \eqref{eq:power system model generator power balance})} %requires that
\icr{we need $\mathbf{1}_{|N|}^T \xi = 0$ at this point}, which \icr{can be} violated if delays are introduced.}
% \todoiny{does this mean we do not have an equilibrium point?or we do not have power balance at equilibrium?}

To solve this problem, we introduce a coordinate transformation again to reformulate the controller.
Notice that the Laplacian $\tilde{L}$ for undirected graphs is positive semidefinite and thus there exists a unique square root $\tilde{L}^{\frac{1}{2}} \geq 0$ such that $\tilde{L} = \tilde{L}^{\frac{1}{2}} \tilde{L}^{\frac{1}{2}}$.
Define a new variable $\zeta$ such that $ \dot{\zeta} =  - \tilde{L}^{\frac{1}{2}} p ^c$. Then, \icl{we have} that $\xi = \tilde{L}^{\frac{1}{2}} \zeta$ and \eqref{eq:controller without delays reformulation-1 primal-dual} becomes
\begin{align}\label{eq:controller without delays reformulation-2 primal-dual}
	\dot{\zeta} =  - \tilde{L}^{\frac{1}{2}}  p^c , \quad
		\dot{p}^c =  - \left( p^M - p^L \right) + \tilde{L}^{\frac{1}{2}}\zeta.  %\label{eq:controller without delays reformulation-2-2 primal-dual}
\end{align}
%To show that \eqref{eq:controller without delays reformulation-2 primal-dual} is equivalent to the primal-dual control algorithm \eqref{eq:controller without delays compact primal-dual}, integrating \eqref{eq:controller without delays reformulation-2-1 primal-dual} and substituting it into $\tilde{L}^{\frac{1}{2}} \zeta$, we obtain that
%$
%	\tilde{L}^{\frac{1}{2}} \int_{0}^{t} \tilde{L}^{\frac{1}{2}} p^{c}(\tau) d\tau
%	= \int_{0}^{t} \tilde{L} p^{c}(\tau) d\tau = \tilde{D}\int_{0}^{t} \tilde{D}^T p^{c}(\tau) d\tau = \tilde{D} \psi
%$, which recovers \eqref{eq:controller without delays compact primal-dual}.
Since $\tilde{L}^{\frac{1}{2}}$ is also a well-defined Laplacian matrix, \il{we
%change $\tilde{L}^{\frac{1}{2}}$ to $\tilde{L}$ and
consider} the following new controller for each bus,
\begin{subequations}\label{eq:controller without delays reformulation primal-dual}
\begin{align}
	\dot{\zeta}_j = & \sum_{i \in \tilde{N}_j} \alpha_{ij}\left( p_i^c - p_j^c \right), ~ j \in N \label{eq:controller without delays reformulation primal-dual 1}\\
		\dot{p}_j^c = & - \left( p_j^M - p_j^L \right) - \sum_{i \in \tilde{N}_j} \alpha_{ij} \left( \zeta_i - \zeta_j \right), ~ j \in N \label{eq:controller without delays reformulation primal-dual 2}
\end{align}
\end{subequations}
where $\alpha_{ij} > 0$ and $\alpha_{ij} = \alpha_{ji}$. \revise{Under this controller,  \ilrr{the \icr{desired} relation $\mathbf{1}_{|N|}^T (p^M - p^L)=0$} %power balance
at the equilibrium point is not affected by delays.}
The only difference between \eqref{eq:controller without delays reformulation primal-dual} and \eqref{eq:controller without delays reformulation-2 primal-dual} is in the Laplacian matrices, which do not affect the analysis results as long as the communication graph is connected.
Therefore, we have the following corollary.
\begin{corollary}\label{cor:convergence primal-dual}
	Let Assumption~\ref{assumption convex cost function} hold and consider an equilibrium of \eqref{eq:power system model}, \eqref{eq:generation dynamics}, \eqref{eq:generation input primal-dual}, \eqref{eq:controller without delays reformulation primal-dual} in which Assumption \ref{assumption angle} holds. Then, there exists an open neighborhood about the equilibrium such that solutions of \eqref{eq:power system model}, \eqref{eq:generation dynamics}, \eqref{eq:generation input primal-dual}, \eqref{eq:controller without delays reformulation primal-dual} \icl{converge} to a set of equilibria that solve the \textup{OGR} problem \eqref{optimal generation regulation problem} with $\omega^* = \mathbf{0}_{|N|}$.
\end{corollary}
The proof is given in Appendix\ref{appendix proof of cor convergence primal-dual}.

\begin{figure}[bhtp]
	\centering
	\includegraphics[width = 1\linewidth]{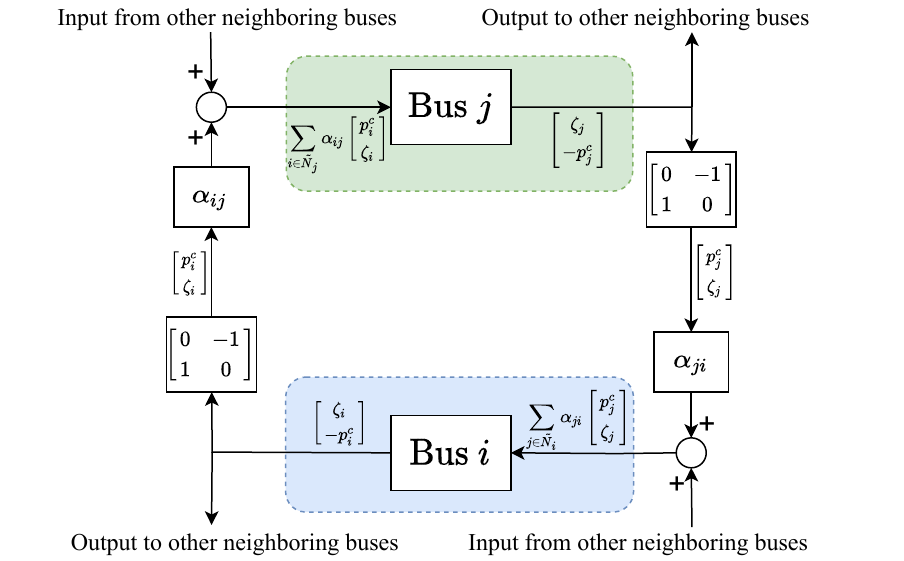}
	\caption{The communication scheme between two buses in the equivalent reformulation \eqref{eq:controller without delays reformulation primal-dual} of the primal-dual control algorithm. The buses are also \icl{coupled} in the physical layer, which is not shown here.}
	\label{fig:communication_scheme_primal-dual}
\end{figure}

\icl{As previously discussed, the} original controller \eqref{eq:controller without delays primal-dual} is incapable of addressing delays \icl{efficiently. The significance of the equivalent reformulation \eqref{eq:controller without delays reformulation primal-dual} that has been derived % mm is that
is that it allows \lmm{to} construct controllers with delay independent stability properties, as it will be shown in the next section.}
The communication scheme between two buses in \eqref{eq:controller without delays reformulation primal-dual} is depicted by \cref{fig:communication_scheme_primal-dual}. It represents the distributed control among buses in the cyber layer of \cref{fig:Schematic overview}. Note that the buses are also physically coupled in the physical layer but we only include the cyber layer in \cref{fig:communication_scheme_primal-dual} to highlight the \icl{structure of the} communication scheme.
\begin{remark}
Compared with \eqref{eq:controller without delays primal-dual}, the reformulations \eqref{eq:controller without delays reformulation-1 primal-dual} and \eqref{eq:controller without delays reformulation primal-dual} are both node-based and do not contain virtual edge dynamics. Thus, no redundant update is carried out as in \eqref{eq:redundant update under delays}.
Compared with \eqref{eq:controller without delays reformulation-1 primal-dual}, \ic{the} reformulation \eqref{eq:controller without delays reformulation primal-dual} requires communication of \icl{the} extra variable $\zeta$, but satisfies \icl{appropriate} passivity properties, which \icl{will be shown} in the next subsection.
The \icl{scheme \eqref{eq:controller without delays reformulation-1 primal-dual}} also reveals the existence of a dynamic average consensus control layer \cite{kia2019tutorial} underneath the primal-dual scheme. From these forms we can easily observe the difference between the DAPI algorithms \cite{zhao2015distributed,trip2016internal,schiffer2017robustness,kasis2020distributed} and the primal-dual ones, i.e., the former applies the dynamic averaging consensus control directly to the generation input (Equation (32) in \cite{zhao2015distributed}) while the latter inherits a layer of dynamic average consensus control \icl{within higher-order dynamics}.
\end{remark}

\subsection{Communication Delays and Scattering Transformation}
\ic{In this subsection we show how %consider the presence of unknown and heterogeneous delays  in
the primal-dual secondary control algorithm described in the previous section can be adapted so as to have convergence guarantees when arbitrary heterogeneous delays are present.}
%We address the unknown and heterogeneous delays in the primal-dual secondary control algorithm in %this subsection.
Suppose that there exist unknown and heterogeneous constant delays in the communication channels between buses. The delay in the communication channel $i \rightarrow j$ is denoted by $T_{ij}$ and the delay in the communication channel $j \rightarrow i$ is denoted by $T_{ji}$, for $(i,j) \in \tilde{E}$.

When there are communication delays in the channel $i \rightarrow j$, $i \in \tilde{N}_j$, bus $j$ receives delayed information of $p_i^c$, $\zeta_i$ from bus~$i$.  Then, controller \eqref{eq:controller without delays reformulation primal-dual} becomes
\begin{align}\label{eq:controller under delays reformulation primal-dual}
\begin{array}{rl}
	\hspace{-3mm}	\dot{\zeta}_j = \hspace{-2mm} & \displaystyle \sum_{i \in \tilde{N}_j} \alpha_{ij}\left( p_i^c (t - T_{ij}) - p_j^c \right), ~ j \in N\\
	\hspace{-3mm}	\dot{p}_j^c =\hspace{-2mm}  & \hspace{-1mm} -  \hspace{-1mm} \left( p_j^M \hspace{-1mm} - \hspace{-0.5mm} p_j^L \right) \hspace{-0.5mm} - \hspace{-1mm} \displaystyle \sum_{i \in \tilde{N}_j} \hspace{-0.5mm} \alpha_{ij} \left( \zeta_i (t - T_{ij}) - \zeta_j \right), ~ j \in N
\end{array}
\end{align}
 which may affect the stability of the power network if delays are large. In fact, the primal-dual controllers are very sensitive to communication delays. Simulations in \Cref{sec:Simulations} show that even some small $T_{ij}$ could \mm{destabilize} the system.

To deal with communication delays, \ic{we start by modifying} \eqref{eq:controller without delays reformulation primal-dual} as
\begin{subequations}\label{eq:control with scattering transformation primal-dual}
	\begin{align}
	\dot{\rho}_j^{\zeta} = & - \rho_j^{\zeta} + \sum_{i \in \tilde{N}_j} \alpha_{ij}\left( r_{ij}^{p} - p_j^c \right), ~ j \in N \label{eq:control with scattering transformation primal-dual 1} \\
	\dot{\zeta}_j = & - \rho_j^{\zeta} + 2 \sum_{i \in \tilde{N}_j} \alpha_{ij}\left( r_{ij}^{p} - p_j^c \right), ~ j \in N \label{eq:control with scattering transformation primal-dual 2}\\
	\dot{\rho}_j^{p} = & \ml{- \rho_j^{p}  - \left( p_j^M \hspace{-0.5mm} - \hspace{-0.5mm} p_j^L \right)  \hspace{-0.5mm} - \hspace{-1mm} \sum_{i \in \tilde{N}_j} \alpha_{ij} \left( r_{ij}^{\zeta}  \hspace{-0.5mm} - \hspace{-0.5mm} \zeta_j \right), ~ j \in N }\label{eq:control with scattering transformation primal-dual 3}\\
	\dot{p}_j^c = &\ml{ - \rho_j^{p}  \hspace{-0.5mm} - \hspace{-0.5mm} 2  \hspace{-0.5mm} \left( p_j^M - p_j^L \right)  \hspace{-0.5mm} - 2 \hspace{-1mm} \sum_{i \in \tilde{N}_j} \hspace{-0.5mm} \alpha_{ij}  \hspace{-0.5mm} \left( r_{ij}^{\zeta}  \hspace{-0.5mm} -  \hspace{-0.5mm}\zeta_j \right)  \hspace{-0.5mm}, ~  \hspace{-0.5mm} j \in N } \label{eq:control with scattering transformation primal-dual 4}
	\end{align}
\end{subequations}
%\todoing{is the first term in the RHS of \eqref{eq:control with scattering transformation primal-dual 2} $\rho^\zeta_j?$}
where $\rho_j^{\ml{\zeta}}$, \ml{$\rho_j^{p}$} are auxiliary states for bus $j$,
\revise{\icl{$r_{ij}^{p}$, $r_{ij}^{\zeta}$  denote %the alternative representations of the
\il{information} bus $j$ receives from bus $i$. }} This information will be appropriately formulated via a communication protocol that will be designed in the rest of this section. Equations \eqref{eq:control with scattering transformation primal-dual 1}, \eqref{eq:control with scattering transformation primal-dual 2} and \ml{\eqref{eq:control with scattering transformation primal-dual 3}, \eqref{eq:control with scattering transformation primal-dual 4}} result from the addition of parallel feedforward compensators to \eqref{eq:controller without delays reformulation primal-dual 1}, \ml{\eqref{eq:controller without delays reformulation primal-dual 2}}, respectively.
%to be designed with bus $i$.
It will \icl{also be shown later within the paper,  %shown later
that} the parallel feedforward compensator provides excess of passivity to guarantee convergence and does not affect the equilibrium since $\rho_j^{\ml{\zeta}}, \ml{\rho_j^{p}}\rightarrow 0$ recovers the original system \cite{li2022parallel}.

\icl{Let $\zeta_j^*$ be an equilibrium value of $\zeta_j$ and similarly let $r_{ij}^{p,*}, p_{i}^{c,*}, r_{ij}^{\zeta,*}, \zeta_i^*$ denote also equilibrium values that satisfy $r_{ij}^{p,*} = p_{i}^{c,*}$, $r_{ij}^{\zeta,*} = \zeta_i^*$.}
%\todoiny{need to define the equilibrium values of all the states as well -  the locality of the passivity property is about an equilibrium point for the states and inputs, rather than on only inputs and outputs.}
%We have the following lemma on system passivity.
\revise{To facilitate stability analysis under delays, we will utilize an interconnection of passive components. In this regard, we present the following lemma on the passivity of the first component.}
\begin{lemma}\label{lem:passivity of reformulated controller}
	The system described by \eqref{eq:power system model}, \eqref{eq:generation dynamics}, \eqref{eq:generation input primal-dual}, \eqref{eq:control with scattering transformation primal-dual} is locally passive\footnote{Strictly speaking, it is incrementally passive, which is \icl{a property that is} independent of \icl{the} equilibrium point. The locality of the result is due to the sinusoidal functions in the system model rather than the nonlinear cost functions. The analytical result becomes global if the sinusoids are linearized.} with respect to input
	$\begin{bmatrix}
		\tilde{r}^{p} & \tilde{r}^{\zeta}
	\end{bmatrix}^T$ and output
	$\begin{bmatrix}
		\tilde{\zeta} & - \tilde{p}^c
	\end{bmatrix}^T$, where the $j^{th}$ element of $\begin{bmatrix}
		\tilde{r}^{p} & \tilde{r}^{\zeta}\end{bmatrix}^T$ and $\begin{bmatrix}
		\tilde{\zeta} & - \tilde{p}^c
	\end{bmatrix}^T$ are defined by
	$
	\begin{bmatrix}
		\tilde{r}_j^p \\ \tilde{r}_j^{\zeta}
	\end{bmatrix} =
		\sum_{i \in \tilde{N}_j} \alpha_{ij} \left(
		\begin{bmatrix}
		r_{ij}^{p} \\ r_{ij}^{\zeta}
	\end{bmatrix} -
	\begin{bmatrix}
		r_{ij}^{p,*} \\ r_{ij}^{\zeta,*}
	\end{bmatrix}
	\right)
	=
	\sum_{i \in \tilde{N}_j} \alpha_{ij}
		\begin{bmatrix}
		\tilde{r}_{ij}^{p} \\ \tilde{r}_{ij}^{\zeta}
	\end{bmatrix}
	$ and
	$\begin{bmatrix}
		\tilde{\zeta}_j \\ - \tilde{p}_j^c
	\end{bmatrix} =
	\begin{bmatrix}
		\zeta_j \\ - p_j^c
	\end{bmatrix} -
	\begin{bmatrix}
		\zeta_j^* \\ - p_j^{c,*}
	\end{bmatrix},$
	 respectively.
%	 Each bus is locally passive from $\begin{bmatrix}
%		\tilde{r}_j^p \\ \tilde{r}_j^{\zeta}
%	\end{bmatrix}$ to $\begin{bmatrix}
%		\tilde{\zeta}_j \\ - \tilde{p}_j^c
%	\end{bmatrix}$.
\end{lemma}
The proof is given in Appendix\ref{appendix proof of lem passivity of reformulated controller}.

%\todoiny{is the term locally passive an established one - clarify in the preliminaries?\\
%-Yes, it has been clarified.}
%\begin{assumption}\label{assumption delays}

\begin{figure}[bhtp]
	\centering
	\includegraphics[width = 1\linewidth, height = 1\linewidth]{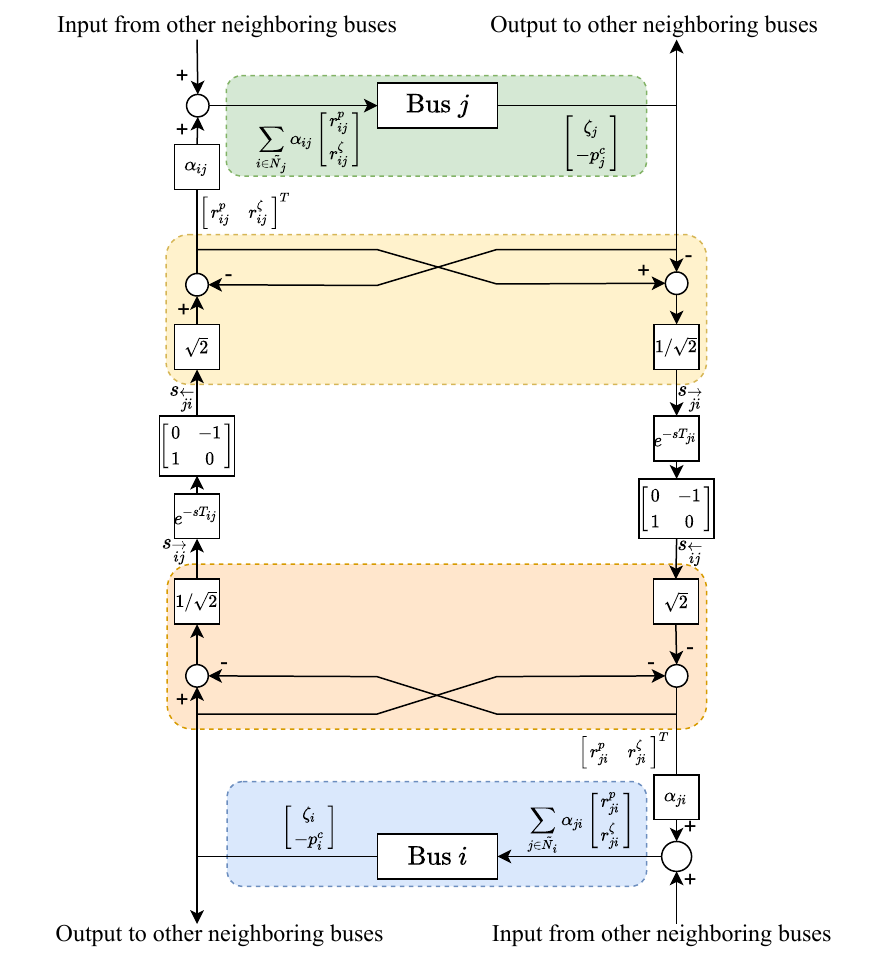}
	\caption{The scattering transformation between two buses under communication delays for the primal-dual scheme. The scattering transformation is passive with respect to the two-port inputs and outputs, as shown in \Cref{lem:passivity of scattering transformation}.}
%	\todoing{top left, above $1/\sqrt{2}$ block: is the minus meant to be a plus?} I have checked it, all the signs should be correct.}
	\label{fig:diagram of scattering_primal-dual}
\end{figure}

To robustify the communication channels against delays, \il{we send} `encoded' information of $\zeta$, $p^c$ instead of their direct information, and then `decode' the information received to obtain $\zeta$, $p^c$.
To this end, we adopt the following scattering \icl{transformation}
%to extract information of
%\icl{to formulate the signals}
%$r_{ij}^{p}$ and $r_{ij}^{\zeta}$,
\begin{equation}\label{eq:scattering transformation primal-dual}
	\begin{aligned}
& s_{\overrightarrow{ij}} = \frac{- 1}{\sqrt{2}}\left( \begin{bmatrix}
		r_{ji}^{p} \\ r_{ji}^{\zeta}
	\end{bmatrix}  - \begin{bmatrix}
		\zeta_i \\ - p_i^c
	\end{bmatrix} \right),
s_{\overleftarrow{ij}} = \frac{ - 1}{\sqrt{2}} \left( \begin{bmatrix}
		r_{ji}^{p} \\ r_{ji}^{\zeta}
	\end{bmatrix}  +\begin{bmatrix}
		\zeta_i \\ - p_i^c
	\end{bmatrix} \right)\\
& s_{\overleftarrow{ji}} = \frac{1}{\sqrt{2}} \left( \begin{bmatrix}
		r_{ij}^{p} \\ r_{ij}^{\zeta}
	\end{bmatrix}  + \begin{bmatrix}
		\zeta_j \\ - p_j^c
	\end{bmatrix} \right),~
s_{\overrightarrow{ji}} = \frac{1}{\sqrt{2}} \left( \begin{bmatrix}
		r_{ij}^{p} \\ r_{ij}^{\zeta}
	\end{bmatrix}  - \begin{bmatrix}
		\zeta_j \\ - p_j^c
	\end{bmatrix} \right)
	\end{aligned}
\end{equation}
where $s_{\overrightarrow{ij}}$ is the scattering variable that bus $i$ sends to bus $j$, and $s_{\overleftarrow{ij}}$ is the scattering variable that bus $i$ \ic{receives} from bus $j$. The other scattering variables above are defined similarly. Let
$E_s = \begin{bmatrix}
	0 & -1\\ 1 & 0
\end{bmatrix}$ be \il{a matrix gain that is applied on} the transmitted scattering variables. %Since
\il{Noting that}
 ($s_{\overrightarrow{ij}}$, $s_{\overleftarrow{ji}}$) and ($s_{\overrightarrow{ji}}$, $s_{\overleftarrow{ij}}$) are the variables in the communication channels $i \rightarrow j$ and $j \rightarrow i$, respectively, it \ic{holds} that
\begin{equation}\label{eq:scattering variables under delays primal-dual}
	s_{\overleftarrow{ij}}(t) = E_s s_{\overrightarrow{ji}}(t - T_{ji}), \quad s_{\overleftarrow{ji}}(t) = E_s s_{\overrightarrow{ij}}(t - T_{ij}).
\end{equation}
where $T_{ji}$, $T_{ij}$ are delays in the communication channels $j \rightarrow i$ and $i \rightarrow j$, respectively, $s_{\overrightarrow{ji}} (t) = \mathbf{0}$ and $s_{\overrightarrow{ij}} (t) = \mathbf{0}$, $\forall t < 0$.
The communication process from bus $j$ to $i$ is summarized as follows.
\begin{enumerate}
\item Bus $j$ encodes its original input $\begin{bmatrix}
		r_{ij}^{p} \\ r_{ij}^{\zeta}
	\end{bmatrix}$ and output $\begin{bmatrix}
		\zeta_j \\ - p_j^c
	\end{bmatrix}$ into the scattering variable $s_{\overrightarrow{ji}}$ based on \eqref{eq:scattering transformation primal-dual}.
\item Scattering variable $s_{\overrightarrow{ji}}$ is transmitted in the communication channel $j \rightarrow i$ under delay $T_{ji}$. Then bus $i$ \icl{receives} $s_{\overleftarrow{ij}}(t)$ given by  \eqref{eq:scattering variables under delays primal-dual}.
\item Bus $i$ decodes the variable $s_{\overleftarrow{ij}}$ and extracts original input $\begin{bmatrix}
		r_{ji}^{p} \\ r_{ji}^{\zeta}
	\end{bmatrix}$ based on \eqref{eq:scattering transformation primal-dual}.
\end{enumerate}
The \ic{communication} from bus $i$ to $j$ is carried out similarly. This communication scheme is depicted in \cref{fig:diagram of scattering_primal-dual}.
If there are no delays, i.e., $T_{ij} = T_{ji} = 0$, we immediately obtain $r_{ij}^{p} = p_i^c$, $r_{ij}^{\zeta} = \zeta_i$, for all $(i,j) \in \tilde{E}$ since \eqref{eq:scattering transformation primal-dual}, \eqref{eq:scattering variables under delays primal-dual} are equivalent to loop transformations in the undelayed case.
\begin{remark}\label{rm:eliminating algebraic loop}
Observe that there exist algebraic loops in the scattering transformation as original inputs and outputs are coupled to construct new inputs and outputs. Eliminating the algebraic loops we obtain
\begin{align}\label{eq: r_ij eliminating algebraic loop}
 \begin{bmatrix} {r}_{ij}^{p} (t) \\  {r}_{ij}^{\zeta} (t) \end{bmatrix} = &
 -  \begin{bmatrix}   {r}_{ij}^{p} ( t - T_{ij} - T_{ji} ) \\ {r}_{ij}^{\zeta} (t - T_{ij} - T_{ji} ) \end{bmatrix} +
 \begin{bmatrix} {\zeta}_j (t - T_{ij} - T_{ji}) \\ - {p}_j^{c} (t - T_{ij} - T_{ji}) \end{bmatrix}
\nonumber \\
& - \begin{bmatrix}  {\zeta}_j (t) \\ - {p}_j^{c} (t)  \end{bmatrix}
 +  2 \begin{bmatrix}  {p}_i^{c} (t - T_{ij}) \\ {\zeta}_i (t - T_{ij}) \end{bmatrix}.
\end{align}
%The term $\begin{bmatrix}   {r}_{ij}^{p} ( t - T_{ij} - T_{ji} ) \\ {r}_{ij}^{\zeta} (t - T_{ij} - T_{ji} ) \end{bmatrix}$ can be eliminated by iterations.
\mm{The system consisting of \eqref{eq:power system model}, \eqref{eq:generation dynamics}, \eqref{eq:generation input primal-dual}, \eqref{eq:control with scattering transformation primal-dual}, \eqref{eq: r_ij eliminating algebraic loop} becomes a delay differential algebraic equation \cite{koerts2021secondary}.}
%\mml{The closed-loop system can also be written in the form $\dot{x} (t) = f (t, x)$ where $f$ is globally Lipschitz on $x$ and continuous on $t$, given any initial conditions that uniquely define $r^{p}(t)$, $r^{\zeta} (t)$.}
\ic{Note that solutions trivially exist and are unique for a given initial condition} %Then, we have that the solution \mm{exists and} is
%uniquely defined by its initial conditions
\cite{khalil1996nonlinear,hale2013introduction}.
\end{remark}
%\todoiny{\st{I am not sure what "iteration" means in the remark above. Perhaps omit the remark as due to the delays it seems we actually do not have algebraic loops?}}

\begin{remark}
Our work differs from the application of scattering transformation to distributed optimization \cite{hatanaka2018passivity},  \icl{due to the presence of the physical power system dynamics that couple the individual buses, in addition to the coupling arising from the communication between the controllers.}
%since \icl{the system considered involves a} coupling of \icl{the} physical and cyber %layers \icl{within the} power network.
%\icl{As discussed within the This also requires an equivalent passive reformulation %beforehand as discussed within the paper}.
%In addition,
\revise{
The scattering transformation depicted in \cref{fig:diagram of scattering_primal-dual} %considers the %transmission
\icr{involves the communication}
of two variables and an extra matrix gain $E_{s}$ is added in \eqref{eq:scattering variables under delays primal-dual}, compared to \icl{the way it is} commonly used in the literature \cite{hokayem2006bilateral,chopra2006passivity,hatanaka2018passivity,li2020smooth}.} The matrix $E_{s}$ \icl{appears} due to the special structure of interconnection in \cref{fig:communication_scheme_primal-dual}.
%Thus, the scattering operator considered here contains the matrix $E_{s}$ in addition to delay terms, %which
\il{This additional matrix does} not disrupt \ic{the} passivity of \ic{the system arising from} the scattering transformation \il{(\Cref{lem:passivity of scattering transformation} below)}  since its norm is not greater than one \cite{hokayem2006bilateral}.
This %\icl{formulation} %
\icl{\ilrr{formulation of the scattering transformation}} %with \icl{a} matrix gain,
\icl{is more broadly} applicable \icl{to} systems having similar distributed control structures, such as the \icl{extensions presented} in \Cref{sec:Extensions}.
\end{remark}

%The constructed scattering transformation is passive with respect to its two-port inputs and outputs as stated in the following
%\il{lemma.}
\revise{The constructed scattering transformation is the second component of the interconnection of passive components, which is passive with respect to its two-port inputs and outputs as stated in the following lemma.}
% In the lemma variables with $\tilde{.}$ denote deviations from equilibrium values as in \Cref{lem:passivity of reformulated controller}.}
%\lmm{Slightly different from \Cref{sec:Preliminaries}, the passivity here is defined with respect to deviations about an %equilibrium point instead of the origin.}
%\subsection{Passivity Properties}
\begin{lemma}\label{lem:passivity of scattering transformation}
	The system \icl{described by} %consisting of
\eqref{eq:scattering transformation primal-dual} and \eqref{eq:scattering variables under delays primal-dual} is passive with respect to input $\begin{bmatrix}
		\begin{bmatrix} \tilde{\zeta}_j & - \tilde{p}_j^c \end{bmatrix} & \begin{bmatrix} \tilde{\zeta}_i & - \tilde{p}_i^c \end{bmatrix}
	\end{bmatrix}^T$
and output
$ - \begin{bmatrix}
		\begin{bmatrix} \tilde{r}_{ij}^{p} & \tilde{r}_{ij}^{\zeta} \end{bmatrix} & \begin{bmatrix} \tilde{r}_{ji}^{p} & \tilde{r}_{ji}^{\zeta} \end{bmatrix}
	\end{bmatrix}^T$, \lmm{where \il{these} variables
%\il{$\tilde{r}_{ij}^{p}, \tilde{r}_{ij}^{\zeta}$}
are \il{as defined} in \Cref{lem:passivity of reformulated controller}.}
\end{lemma}
%\todoiny{We should state that the passivity property is for deviations about an equilibrium point. In the introduction passivity is defined w.r.t. the origin.}
The proof is given in Appendix\ref{appendix proof of lem passivity of scattering transformation}.

\subsection{Convergence and Optimality}
After showing passivity of the component \eqref{eq:power system model}, \eqref{eq:generation dynamics}, \eqref{eq:generation input primal-dual}, \eqref{eq:control with scattering transformation primal-dual} \revise{via \Cref{lem:passivity of reformulated controller}}, and the component \eqref{eq:scattering transformation primal-dual}, \eqref{eq:scattering variables under delays primal-dual} \revise{\icr{via \Cref{lem:passivity of scattering transformation}}}, \icl{respectively,} we are ready to \icl{prove} %present
asymptotic stability of the closed-loop system under \icl{arbitrary constant} delays.

\begin{theorem}\label{thm:convergence under delays primal-dual}
Let Assumption~\ref{assumption convex cost function} hold,  and consider an equilibrium \il{point} of system \eqref{eq:power system model}, \eqref{eq:generation dynamics}, \eqref{eq:generation input primal-dual}, \eqref{eq:control with scattering transformation primal-dual}, \eqref{eq:scattering transformation primal-dual}, \eqref{eq:scattering variables under delays primal-dual} in which Assumption \ref{assumption angle} holds.
%Suppose the constant
\icl{The delays $T_{ij}\geq0$, $\forall (i, j) \in \tilde{E}$ are assumed to be constant, and are allowed to take arbitrary bounded values and be heterogeneous}.
Then, there exists an open neighborhood $\Omega$ about this equilibrium \il{point} such that the solutions of \eqref{eq:power system model}, \eqref{eq:generation dynamics}, \eqref{eq:generation input primal-dual}, \eqref{eq:control with scattering transformation primal-dual}, \eqref{eq:scattering transformation primal-dual}, \eqref{eq:scattering variables under delays primal-dual} with initial conditions in
%$\mathcal{C} \left( \left[ -\max_{\forall i,j}\{T_{ij}\}, 0 \right], \Omega \right)$,
\icl{$\Omega$
%asymptotically
converge} \il{to an equilibrium point that solves the} %to the optimal solution of the
\textup{OGR} problem \eqref{optimal generation regulation problem} with $\omega^* = \mathbf{0}_{|N|}$.
\end{theorem}
%The proof is given in Appendix\ref{appendix proof of thm convergence under delays primal-dual}.
\mml{The \ic{proof
%involves a nontrivial application of the invariance principle to the delay differential algebraic %equation, which
is} given in Appendix\ref{appendix proof of thm convergence under delays primal-dual}.}

\section{Extensions}\label{sec:Extensions}
As previously mentioned, the conventional primal-dual controllers can be extended to consider more complex scenarios such as tie-line power flow constraints \cite{mallada2017optimal} and generation boundedness constraints \cite{kasis2019stability}.
\ic{We show in this section that the proposed delay} independent primal-dual controllers \icl{can \ic{also} adopt such additional requirements.}
%have good extensibility as well.
%In the following subsections,
\ic{In particular, we} consider \ic{additional} tie-line power flow constraints, generation boundedness constraints, and then relax the requirement of demand measurements using an observer layer.

\subsection{Tie-line Power Flow Constraints}\label{subsec:Tie-line Power Flow Constraints}
Recall that ${C}_{k}$ is the set of buses in the control area \icl{$k\in K$.} The following \textup{OGR-2} problem considers additional tie-line power flow constraints,
 \begin{align}\label{optimal generation regulation problem 2}
\begin{array}{rl}
	\hspace{-2mm}\textbf{OGR-2:} \quad & \underset{p^M, p}{\min} ~ \displaystyle \sum_{j \in G} Q_j (p_j^M),\\
	& \text{subject to } \displaystyle  \sum_{j \in G} p_j^M = \sum_{j \in N} p_j^L,\\
	& \qquad \qquad \displaystyle  \sum_{(i,j) \in {B}_{k}} \hat{D}_{k,ij} p_{ij} = \hat{P}_{k}, ~ k \in {K}
\end{array}
\end{align}
where $\hat{P}_k$ is the net power injection of area $k$,  ${B}_{k}$ is the set of physical lines that connect area $k$ to other areas, and $\hat{D}_{k,ij} = 1$ if $i \in C_{k}$, $\hat{D}_{k,ij} = -1$ if $j \in {C}_{k}$, and $\hat{D}_{k,ij} = 0$, otherwise.
The second constraint in \eqref{optimal generation regulation problem 2} specifies power transfer between areas.
Moreover, let $\hat{D} = [\hat{D}_{k,ij}] \in \mathbb{R}^{|{K}| \times |E|}$
%\todoiny{do you mean $E$ or $\tilde E$}
and $E_{{K}} = [e_{1} ~\ldots~e_{|K|}]^T \in \mathbb{R}^{|{K}| \times |N|}$, where $e_{k} \in \mathbb{R}^{|N|}$, $k \in K$ is a vector with elements $(e_{k})_{j} = 1$ if $j \in C_{k}$ and $(e_{k})_{j} =0$ otherwise.
Then, $\hat{D}$ can be related with the incidence matrix ${D}$ by \icl{(\hspace{1sp}\cite{mallada2017optimal})}
\begin{align}\label{eq:D=E_K D}
	\hat{D} = E_{{K}} {D}.
\end{align}
%Assume that problem \eqref{optimal generation regulation problem 2} is feasible, i.e., there exists an equilibrium point of the system \eqref{eq:power system model}, \eqref{eq:generation dynamics} such that \eqref{optimal generation regulation problem 2} is solved.

Note that solving the OGR-2 problem \eqref{optimal generation regulation problem 2} is challenging since we cannot manipulate the variables $p_j^M$ or $p_{ij}$ directly but can only rely on controllers \ic{for the} generation input $u$ to alter their values.
Since we are aiming \ic{in this work} to design controllers with delay independent \ic{stability,} we adopt the reformulated scheme in \Cref{sec:Primal-dual Scheme} and propose an extension of the primal-dual controller to solve \eqref{optimal generation regulation problem 2}.
For bus $j$ in the $k^{th}$ control area, we replace the controller \eqref{eq:controller without delays reformulation primal-dual} with
\begin{subequations}\label{eq:controller without delays tie-line primal-dual}
\begin{align}
	\dot{\zeta}_j = & \sum_{i \in \tilde{N}_j} \alpha_{ij}\left( p_i^c - p_j^c \right) - \hspace{-1mm} \sum_{i \in \tilde{N}_j} \mm{\alpha_{ij}} \left( \pi_i - \pi_j \right),\label{eq:controller without delays tie-line primal-dual 1}\\[-0.5mm]
		\dot{p}_j^c = & - \left( p_j^M - p_j^L \right) - \sum_{i \in \tilde{N}_j} \alpha_{ij} \left( \zeta_i - \zeta_j \right), \label{eq:controller without delays tie-line primal-dual 2}\\[-0.5mm]
		\dot{\pi}_j = & \sum_{i \in \tilde{N}_j } \hspace{-1mm} \mm{\alpha_{ij}} \left( \zeta_i \hspace{-0.5mm} - \hspace{-0.5mm} \zeta_j \right) -  \hspace{-3mm} \sum_{i \in \tilde{N}_j \cap {C}_{k}}  \hspace{-3mm} \mm{\alpha_{ij}} \left(\phi_i - \phi_j \right) \hspace{-0.5mm} - \hspace{-0.5mm} J_j \hat{P}_{k}, \label{eq:controller without delays tie-line primal-dual 3}\\[-0.5mm]
		\dot{\phi}_j = & \sum_{i \in \tilde{N}_j \cap {C}_k} \hspace{-2mm} \mm{\alpha_{ij}} \left( \pi_i - \pi_j \right), \label{eq:controller without delays tie-line primal-dual 4}
\end{align}
\end{subequations}
where $\pi_j$, $\phi$ are auxiliary variables associated with the tie-line power flow constraints, and $J_j \in \{0, 1 \}$, \ic{is a} constant representing the knowledge of $\hat{P}_{k}$. \icl{In particular, we} assume that in the area $k$  only one bus $j_{k}$ has access to the value of $\hat{P}_{k}$ such that $J_j = 1$ if $j = j_k$ and $J_j = 0$, otherwise.
In \eqref{eq:controller without delays tie-line primal-dual 3} and \eqref{eq:controller without delays tie-line primal-dual 4}, the \icl{use of the set} $\tilde{N}_j \cap {C}_k$ implies that variables $\phi_j$ and $\pi_j$ are exchanged only within the same control area.
\begin{remark}
Compared with \cite{mallada2017optimal}, controller \eqref{eq:controller without delays tie-line primal-dual} is proposed without centralized control or virtual edge dynamics such that communication delays can be similarly addressed without the obstacles discussed in \Cref{subsection:Equivalent Reformulation of the Primal-dual Control}. Similar controllers can be found in \cite{yang2020distributed}, \icl{where an undelayed scheme is analyzed that} relaxes the requirement \icl{for} demand measurements. Compared to \cite{yang2020distributed}, our proposed controller \eqref{eq:controller without delays tie-line primal-dual} \icl{does not  need to verify conditions that involve global parameters to guarantee stability, and can also be used to achieve stability for arbitrary communication delays as it will be shown in this section.}
%is fully distributed without the need to verify global stability conditions.
\end{remark}

Let $\tilde{L}_{{K}}$ be the Laplacian matrix for the graph $\mathcal{G}(N,\tilde{E} / \tilde{B})$, where $\tilde{B}$ is the set of inter-area communication lines. Then, algorithm \eqref{eq:controller without delays tie-line primal-dual} can be written in the compact form
\begin{subequations}\label{eq:controller without delays tie-line compact form primal-dual}
\begin{align}
	\dot{\zeta} = & -\tilde{L} p^c + \tilde{L} \pi, \label{eq:controller without delays tie-line compact form primal-dual 1}\\
		\dot{p}^c = & - \left( p^M - p^L \right) + \tilde{L} \zeta, \label{eq:controller without delays tie-line compact form primal-dual 2}\\
		\dot{\pi} = & - \tilde{L} \zeta  + \tilde{L}_{{K}} \phi - J \hat{P}, \label{eq:controller without delays tie-line compact form primal-dual 3}\\
		\dot{\phi} = & - \tilde{L}_{{K}} \pi,\label{eq:controller without delays tie-line compact form primal-dual 4}
\end{align}
\end{subequations}
where $J = [J_{i,j}]\in \mathbb{R}^{|N| \times |{K}|}$ with $J_{j_k, k} = 1$ and $J_{j_k, k} = 0$ otherwise, and $\hat{P} = \left( \hat{P}_1,\ldots, \hat{P}_{|{K}|} \right)^T$.

%Let us
\il{We first} show the optimality of the control algorithm.
\begin{lemma}[Optimality]\label{lem:optimality tie-line primal-dual}
Let Assumption \ref{assumption convex cost function} hold.
Any equilibrium of system \eqref{eq:power system model}, \eqref{eq:generation dynamics}, \eqref{eq:generation input primal-dual}, \eqref{eq:controller without delays tie-line primal-dual} is an optimal solution to the \textup{OGR-2} problem \eqref{optimal generation regulation problem 2} with $\omega^* = \mathbf{0}_{|N|}$.
\end{lemma}
The proof is provided in Appendix\ref{appendix proof of lemma optimality}.

When there are communication delays in the channel $i \rightarrow j$, $i \in \tilde{N}_j$, bus $j$ receives delayed information of $p_i^c$, $\zeta_i$, $\pi_i$ and $\phi_i$ from bus~$i$. Then, controller \eqref{eq:controller without delays tie-line primal-dual} becomes
\begin{subequations}\label{eq:controller under delays tie-line primal-dual}
\begin{align}
\dot{\zeta}_j \hspace{-1mm} =  & \hspace{-0.5mm} \sum_{i \in \tilde{N}_j}  \hspace{-0.5mm} \alpha_{ij}  \hspace{-0.5mm} \left( p_i^c( t \hspace{-0.5mm} \hspace{-0.5mm} - \hspace{-0.5mm} T_{ij}) \hspace{-0.5mm} - \hspace{-0.5mm} p_j^c \right)  \hspace{-1mm} - \hspace{-1.0mm} \sum_{i \in \tilde{N}_j}  \hspace{-0.5mm} \mm{\alpha_{ij}}  \hspace{-0.5mm} \left( \pi_i(t \hspace{-0.5mm} - \hspace{-0.5mm}T_{ij}) \hspace{-1mm} - \hspace{-0.5mm} \pi_j \right),\\[-1mm]
		\dot{p}_j^c \hspace{-1mm}  = & - \left( p_j^M - p_j^L \right) - \sum_{i \in \tilde{N}_j} \alpha_{ij} \left( \zeta_i (t - T_{ij})- \zeta_j \right), \\
		\dot{\pi}_j  \hspace{-1mm} = & \hspace{-1.5mm} \sum_{i \in \tilde{N}_j } \hspace{-1mm}  \alpha_{ij} \left( \zeta_i (t  \hspace{-1mm} - \hspace{-1mm} T_{ij})- \zeta_j \right) \hspace{-1mm} - \hspace{-4 mm}  \sum_{i \in \tilde{N}_j \cap {C}_{k}} \hspace{-3mm} \mm{\alpha_{ij}}  \left( \phi_i(t \hspace{-1mm}  -  \hspace{-1mm}  T_{ij})  \hspace{-1mm}  -  \hspace{-1mm}  \phi_j \right)  \hspace{-1mm}  -  \hspace{-1mm}  J_j \hat{P}_{k},\\[-1mm]
		\dot{\phi}_j \hspace{-1mm} = & \sum_{i \in \tilde{N}_j \cap {C}_k} \mm{\alpha_{ij}}  \left( \pi_i(t - T_{ij}) - \pi_j \right),
\end{align}
\end{subequations}
for all $j \in N$, if variables are directly communicated. The delays may \mm{destabilize} the system if they are large.
%\todoiny{The sentence above does not read well grammatically. What `may destabilize' the system? Is it the delays, or the direct communication of the variables?}

To deal with communication delays, \il{we %let us
modify} the controllers \eqref{eq:controller without delays tie-line primal-dual} as
\begin{align}\label{eq:control with scattering transformation tie-line primal-dual}
\begin{array}{rl}
\dot{\rho}_j^{\ml{\zeta}} \hspace{-1mm} = & \hspace{-1mm} - \rho_j^{\ml{\zeta}} + \hspace{-1mm} \displaystyle \sum_{i \in \tilde{N}_j} \hspace{-1mm} \alpha_{ij} \hspace{-0.5mm} \left( r_{ij}^{p} - p_j^c \right) - \hspace{-1mm} \displaystyle \sum_{i \in \tilde{N}_j} \hspace{-1mm} \mm{\alpha_{ij}} \hspace{-0.5mm} \left( r_{ij}^{\pi'} \hspace{-0.5mm} - \hspace{-0.5mm} \pi_j \right), \\
\dot{\zeta}_j \hspace{-1mm} = & \hspace{-1mm} - \rho_j^{\ml{\zeta}} + 2 \hspace{-1mm} \displaystyle \sum_{i \in \tilde{N}_j} \alpha_{ij} \hspace{-0.5mm} \left( r_{ij}^{p} \hspace{-0.5mm} - \hspace{-0.5mm} p_j^c \right) \hspace{-0.5mm} - 2 \hspace{-1mm} \displaystyle \sum_{i \in \tilde{N}_j} \hspace{-0.5mm} \mm{\alpha_{ij}} \hspace{-0.5mm} \left( r_{ij}^{\pi'} \hspace{-0.5mm} - \hspace{-0.5mm} \pi_j \right),\\
\dot{\rho}_j^{p} \hspace{-1mm} = & \ml{- \rho_j^{p} - \hspace{-1mm} \left( p_j^M - p_j^L \right) - \displaystyle \sum_{i \in \tilde{N}_j} \alpha_{ij} \left( r_{ij}^{\zeta}- \zeta_j \right),}\\
		\dot{p}_j^c \hspace{-1mm} = & \ml{- \rho_j^{p} - \hspace{-0.5mm} 2 \left( p_j^M - p_j^L \right) - 2 \displaystyle\sum_{i \in \tilde{N}_j} \alpha_{ij} \left( r_{ij}^{\zeta}  - \zeta_j \right),}\\
\dot{\rho}_j^{\pi} \hspace{-1mm} = & \hspace{-1mm} - \hspace{-0.5mm} \rho_j^{\pi} \hspace{-0.5mm} + \hspace{-0.5mm}\hspace{-1mm} \displaystyle\sum_{i \in \tilde{N}_j} \hspace{-0.5mm} \mm{\alpha_{ij}} \hspace{-0.5mm} \left( \hspace{-0.5mm} r_{ij}^{\zeta'} \hspace{-0.5mm} - \hspace{-0.5mm} \zeta_j \hspace{-0.5mm} \right) \hspace{-1mm} - \hspace{-2mm} \displaystyle \sum_{i \in \tilde{N}_j \cap {C}_{k}} \hspace{-3.5mm} \mm{\alpha_{ij}} \hspace{-0.5mm} \left( \hspace{-0.5mm} r_{ij}^{\phi''} \hspace{-0.5mm} - \hspace{-0.5mm} \phi_j \hspace{-0.5mm} \right) \hspace{-1mm} - \hspace{-1mm} J_j \hat{P}_{k},\\
\dot{\pi}_j \hspace{-1mm} = & \hspace{-1mm} - \hspace{-0.5mm} \rho_j^{\pi} \hspace{-0.5mm} + 2 \hspace{-1mm} \displaystyle \sum_{i \in \tilde{N}_j} \hspace{-0.5mm} \mm{\alpha_{ij}} \hspace{-0.5mm} \left( \hspace{-0.5mm}  r_{ij}^{\zeta'} \hspace{-0.5mm} - \hspace{-0.5mm} \zeta_j \hspace{-0.5mm} \right) \hspace{-1mm} - \hspace{-0.5mm} 2 \hspace{-3.5mm} \displaystyle\sum_{i \in \tilde{N}_j \cap {C}_{k}} \hspace{-3.5mm} \mm{\alpha_{ij}} \hspace{-0.5mm} \left( \hspace{-0.5mm} r_{ij}^{\phi''} \hspace{-0.5mm} - \hspace{-0.5mm} \phi_j \hspace{-1mm} \right) \hspace{-1mm} - \hspace{-1mm} 2 J_j \hat{P}_{k},\\
\dot{\rho}_j^{\phi} \hspace{-1mm} = & \hspace{-1mm} -\rho_j^{\phi} + \hspace{-3.5mm} \displaystyle\sum_{i \in \tilde{N}_j \cap {C}_k} \hspace{-3.5mm} \mm{\alpha_{ij}} \left( r_{ij}^{\pi''} - \pi_j \right),\\
		\dot{\phi}_j \hspace{-1mm} = & \hspace{-1mm} - \rho_j^{\phi} + 2 \hspace{-3.5mm} \displaystyle \sum_{i \in \tilde{N}_j \cap {C}_k} \hspace{-3.5mm} \mm{\alpha_{ij}} \left(r_{ij}^{\pi''} - \pi_j \right),\\[-2mm]
\end{array}
\end{align}
for all $j \in N$ and $j \in C_{k}$, where $\rho_j^{\zeta}$, $\rho_j^{p}$, $\rho_j^{\pi}$, $\rho_j^{\phi}$ are auxiliary states resulting from parallel feedforward compensation on \eqref{eq:controller without delays tie-line primal-dual}, \il{similarly to \eqref{eq:control with scattering transformation primal-dual}. As} \eqref{eq:controller without delays tie-line primal-dual} has three pairs of communication variables, there are three sets of scattering \il{transformations,} and
$r_{ij}^{p}$, $r_{ij}^{\zeta}$, $r_{ij}^{\zeta'}$, $r_{ij}^{\pi'}$, $r_{ij}^{\pi''}$, $r_{ij}^{\phi''}$ represent variables bus $j$ \icl{formulates using the information
%extracted from scattering variables in the
communicated \ic{from}} bus $i$.
The above extended algorithm also inherits passivity properties. \ic{As in Lemma \ref{lem:passivity of reformulated controller} the superscript $^*$ in Lemma \ref{lem:passivity of controller tie-line primal-dual} below is used to denote values at an equilibrium point.}
\begin{lemma}\label{lem:passivity of controller tie-line primal-dual}
	The system described by \eqref{eq:power system model}, \eqref{eq:generation dynamics}, \eqref{eq:generation input primal-dual}, \eqref{eq:control with scattering transformation tie-line primal-dual} is locally passive from input
	$\begin{bmatrix}
		\tilde{r}^{p} & \tilde{r}^{\zeta} & \tilde{r}^{\zeta'} & \tilde{r}^{\pi'} & \tilde{r}^{\pi''} & \tilde{r}^{\phi''}
	\end{bmatrix}^T$ to output
	$\begin{bmatrix}
		\tilde{\zeta} & - \tilde{p}^c & \tilde{\pi} & - \tilde{\zeta} & \tilde{\phi} & - \tilde{\pi}
	\end{bmatrix}^T$, where the $j^{th}$ element of the variables are defined by
	\begin{align*}
	& \left[ \begin{smallmatrix}
		\tilde{r}_j^p & \tilde{r}_j^{\zeta} & \tilde{r}_j^{\zeta'} & \tilde{r}_j^{\pi'}
	\end{smallmatrix} \right]^T
	\hspace{-2mm} = \hspace{-1.5mm} \sum_{i \in \tilde{N}_j} \hspace{-1mm} \alpha_{ij}
	\hspace{-0.5mm}	\left[ \begin{smallmatrix}
		r_{ij}^{p} - r_{ij}^{p,*} & r_{ij}^{\zeta} -r_{ij}^{\zeta,*}  & r_{ij}^{\zeta'} - r_{ij}^{\zeta,*} & r_{ij}^{\pi'} - r_{ij}^{\pi,*}
	\end{smallmatrix} \right]^T \hspace{-1mm},\\
& \left[
	\begin{smallmatrix}
		\tilde{r}_j^{\pi''} & \tilde{r}_j^{\phi''}
	\end{smallmatrix} \right]^T = \hspace{-2mm} \sum_{i \in \tilde{N}_j \cap {C}_{k}} \alpha_{ij}
	\left[ \begin{smallmatrix}
		r_{ij}^{\pi''} - r_{ij}^{\pi,*} & r_{ij}^{\phi''} - r_{ij}^{\phi,*}
	\end{smallmatrix} \right]^T \hspace{-1mm},  ~\lmm{j \in C_{k}},\\
& \left[
\begin{smallmatrix}
		\tilde{\zeta}_j & \tilde{p}_j^c & \tilde{\pi}_j & \tilde{\phi}_j
	\end{smallmatrix} \right]^T =
	\left[
	\begin{smallmatrix}
		\zeta_j - \zeta_j^* & p_j^c - p_j^{c,*}  & \pi_j - \pi_j^* & \pi_j - \phi_j^*
	\end{smallmatrix} \right]^T
\end{align*}
%	\todoiny{what is $K$ above in the LHS? It does not seem to appear in the RHS, Similarly in the previous expression. \lmm{I intended to say that $\tilde{r}_j^{\pi, {K}}$, $\tilde{r}_j^{\phi, {K}}$ are variables of information received from neighboring buses only in the same area. }}
\end{lemma}
The proof is given in Appendix\ref{appendix proof of lem passivity of controller tie-line primal-dual}.
\begin{remark}
	The input and output in \Cref{lem:passivity of controller tie-line primal-dual} contain repeated variables.  This may \ic{appear %seem
 redundant,} %at the first glance
but they are formulated in this way such that one can construct a passive scattering transformation for the communication channel.
\end{remark}

We \ic{now} %then
define a new set of scattering variables, \icl{which are the \ic{variables} that get communicated between buses.} %$i, j$.}
%to be communicated and extract the information of ${r}^{p}_{ij}$, ${r}^{\zeta}_{ij}$, %${r}^{\pi}_{ij}$, ${r}^{\phi}_{ij}$ for bus $j$ from them.
Let $r_{ji} = \left({r}^{p}_{ji}, {r}^{\zeta}_{ji}, {r}^{\zeta'}_{ji}, {r}^{\pi'}_{ji},  \sigma_{ji} {r}^{\pi''}_{ji},  \sigma_{ji} {r}^{\phi''}_{ji} \right)^T$  and $y_{ji} = \left( \zeta_i, - p_i^c, \pi_i, - \zeta_i, \sigma_{ji}\phi_i, - \sigma_{ji} \pi_i \right)^T$, where $\sigma_{ji} = 0$ if $(i,j) \in \tilde{B}$ or $(j,i) \in \tilde{B}$, and $\sigma_{ji} = 1$ otherwise. %With an abuse of notation,
\icl{The new} scattering variables are \icl{given by\footnote{\icl{For convenience in the presentation, we use the notation $s_{\overrightarrow{ij}}, s_{\overleftarrow{ij}}, s_{\overleftarrow{ji}}, s_{\overrightarrow{ji}}, E_s$ as in the previous section, even though these are defined below in a different way as more variables get communicated.}}
%\todoiny{\st{why did you say abuse of notation?}}
}
\begin{equation}\label{eq:scattering transformation tie-line primal-dual}
\small
\begin{aligned}
%& s_{\overrightarrow{ij}} = \frac{- 1}{\sqrt{2}}\left( \begin{bmatrix}
%		{r}^{p}_{ji} \\ {r}^{\zeta}_{ji} \\ {r}^{\zeta}_{ji} \\ {r}^{\pi}_{ji} \\ {r}^{\pi}_{ji} \\ {r}^{\phi}_{ji}
%	\end{bmatrix}  - \begin{bmatrix}
%		\zeta_i \\ - p_i^c \\ \pi_i \\ - \zeta_i \\ \phi_i \\ - \pi_i
%	\end{bmatrix} \right),
& s_{\overrightarrow{ij}} = -\frac{1}{\sqrt{2}}\left( r_{ji} - y_{ji} \right), ~
s_{\overleftarrow{ij}} = - \frac{1}{\sqrt{2}} \left( r_{ji}  + y_{ji}\right)\\
& s_{\overleftarrow{ji}} = \frac{1}{\sqrt{2}} \left( r_{ij} + y_{ij} \right),~
s_{\overrightarrow{ji}} = \frac{1}{\sqrt{2}} \left(r_{ij} - y_{ij} \right).
	\end{aligned}
\end{equation}
Since $(s_{\overrightarrow{ij}},~s_{\overleftarrow{ji}})$ and $(s_{\overrightarrow{ji}},~s_{\overleftarrow{ij}})$ are the variables in the communication channels $i \rightarrow j$ and $j \rightarrow i$, respectively, it \icl{holds} that
\begin{equation}\label{eq:scattering variables under delays tie-line primal-dual}
	s_{\overleftarrow{ij}}(t) = E_s s_{\overrightarrow{ji}}(t - T_{ji}), \quad s_{\overleftarrow{ji}}(t) = E_s s_{\overrightarrow{ij}}(t - T_{ij})
\end{equation}
\icl{where $T_{ji}$, $T_{ij}\geq0$ are delays} in the communication channels $j \rightarrow i$ and $i \rightarrow j$, respectively. \ic{We also define the matrix gain} $E_s = I_3 \otimes \begin{bmatrix}
	0 & -1 \\ 1 & 0
\end{bmatrix}$. % with an abuse of notation.
The new scattering \ic{transformation %basically
contains} three groups of input and output variables, \icl{and is analogous to the one used in the previous section.} %similar to the previous ones.
\ic{Therefore,} it is also passive, as \icl{stated} below.
\begin{lemma}\label{lem:passivity of scattering transformation tie-line}
	The system \ic{described by} %consisting of
\eqref{eq:scattering transformation tie-line primal-dual} and \eqref{eq:scattering variables under delays tie-line primal-dual} is passive from input $
\begin{bmatrix}
\tilde{y}_{ij} \\ \tilde{y}_{ji}
\end{bmatrix}$
to output
$- \begin{bmatrix}
\tilde{r}_{ij}\\ \tilde{r}_{ji}
 	\end{bmatrix}$
where $\tilde{r}_{ij} = r_{ij} - r_{ij}^*$ \il{and $\tilde{y}_{ij} = y_{ij} - y_{ij}^*$.}
%are defined in \Cref{lem:passivity of controller tie-line primal-dual} \il{and the text above}.
\end{lemma}
The proof is given in Appendix\ref{appendix proof of lem passivity of scattering transformation tie-line}.

%\todoiny{I think move the proof in the appendix}

\ic{Using} the passivity of the two components \eqref{eq:power system model}, \eqref{eq:generation dynamics}, \eqref{eq:generation input primal-dual}, \eqref{eq:control with scattering transformation tie-line primal-dual}, and  \eqref{eq:scattering transformation tie-line primal-dual}, \eqref{eq:scattering variables under delays tie-line primal-dual} \ic{stated in} \Cref{lem:passivity of controller tie-line primal-dual} and \Cref{lem:passivity of scattering transformation tie-line}, \icl{respectively,
we are ready to deduce} asymptotic stability of the closed-loop system under \icl{arbitrary constant} delays.
\begin{theorem}\label{thm:Convergence under delays tie-line primal-dual}
Let Assumption~\ref{assumption convex cost function} hold and consider an equilibrium \il{point} of system \eqref{eq:power system model}, \eqref{eq:generation dynamics}, \eqref{eq:generation input primal-dual}, \eqref{eq:control with scattering transformation tie-line primal-dual}, \eqref{eq:scattering transformation tie-line primal-dual}, \eqref{eq:scattering variables under delays tie-line primal-dual} in which Assumption \ref{assumption angle} holds.
\icl{The delays $T_{ij}\geq0$, $\forall (i, j) \in \tilde{E}$ are assumed to be constant, and are allowed to take arbitrary bounded values and be heterogeneous.}
%unknown and heterogeneous.
Then, there exists an open neighborhood $\Omega$ about this equilibrium \il{point} such that the solutions of \eqref{eq:power system model}, \eqref{eq:generation dynamics}, \eqref{eq:generation input primal-dual}, \eqref{eq:control with scattering transformation tie-line primal-dual}, \eqref{eq:scattering transformation tie-line primal-dual}, \eqref{eq:scattering variables under delays tie-line primal-dual} with initial conditions in \icl{$\Omega$
%$\mathcal{C} \left( \left[ -\max_{\forall i,j}\{T_{ij}\}, 0 \right], \Omega \right)$
converge} \il{to an equilibrium point that solves} %to the optimal solution of
the
\textup{OGR-2} problem \eqref{optimal generation regulation problem 2} with $\omega^* = \mathbf{0}_{|N|}$.
\end{theorem}
The proof is given in Appendix\ref{appendix proof of theorem convergence under delays tie-line primal-dual}.

\subsection{Generation Boundedness Constraints}\label{subsec:Generation Boundedness Constraints}
In this subsection, we consider bounds for the minimum and maximum values of generation to account for a more realistic operating condition.
\begin{align}\label{optimal generation regulation problem 3}
\begin{array}{rl}
	\textbf{OGR-3:} \quad & \underset{p^M}{\min} \displaystyle \sum_{j \in G} Q_j (p_j^M),\\
	& \text{subject to } \displaystyle \sum_{j \in G} p_j^M = \sum_{j \in N} p_j^L,\\
	& p_j^{M, \min} \leq p_j^M \leq p_j^{M,\max}, \forall j \in G,
\end{array}
\end{align}
where $p_j^{M,\min}$, $p_j^{M,\max}$ are lower and upper bounds for the for generation at bus $j$.

%\todoing{Looking at it again, and the assumption in the Hale and Lunel book, saturation %functions can be handled. The only condition needed is continuity of the vector field and }
In the undelayed case, a \icl{function with saturation} can be used to modify the generation input \il{so as to satisfy thee generation bounds at equilibrium} \cite{kasis2019stability}.
\il{In order, however, to avoid complications associated with a non-smooth vector field in an infinite dimensional setting associated with delays, we use below an alternative scheme that leads to smooth dynamics.}
%, for satisfying the generation constraints at equilibrium .}
%However, \il{the non-smoothness of the dynamics in conjunction with the infinite dimensional character of the system complicated the analysis and a different }
%this does not apply to the delayed case since the non-smoothness of the dynamics \mm{complicates the direct application of an} invariance principle \mml{\ic{to} the delay differential algebraic equation}.
%%\todoiny{Why was this the case?}
%Therefore,  \ic{a different approach is used below.}

\il{In particular, to} solve the \textup{OGR-3} problem \eqref{optimal generation regulation problem 3}, the controller \mm{\eqref{eq:control with scattering transformation primal-dual}} is \ic{left} unchanged, while the generation input \eqref{eq:generation input primal-dual} is changed to
\begin{equation}\label{eq:generation input generation bounds}
	u_j \hspace{-1mm} =  k_{c,j} \left( p_j^c - \omega_j \right) +  \frac{p_j^M}{k_{g,j}} - k_{c,j} \left( Q_j'\left(p_j^M\right) - \lambda_{j}^{2}  + \mu_{j}^{2}  \right)  \hspace{-0.5mm}, ~ j \in G
\end{equation}
where $\lambda_j^2$, $\mu_j^2$ can be viewed as non-negative Lagrange multipliers for the inequalities.
%\todoiny{ is "non-negative" referring to the nonlinearities or the Lagrange multipliers?}
The variables $\lambda_j$, $\mu_j$ follow the \ic{dynamics}
\begin{subequations}\label{eq:local inequality multipliers}
\begin{align}
\dot{\lambda}_j = 2 \lambda_j \left(p_j^{M,\min} - p_j^M \right), ~ \lambda_j(0) > 0,\\
\dot{\mu}_j = 2 \mu_j \left( p_j^M- p_j^{M, \max} \right), ~ \mu_j (0) > 0.
\end{align}
\end{subequations}
The update rule \eqref{eq:local inequality multipliers} is \ic{used here as a means of %proposed to
solving an inequality constrained optimization problem} with smooth dynamics \cite{li2018generalized}.
We show \ic{below that optimality} and convergence are guaranteed by the augmented \ic{dynamics described above.} %s in the following \ic{Lemma}.
\begin{lemma}[Optimality]\label{lem:optimality generation bounds}
Let Assumption \ref{assumption convex cost function} hold.
Any equilibrium of system \eqref{eq:power system model}, \eqref{eq:generation dynamics}, \eqref{eq:controller without delays primal-dual}, \eqref{eq:generation input generation bounds}, \eqref{eq:local inequality multipliers} is an optimal solution to the \textup{OGR-3} problem \eqref{optimal generation regulation problem 3} with $\omega^* = \mathbf{0}_{|N|}$ and $p^{c,*} = \text{Im}(\mathbf{1}_{|N|})$.
\end{lemma}
The proof is given in Appendix\ref{appendix proof of lemma optimality generation bounds}.

The stability of the primal-dual controlled \ic{system considered that incorporates generation bounds and delays is stated in the Theorem below.}  %as follows.in
%\todoiny{I guess an issue with constraints as formulated is that they are satisfied only at equilibrium, whereas the non-smooth system would also satisfy them transiently?}
\begin{theorem}\label{thm:convergence under delays generation bounds}
Let Assumption~\ref{assumption convex cost function} hold and
consider an equilibrium \il{point} of system \eqref{eq:power system model}, \eqref{eq:generation dynamics}, \eqref{eq:control with scattering transformation primal-dual}, \eqref{eq:scattering transformation primal-dual}, \eqref{eq:scattering variables under delays primal-dual}, \eqref{eq:generation input generation bounds}, \eqref{eq:local inequality multipliers}, in which Assumption \ref{assumption angle} holds.
%Suppose the constant delays $T_{ij}$, $\forall (i, j) \in \tilde{E}$ are bounded, unknown and %heterogeneous.
\ic{The delays $T_{ij}\geq0$, $\forall (i, j) \in \tilde{E}$ are assumed to be constant, and are allowed to take arbitrary bounded values and be heterogeneous.}
\il{Then, there} exists an open neighborhood $\Omega$ about this equilibrium \il{point} such that the solutions of  \eqref{eq:power system model}, \eqref{eq:generation dynamics}, \eqref{eq:control with scattering transformation primal-dual}, \eqref{eq:scattering transformation primal-dual}, \eqref{eq:scattering variables under delays primal-dual}, \eqref{eq:generation input generation bounds}, \eqref{eq:local inequality multipliers} with initial conditions in
\mm{$\Omega$}
%$\mathcal{C} \left( \left[ -\max_{\forall i,j}\{T_{ij}\}, 0 \right], \Omega \right)$
%asymptotically
\il{converge to an equilibrium point that solves the} \textup{OGR-3} problem \eqref{optimal generation regulation problem 3} with $\omega^* = \mathbf{0}_{|N|}$.
\end{theorem}
The proof is given in Appendix\ref{appendix proof of theorem convergence under delays generation bounds}.

\begin{remark}
\il{It should be noted that the} \ic{dynamics in \eqref{eq:local inequality multipliers} follow} from the dual gradient ascent of the generalized Lagrangian
\begin{align}\label{eq:generalized Lagrangian}
\mathcal{L} = \hspace{-1mm} \sum_{j \in G} Q_{j} \left( p_j^M \right) + \lambda_j^2 \left( p_j^{M,\min} - p_j^M \right) + \mu_j^2 \left( p_j^{M} - p_j^{M,\max} \right)
\end{align}
%\mm{We resort to this method as the smoothness of \eqref{eq:local inequality multipliers}  %ensures the  invariance principle is applicable \mm{to the delay differential algebraic %equation}.}
%\todoing{Removed the sentence to avoid repetition.}
%\todoiny{As before let's discuss again what the issue was as we have now changed the invariance principle arguments.}
%following from the fact that
%the delayed system is infinite dimensional, and establishing compactness of solutions via the level sets of the Lyapunov functional is not feasible for our controllers.
%Our proof instead
%\ml{focuses on the convergence to the positive limit set,}
%requires smooth trajectories.
\end{remark}

\subsection{Observer-based Estimation for \il{Demand}}\label{subsec:Observer-based Estimation for Demands}
The control strategy proposed in \eqref{eq:control with scattering transformation primal-dual}, \eqref{eq:scattering transformation primal-dual}, \eqref{eq:scattering variables under delays primal-dual} requires the explicit knowledge of the uncontrollable frequency independent demand $p^L$.
In this subsection, we include observer dynamics borrowed from \cite{kasis2019stability} to relax this requirement without affecting the stability and optimality presented in Theorem \ref{thm:convergence under delays primal-dual}. The controller under observer dynamics is
\begin{subequations}\label{eq:controller observer-based primal-dual}
\begin{align}
\hspace{-2mm} & \dot{\rho}_j^{\ml{\zeta}} = - \rho_j^{\ml{\zeta}} +  \sum_{i \in \tilde{N}_j} \alpha_{ij}\left( r_{ij}^{p} - p_j^c \right), ~ j \in N\\
	&\dot{\zeta}_j = - \ml{\rho_j^{\zeta}} + 2 \sum_{i \in \tilde{N}_{j}} \alpha_{ij}\left( r_{ij}^p - p_j^c \right), ~ j \in N\\
	&\ml{\dot{\rho}_j^{p} = - {\rho}_j^{p} - \left( p_j^M - \chi_j \right)\hspace{-0.5mm}  - \hspace{-1mm}  \sum_{i \in \tilde{N}_{j}} \alpha_{ij} \left( r_{ij}^{\zeta} - \zeta_j \right) \hspace{-0.5mm} , ~ j \in N }\\
	&\dot{p}_j^c = - \ml{{\rho}_j^{p} \hspace{-1mm} - \hspace{-0.5mm} 2 \hspace{-0.5mm} \left( p_j^M \hspace{-0.5mm} - \hspace{-0.5mm} \chi_j \right) \hspace{-0.5mm} - \hspace{-1mm}  2 \hspace{-0.5mm} \sum_{i \in \tilde{N}_{j}} \alpha_{ij} \left( r_{ij}^{\zeta} \hspace{-0.5mm}- \hspace{-0.5mm} \zeta_j \right) \hspace{-0.5mm}, ~ j \in N}\\
	&\tau_{\chi,j} \dot{\chi}_j = b_j - \omega_j - p_j^c - \chi_j, ~ j \in G \label{eq:controller observer-based primal-dual_3}\\
	&M_j \dot{b}_j \hspace{-0.5mm} = \hspace{-0.5mm} - \chi_j \hspace{-0.5mm}  + \hspace{-0.5mm}  p_j^M \hspace{-0.5mm}  - \hspace{-0.5mm}  \Lambda_j \omega_j \hspace{-0.5mm} - \hspace{-2mm}  \sum_{k: j \rightarrow k} p_{jk} \hspace{-0.5mm}  + \hspace{-1mm}  \sum_{i: i \rightarrow j}\hspace{-0.5mm}  p_{ij}, ~ j \in G \label{eq:controller observer-based primal-dual_4}\\
	&0 = -\chi_j - \Lambda_j \omega_j - \sum_{k: j \rightarrow k} p_{jk} + \sum_{i: i \rightarrow j} p_{ij}, ~ j \in L \label{eq:controller observer-based primal-dual_5}
\end{align}
\end{subequations}

The stability of the primal-dual controlled system relaxing demand measurement under delays is given as follows.
\begin{theorem}\label{thm:convergence observer-based primal-dual}
Let Assumption~\ref{assumption convex cost function} hold and consider an equilibrium \il{point} of system \eqref{eq:power system model}, \eqref{eq:generation dynamics}, \eqref{eq:generation input primal-dual}, \eqref{eq:controller observer-based primal-dual}, \eqref{eq:scattering transformation primal-dual}, \eqref{eq:scattering variables under delays primal-dual} in which Assumption \ref{assumption angle} holds.
\ic{The delays $T_{ij}\geq0$, $\forall (i, j) \in \tilde{E}$ are assumed to be constant, and are allowed to take arbitrary bounded values and be heterogeneous.}
%Suppose the constant delays $T_{ij}$, $\forall (i, j) \in \tilde{E}$ are bounded, unknown and heterogeneous.
Then, there exists an open neighborhood $\Omega$ about this equilibrium \il{point} such that the solutions of \eqref{eq:power system model}, \eqref{eq:generation dynamics}, \eqref{eq:generation input primal-dual}, \eqref{eq:controller observer-based primal-dual} , \eqref{eq:scattering transformation primal-dual}, \eqref{eq:scattering variables under delays primal-dual} with initial conditions in %$\mathcal{C} \left( \left[ -\max_{\forall i,j}\{T_{ij}\}, 0 \right], \Omega \right)$
\mm{$\Omega$} \il{converge to an equilibrium point that solves} %to the optimal solution of
the \textup{OGR} problem \eqref{optimal generation regulation problem} with $\omega^* = \mathbf{0}_{|N|}$.
\end{theorem}
The proof is given in Appendix\ref{appendix proof of theorem convergence observer-based primal-dual}.

Note that we do not include delays in the observer layer \eqref{eq:controller observer-based primal-dual_3}, \eqref{eq:controller observer-based primal-dual_4},  \eqref{eq:controller observer-based primal-dual_5} for simplicity in the presentation.
If delays in the estimation of $p_{ij}$ are considered,  a similar controller incorporating \ic{a scattering} transformation can be readily constructed.
It is also worth noting that the three extensions in \Cref{subsec:Tie-line Power Flow Constraints,subsec:Generation Boundedness Constraints,subsec:Observer-based Estimation for Demands} can \ic{be easily combined together, while still maintaining the optimality and convergence properties that have been derived in the paper}.  They are \ic{presented} separately here for brevity and \ic{to facilitate} readability.
%More operational constraints can be considered, for example, the robust stability results in this paper can be readily extended to systems with higher-order turbine governor models by similar arguments to  \cite{kasis2019stability}.

\section{Simulations}\label{sec:Simulations}
In this section, we illustrate our results with simulations on a 5-bus example and then on the IEEE-39 test system.

\subsection*{Example 1}
\begin{figure}[htbp]
	\centering
	\includegraphics[width = .5 \linewidth]{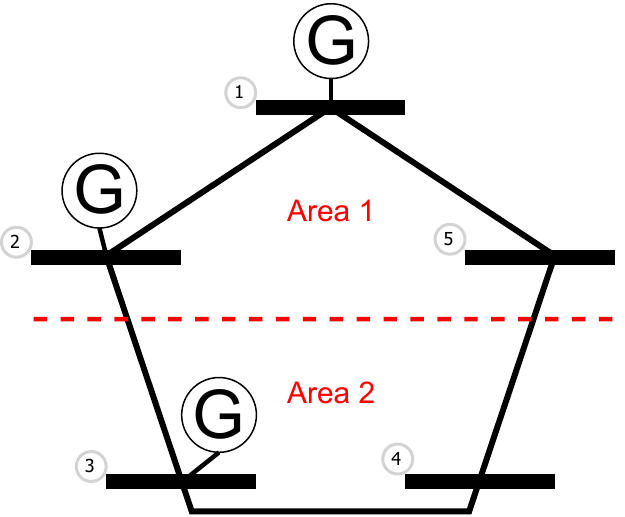}
	\caption{A $5$-bus power network example with three generators and two control areas.}
	\label{fig:5bus_example}
\end{figure}

Consider a $5$-bus power network example with three generators and two control areas, as shown in \cref{fig:5bus_example}.
The cost functions are $Q_j (p_j^M) = \frac{1}{2} q_j \left(p_j^M  - c_j \right)^2$, $j \in G$. The bus parameters are given in \Cref{table:2}. Let $k_{g,j} = 1$, $k_{c,j} = 1$, $\forall j \in N$.
\begin{table}[h!]
\centering
\begin{tabular}{c| c c c c c}
\hline
Bus & 1 & 2 & 3 & 4 & 5\\
\hline
$q_j$ & 2.4 & 4 & 3.4 & -- & --\\
\hline
$c_j$ & 0.3 & 0.1 & 0.2 & -- & --\\
\hline
$M_j$ & 13 & 12.1 & 14.3 & -- & --\\
\hline
$\tau_j$ & 0.3 & 0.4 & 0.35 & -- & --\\
\hline
$\Lambda_j$ & 1 & 0.8 & 1.1& 1 & 0.9\\
\hline
$p^L_j$ & 0.1 & 0.2 & 0.3 & 0.4 & 0.5\\
\hline
\end{tabular}
\caption{Bus parameters.}
\label{table:2}
\end{table}
The loads $p^L$ are added at time $t = 5 s$. The tie-line power flow is scheduled to be $\hat{P} = [-0.5, 0.5]$, which is only known to bus $2$ and $3$.

\begin{figure}[htbp]
	\centering
	\subfigure[Frequencies under controller \eqref{eq:redundant update under delays},~\eqref{eq:controller without delays bus dynamics primal-dual}.]{\includegraphics[width = 0.49\linewidth]{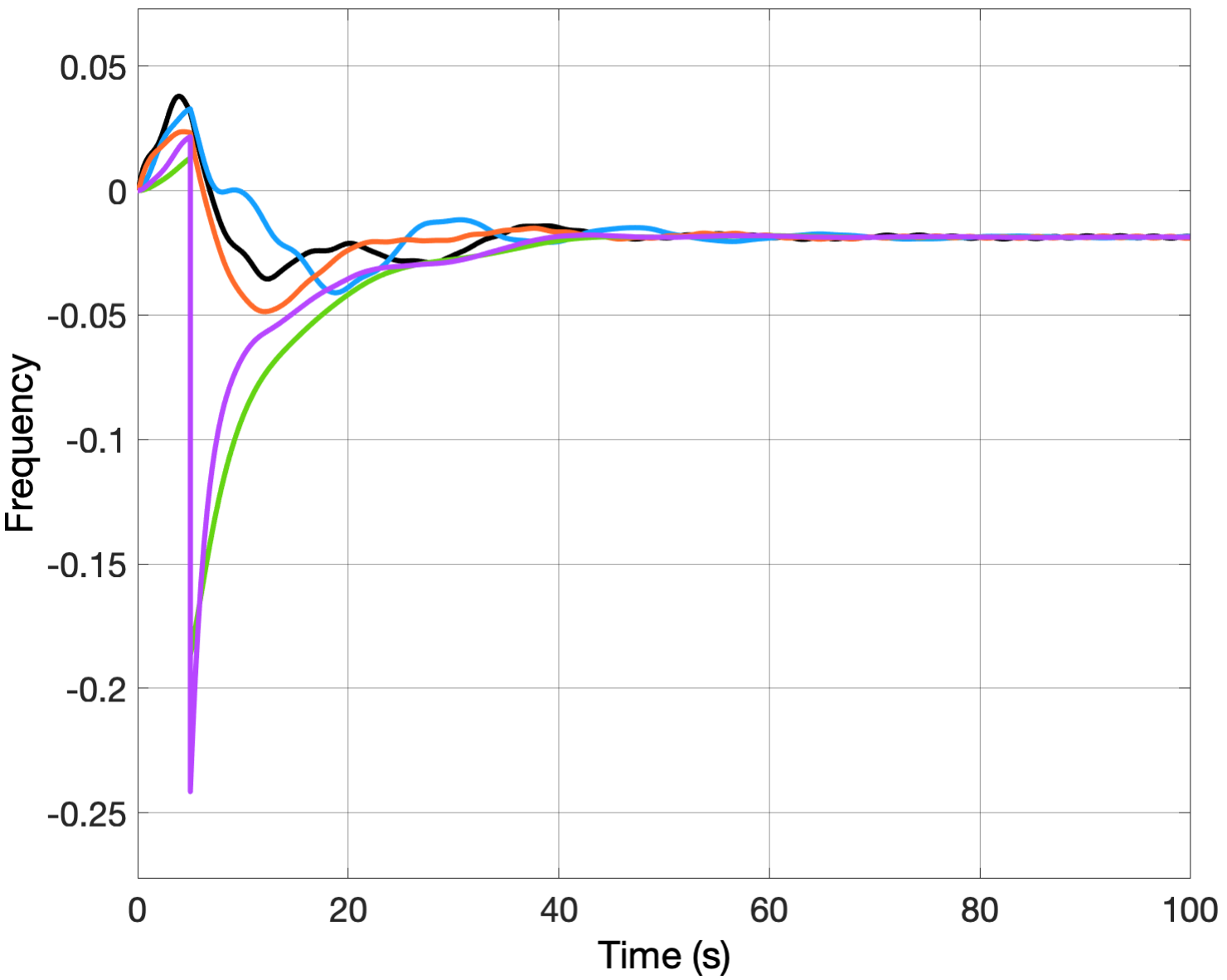}\label{fig:original controllers failed under delays_1}}
	\subfigure[Frequencies under controller \eqref{eq:controller under delays reformulated-1 primal-dual}.]{\includegraphics[width = 0.49\linewidth]{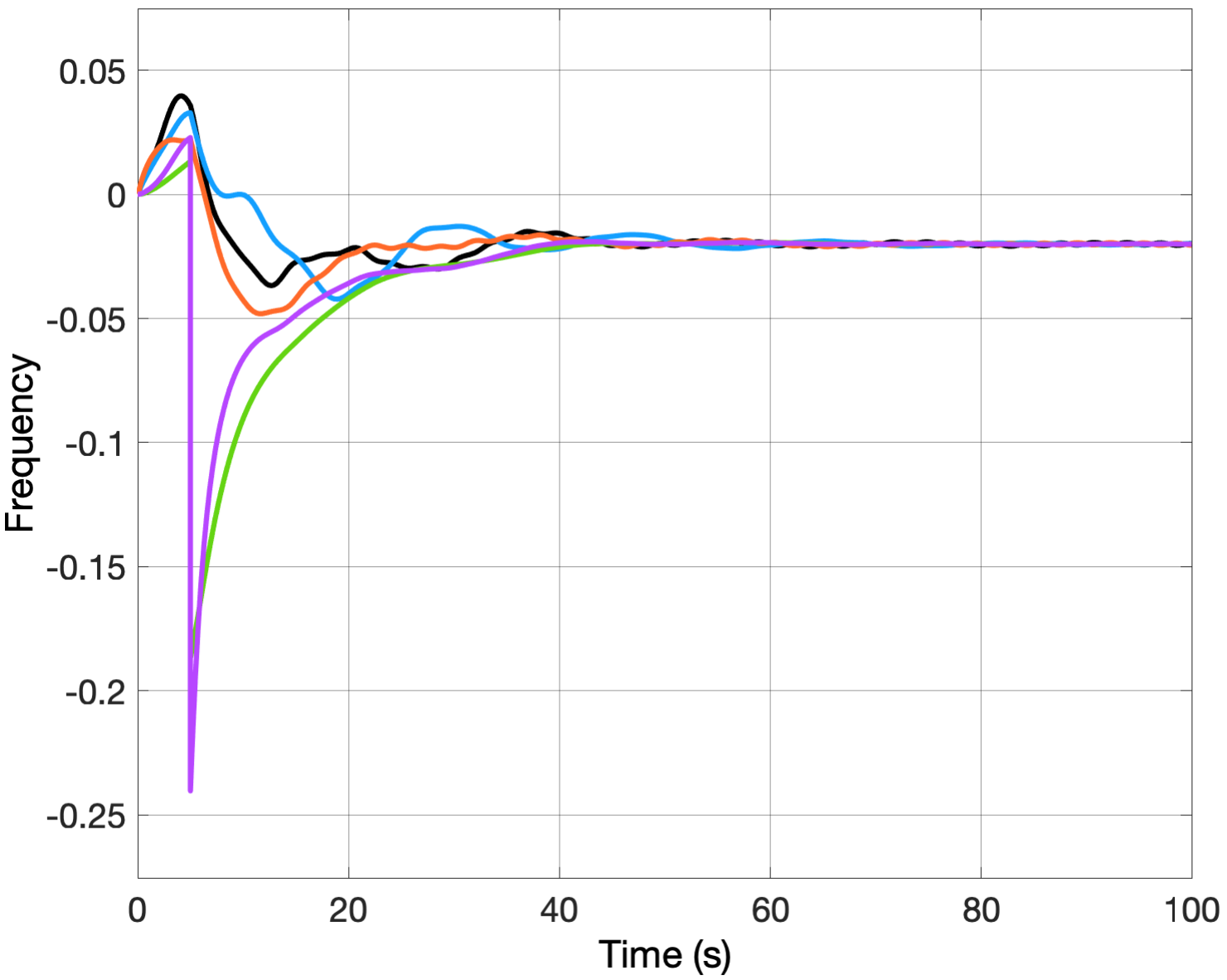}\label{fig:original controllers failed under delays_2}}
	\caption{\il{The frequency at each of the buses in Example 1 when the primal-dual control policies \eqref{eq:redundant update under delays},~\eqref{eq:controller without delays bus dynamics primal-dual} and \eqref{eq:controller under delays reformulated-1 primal-dual}, respectively, are implemented} with $T_{ij} = 0.01s$.}
	\label{fig:original controllers failed under delays}
\end{figure}

We first show that the primal-dual controllers in the original formulation \eqref{eq:controller without delays primal-dual} and the form \eqref{eq:controller without delays reformulation-1 primal-dual} are incapable of %addressing
\il{incorporating} communication delays \il{efficiently}, though they only require communication of the power command signal.
As discussed in \Cref{subsection:Equivalent Reformulation of the Primal-dual Control}, the controller \eqref{eq:controller without delays edge dynamics primal-dual} under communication delays becomes \revise{\eqref{eq:redundant update under delays}}
%\begin{subequations}\label{eq:controller under delays primal-dual}
%\begin{align}
%\gamma_{ij} \dot{\psi}_{ij}^{j} = & p_i^c (t - T_{ij}) - p_j^c, ~ (i,j) \in \tilde{E}\\
%\gamma_{kj} \dot{\psi}_{kj}^{j} = & p_j^c - p_k^c (t - T_{kj}), ~ (j,k) \in \tilde{E}\\
%		\gamma_{j} \dot{p}_j^c = & - \left( p_j^M - p_j^L \right) - \hspace{-2mm} \sum_{k: j \rightarrow k} \hspace{-1mm} \psi_{jk}^{j} \hspace{-0.5mm} + \hspace{-1.5mm} \sum_{i: i\rightarrow j} \psi_{ij}^{j}, ~ j \in N
%\end{align}
%\end{subequations}
with $T_{ij}$ representing the delay in the communication channel $i \rightarrow j$,
while the controller \eqref{eq:controller without delays reformulation-1 primal-dual} under communication delays becomes
\begin{align}\label{eq:controller under delays reformulated-1 primal-dual}
	\dot{\xi}_j \hspace{-1mm} = \hspace{-2mm} \sum_{i \in \tilde{N}_j} \alpha_{ij}\left( p_i^c (t - T_{ij}) - p_j^c \right),  ~
		\dot{p}_j^c \hspace{-1mm} = \hspace{-1mm} - \hspace{-1mm} \left( p_j^M \hspace{-0.5mm} - \hspace{-0.5mm} p_j^L \right) \hspace{-0.5mm} + \hspace{-0.5mm} \xi_j,
\end{align}
for all $j \in N$, with initial conditions $ \sum_{j \in N}\zeta_j(0) =  0$. The performance under controllers \revise{\eqref{eq:redundant update under delays}, \eqref{eq:controller without delays bus dynamics primal-dual}} and \eqref{eq:controller under delays reformulated-1 primal-dual} with $T_{ij} = 0.01s$, $\forall i,~j$, are shown in \cref{fig:original controllers failed under delays}. We can observe that small delays affect the equilibrium and the frequency is not restored to the nominal value.

\begin{figure}[htbp]
	\centering
	\subfigure[Controller \eqref{eq:controller under delays reformulation primal-dual} with delays $T_{ij} = 0.03s$.]{\includegraphics[width = 0.49\linewidth]{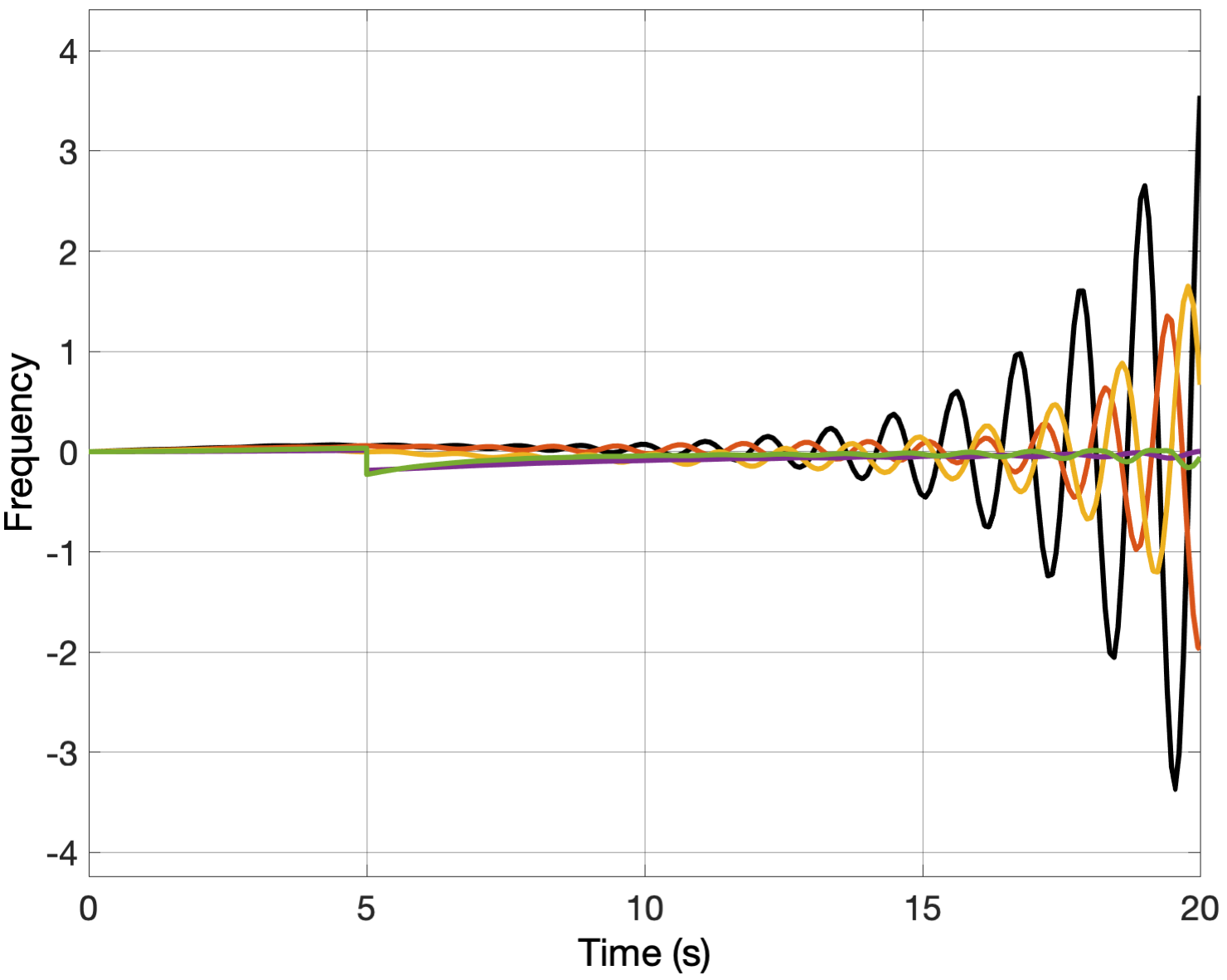}\label{fig:example 7}}
	\subfigure[Controller \eqref{eq:controller under delays tie-line primal-dual} with delays $T_{ij} = 0.03s$.]{\includegraphics[width = 0.49\linewidth]{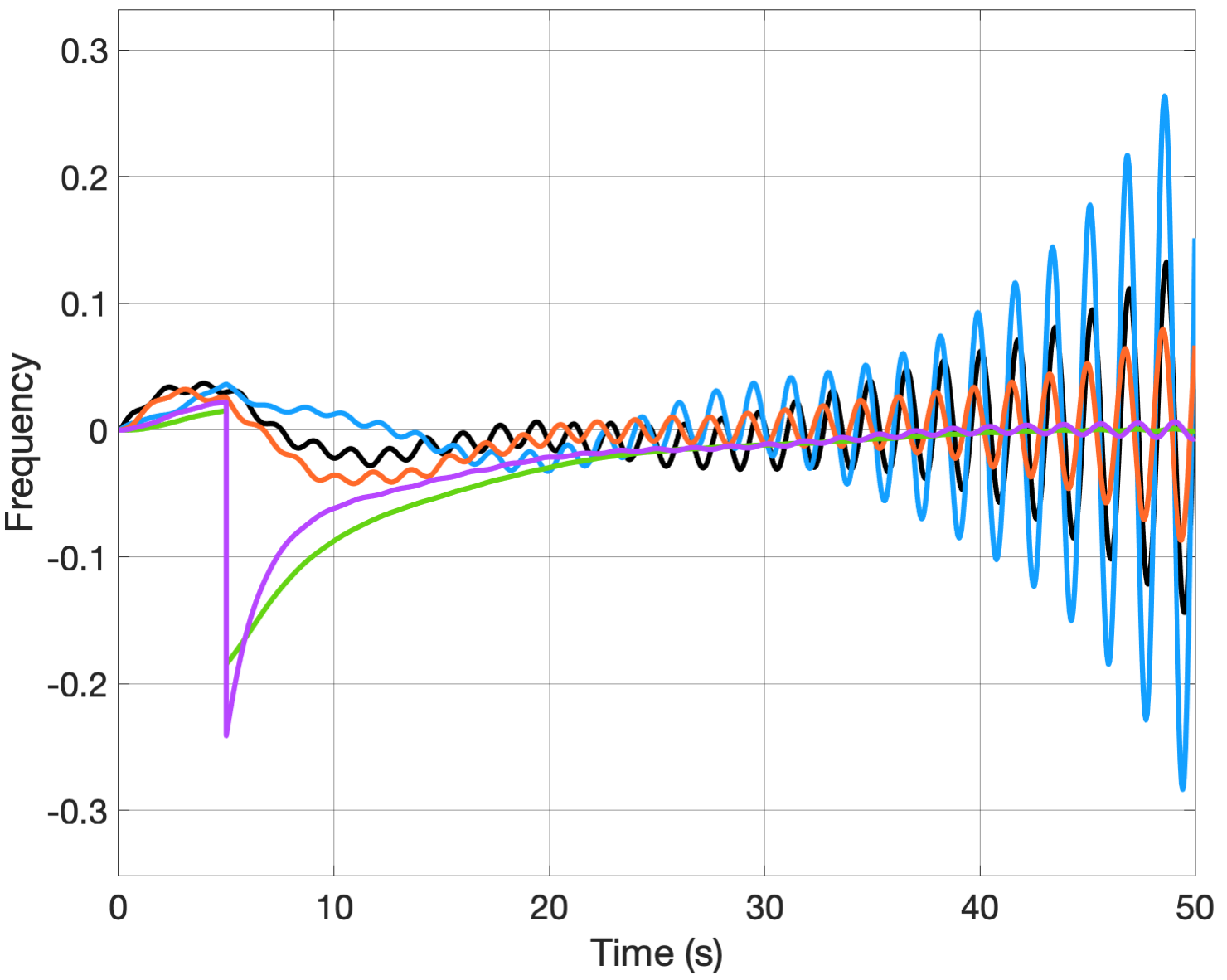}\label{fig:example 8}}
	\subfigure[Controller \eqref{eq:controller without delays reformulation primal-dual} without delays.]{\includegraphics[width = 0.49\linewidth]{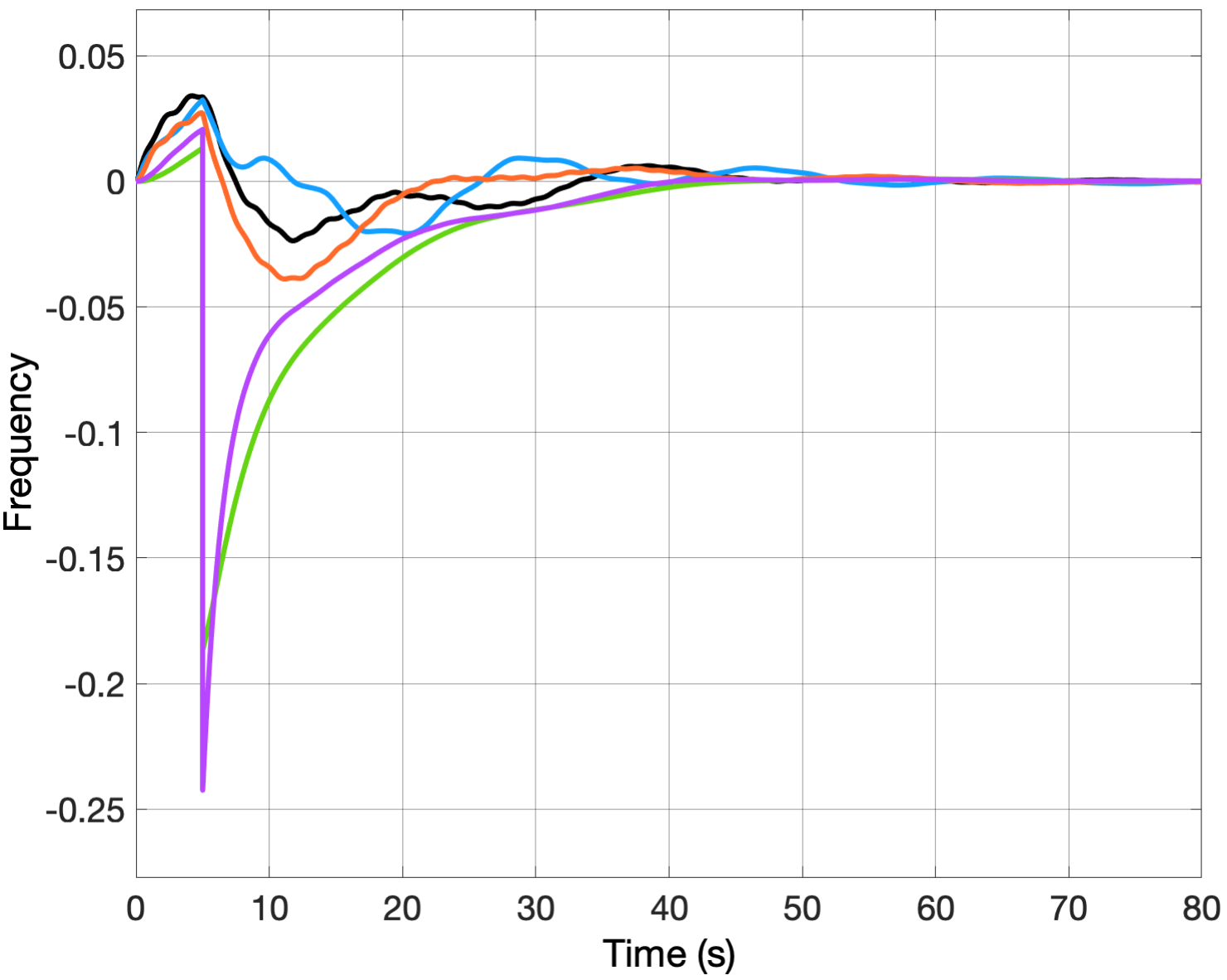}\label{fig:example 1}}
	\subfigure[Controller \eqref{eq:control with scattering transformation primal-dual} with scattering transformation \eqref{eq:scattering transformation primal-dual}, \eqref{eq:scattering variables under delays primal-dual} under delays $T_{ij} \in (0.1, 1)s$.]{\includegraphics[width = 0.49\linewidth]{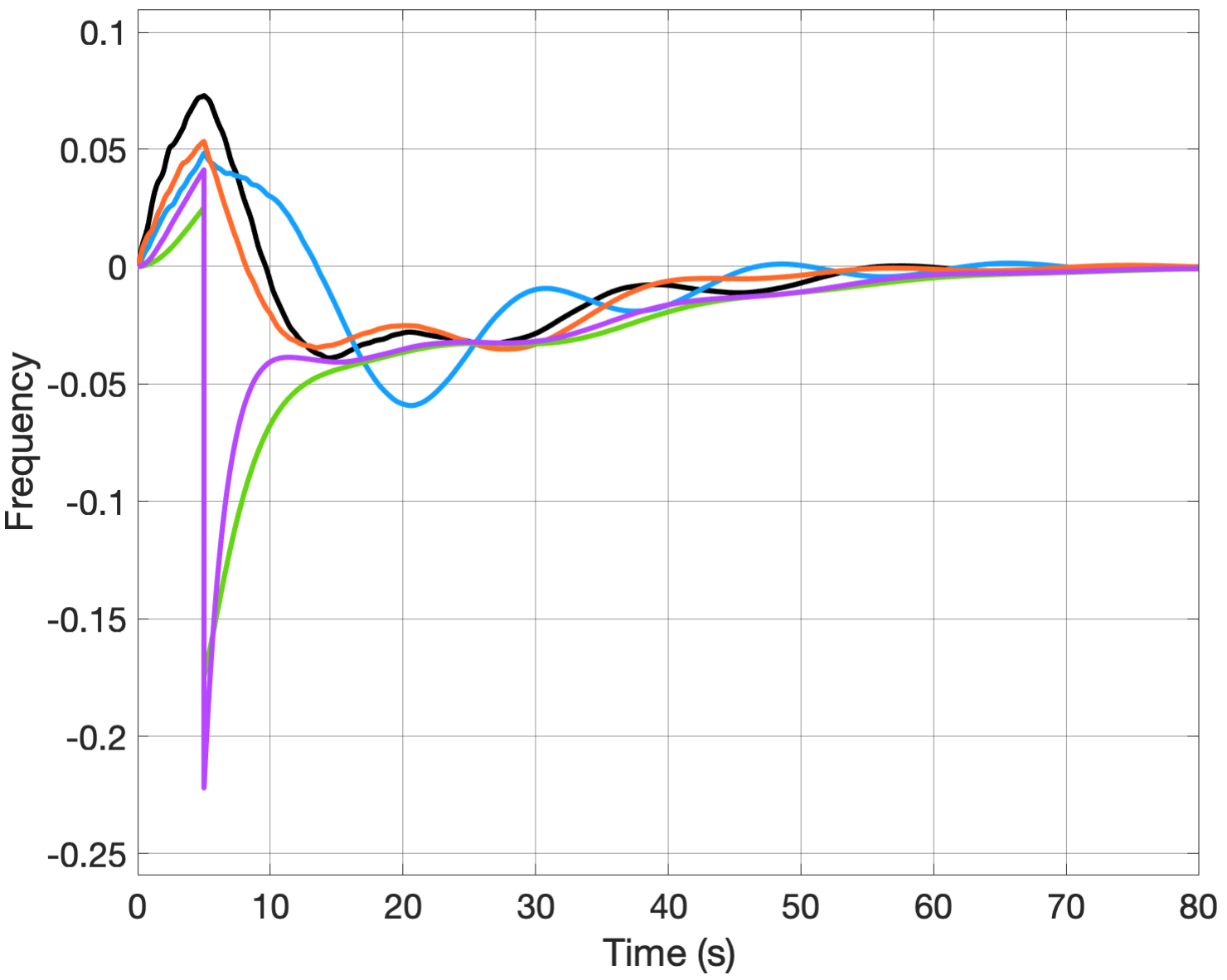}\label{fig:example 2}}
	\subfigure[Controller \eqref{eq:controller without delays tie-line primal-dual} without delays.]{\includegraphics[width = 0.49\linewidth]{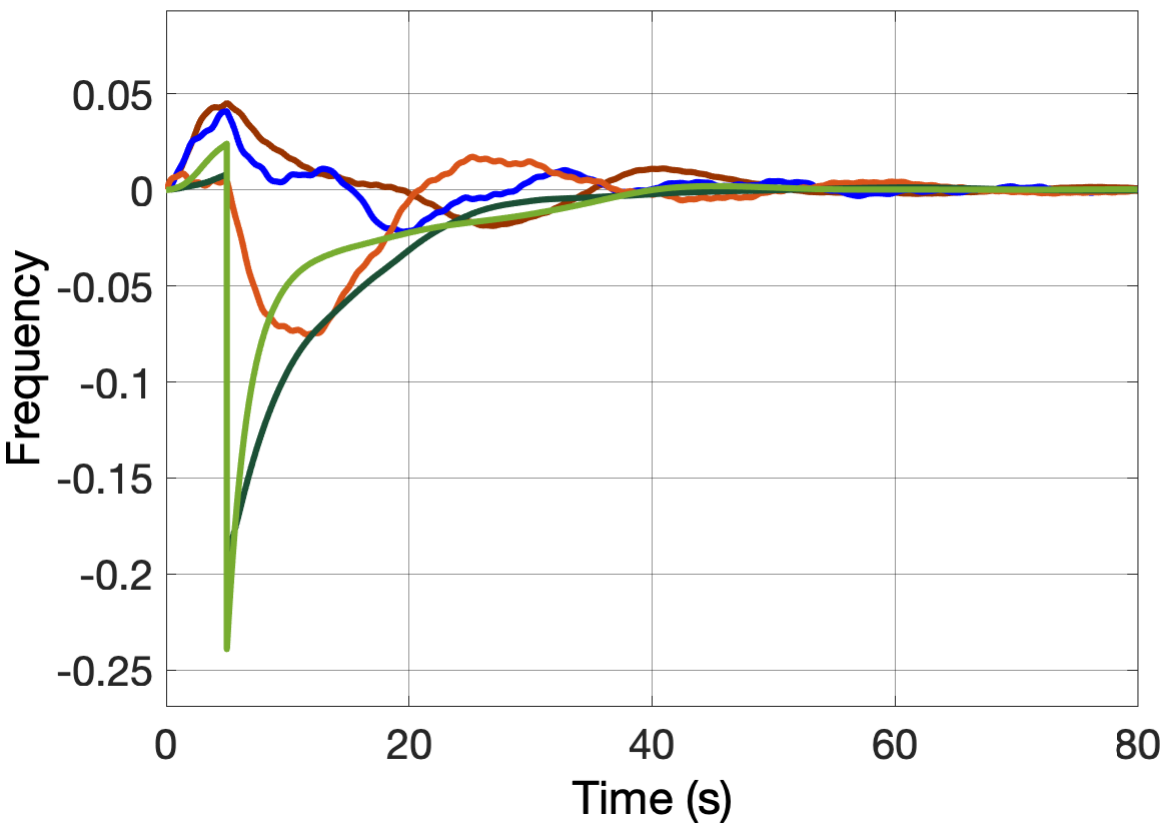}\label{fig:example 3}}
	\subfigure[Controller \eqref{eq:controller without delays tie-line primal-dual} without delays.]{\includegraphics[width = 0.49\linewidth]{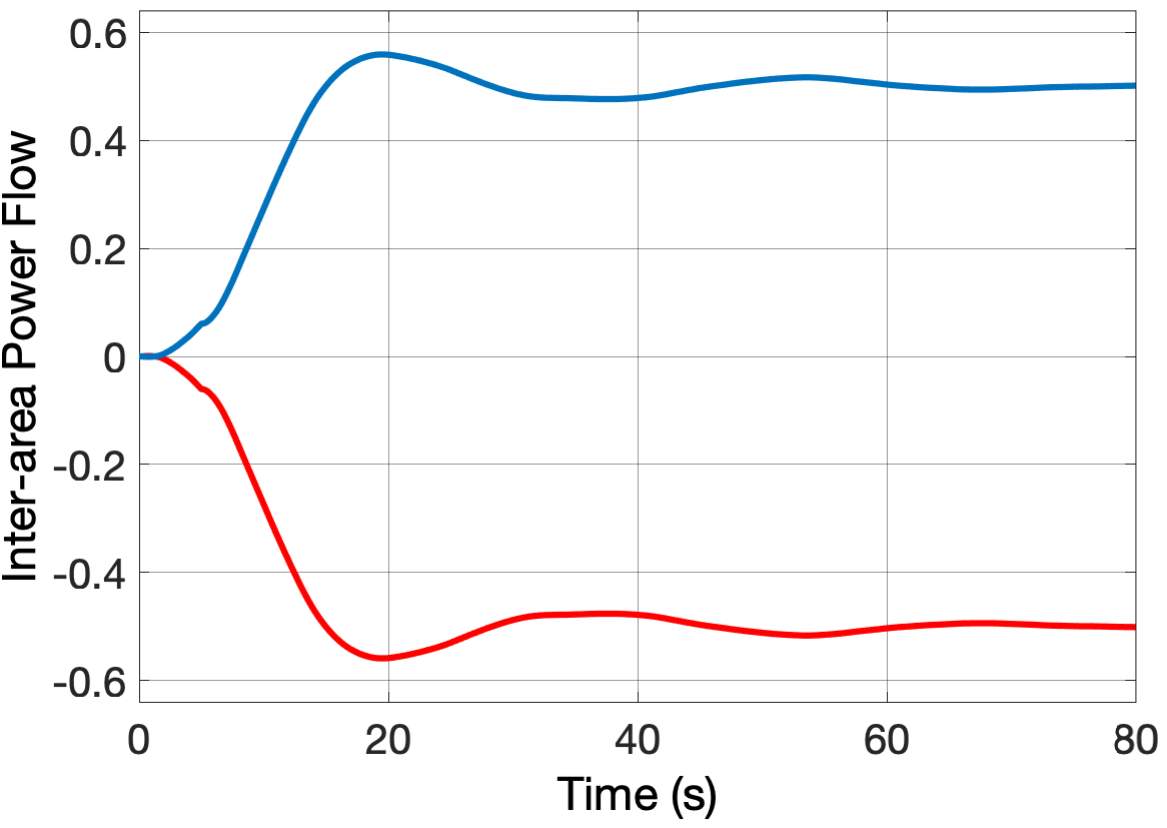}\label{fig:example 4}}
	\subfigure[Extended controller  \eqref{eq:control with scattering transformation tie-line primal-dual} with scattering transformation \eqref{eq:scattering transformation tie-line primal-dual}, \eqref{eq:scattering variables under delays tie-line primal-dual} under delays $T_{ij} \in (0.1, 1) s$.]{\includegraphics[width = 0.49\linewidth]{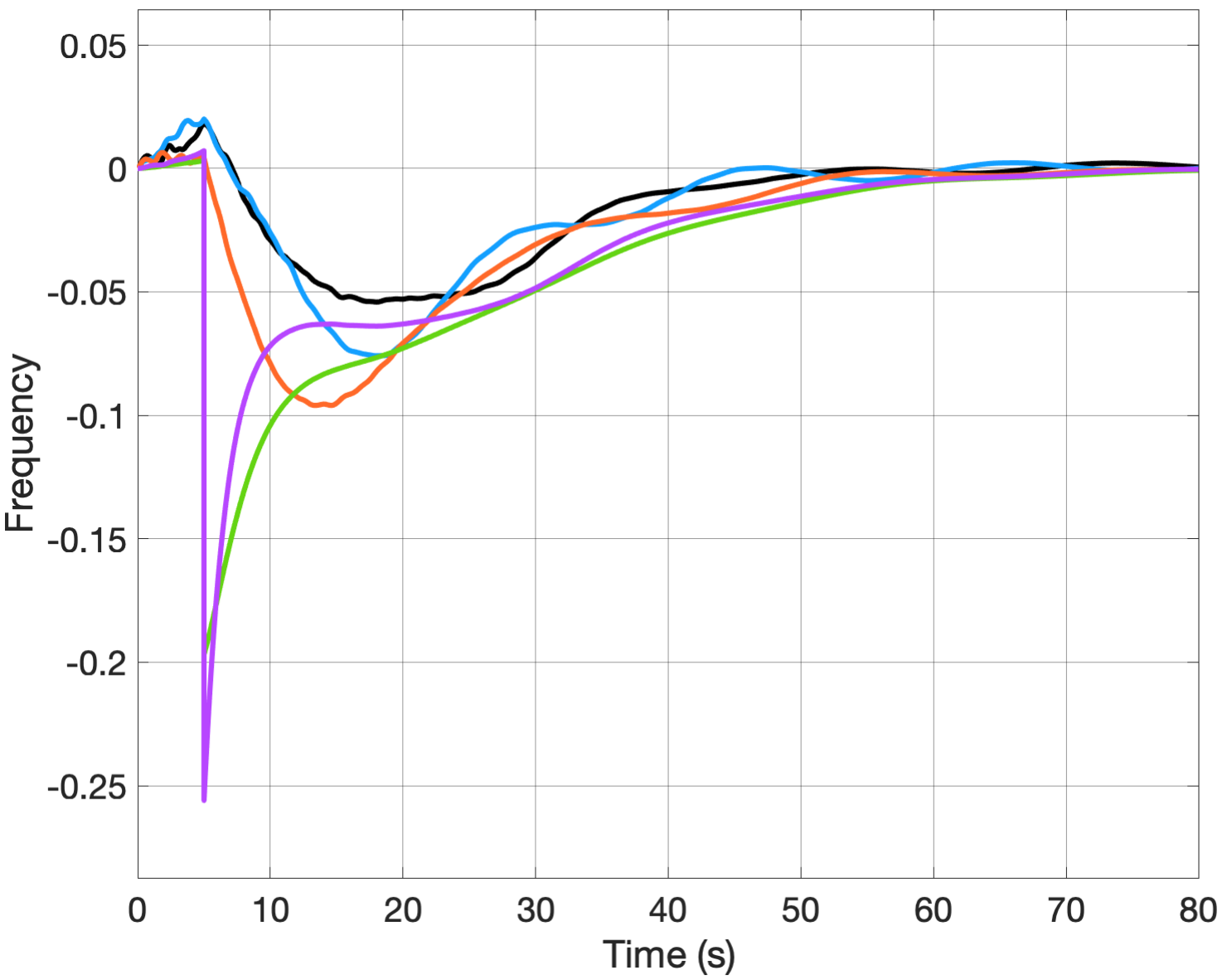}\label{fig:example 5}}
	\subfigure[Extended controller \eqref{eq:control with scattering transformation tie-line primal-dual} with scattering transformation \eqref{eq:scattering transformation tie-line primal-dual}, \eqref{eq:scattering variables under delays tie-line primal-dual} under delays $T_{ij} \in (0.1, 1) s$.]{\includegraphics[width = 0.49\linewidth]{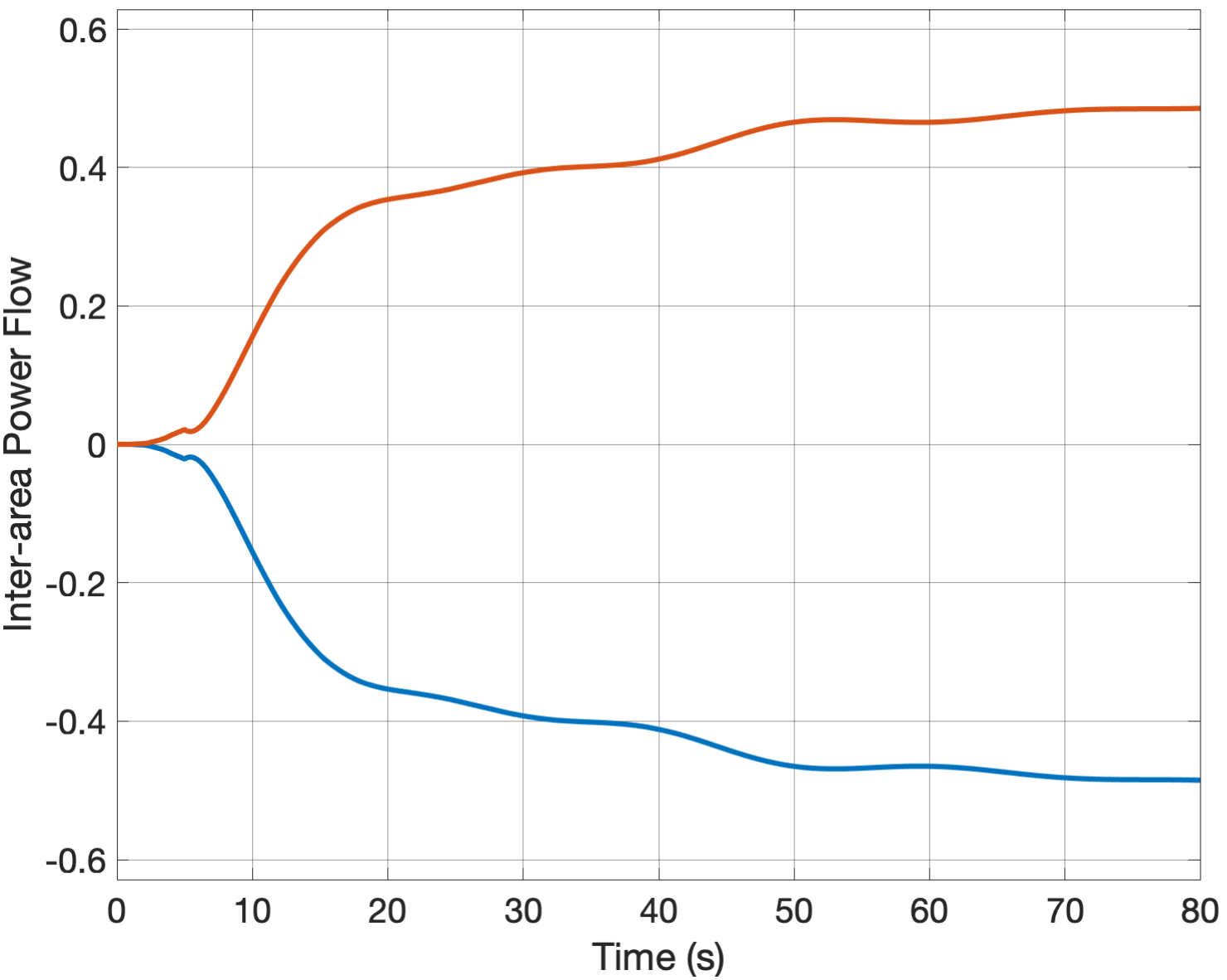}\label{fig:example 6}}
	\caption{Frequencies and inter-area power flows for the primal-dual controllers.}
	\label{fig:example}
\end{figure}

Next, the \ic{performance} of the reformulated primal-dual control algorithm \eqref{eq:controller without delays reformulation primal-dual} and its extension \eqref{eq:controller without delays tie-line primal-dual}  considering tie-line power flow \il{constraints} are shown in \cref{fig:example}.
It can be observed from \cref{fig:example 7,fig:example 8} that the original controllers \eqref{eq:controller without delays reformulation primal-dual} and \eqref{eq:controller without delays tie-line primal-dual} are sensitive to even small delays. \il{On the other hand, when combined} with the scattering transformation, the primal-dual controllers can deal with arbitrary bounded, heterogeneous, and unknown delays.

\begin{figure}[htbp]
	\centering
	\subfigure[$d_{ij} = \frac{1}{e_{ij} + t}$ without scattering transformation.]{\includegraphics[width = 0.49\linewidth]{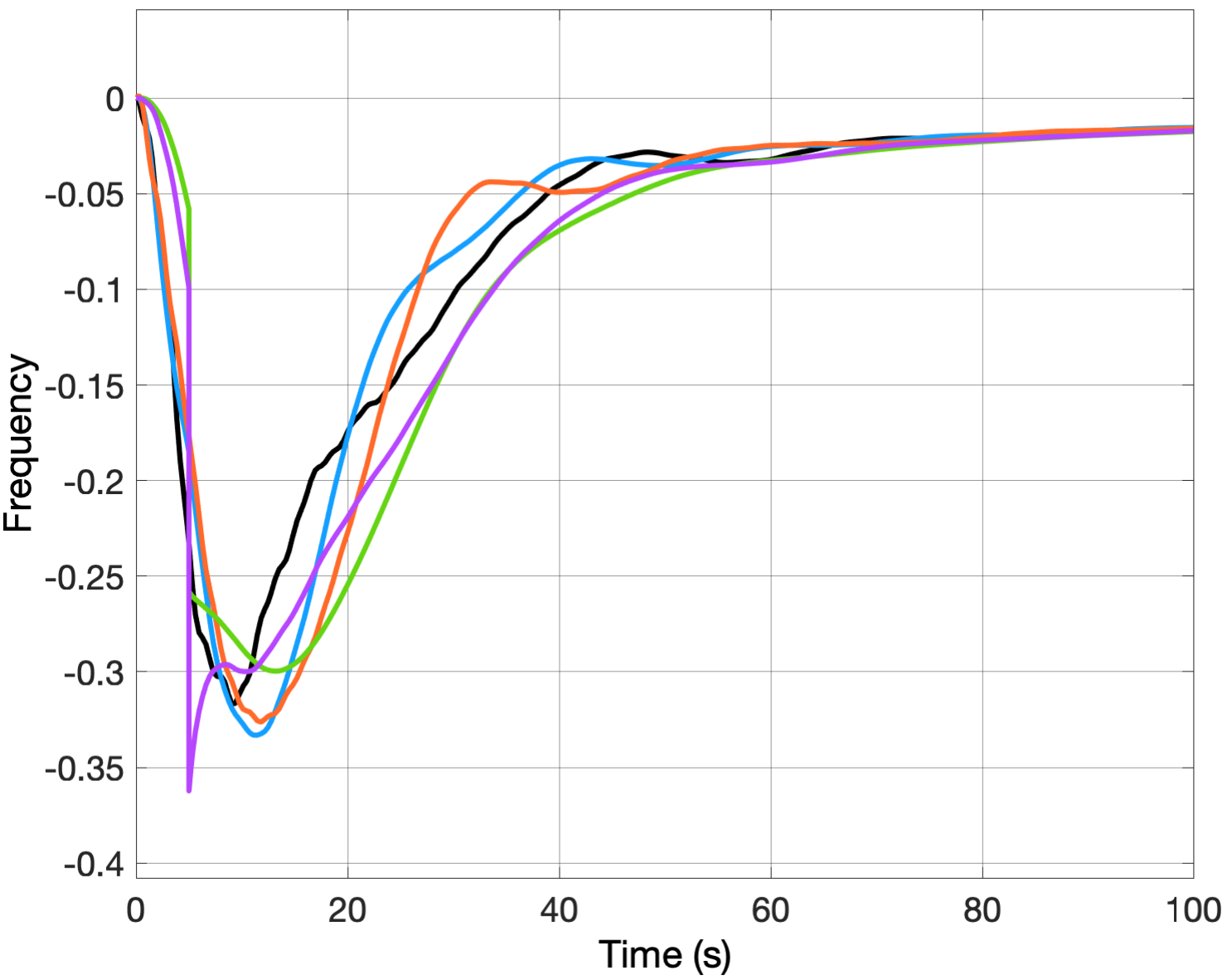}\label{fig:example_disturbance_1}}
	\subfigure[$d_{ij} = \frac{1}{e_{ij} + t}$ with scattering transformation.]{\includegraphics[width = 0.49\linewidth]{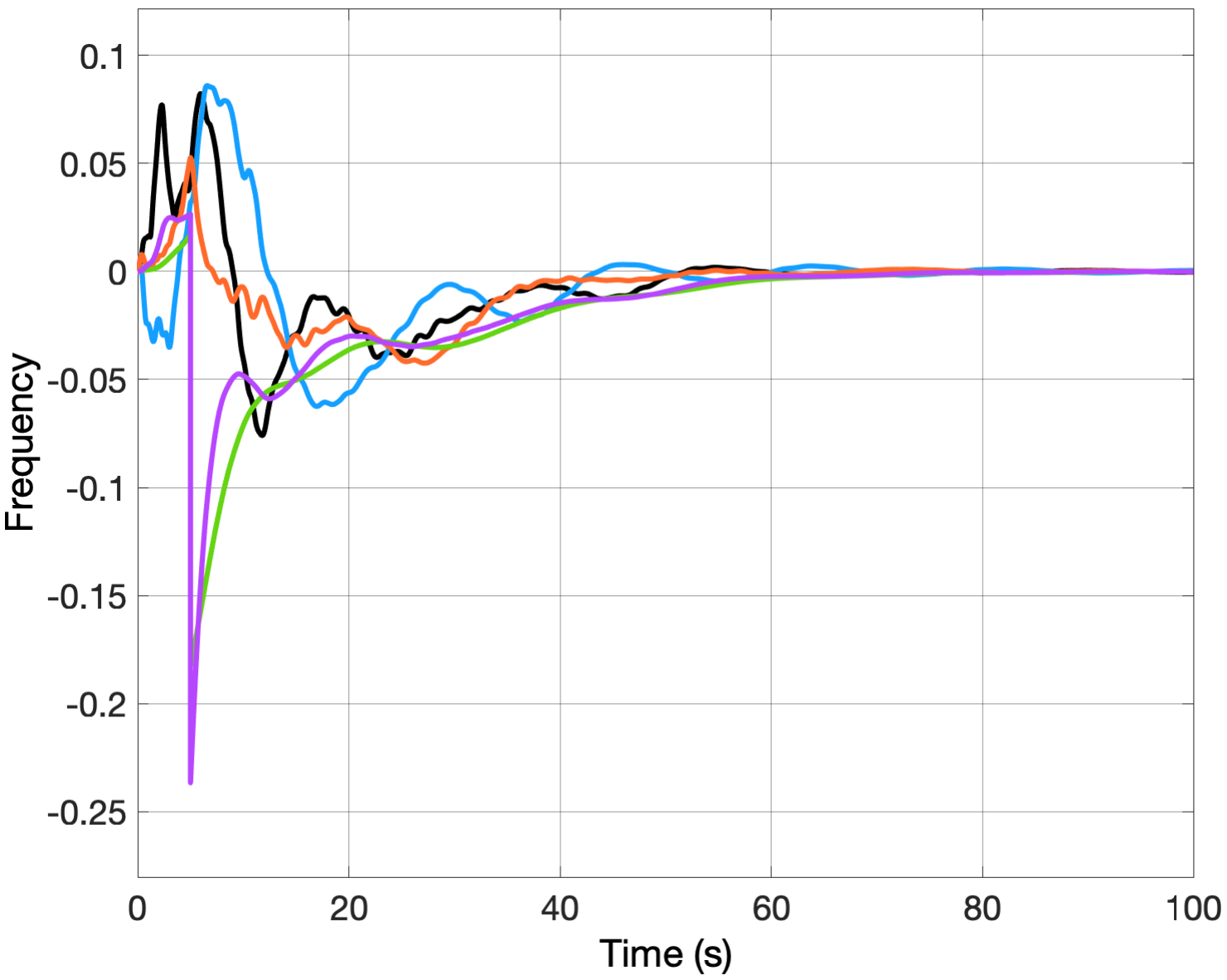}\label{fig:example_disturbance_2}}
	\subfigure[$d_{ij}$ is AWGN with noise power $0.01$, without scattering transformation.]{\includegraphics[width = 0.49\linewidth]{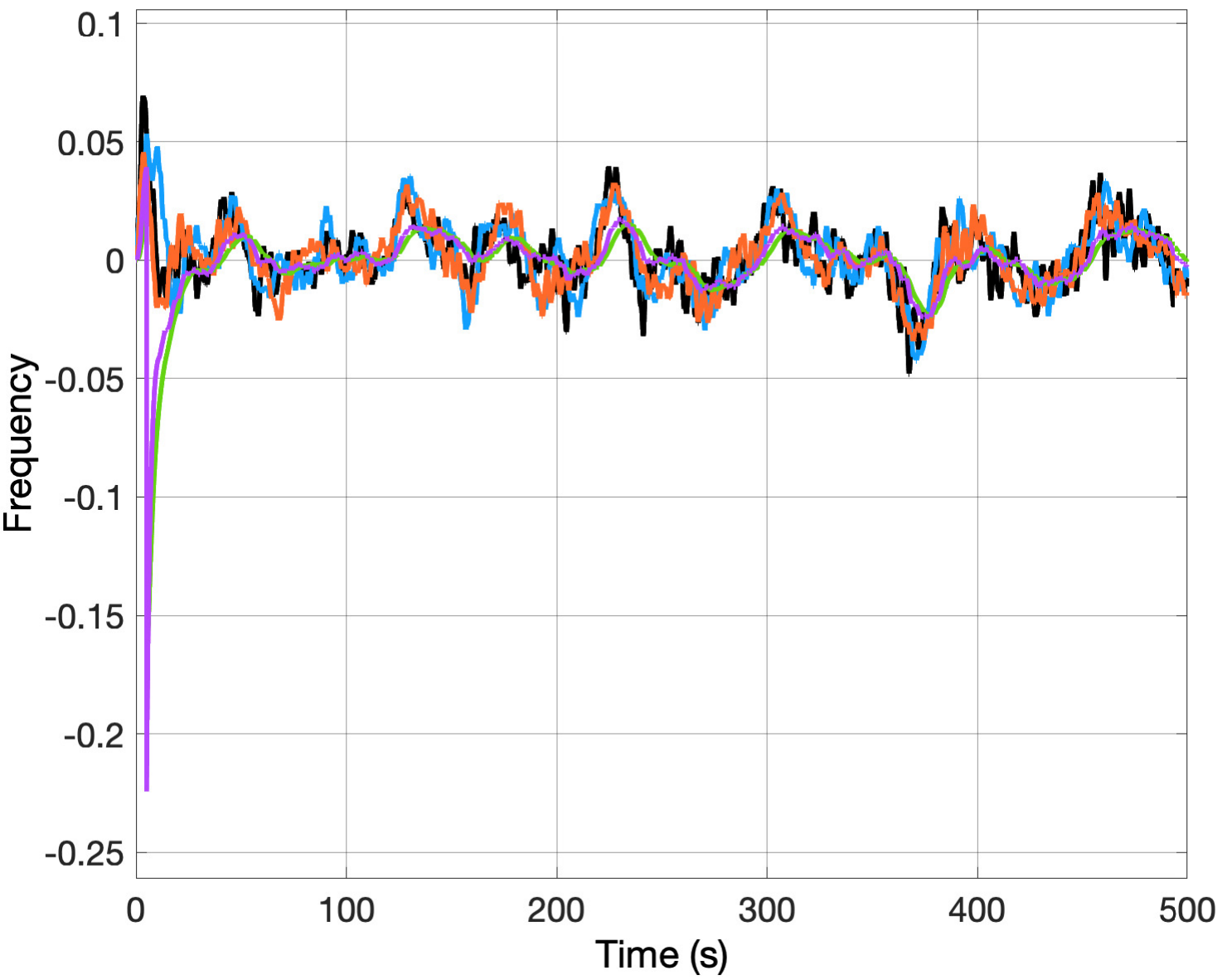}\label{fig:example_disturbance_3}}
	\subfigure[$d_{ij}$ is AWGN with noise power $0.01$, with scattering transformation.]{\includegraphics[width = 0.49\linewidth]{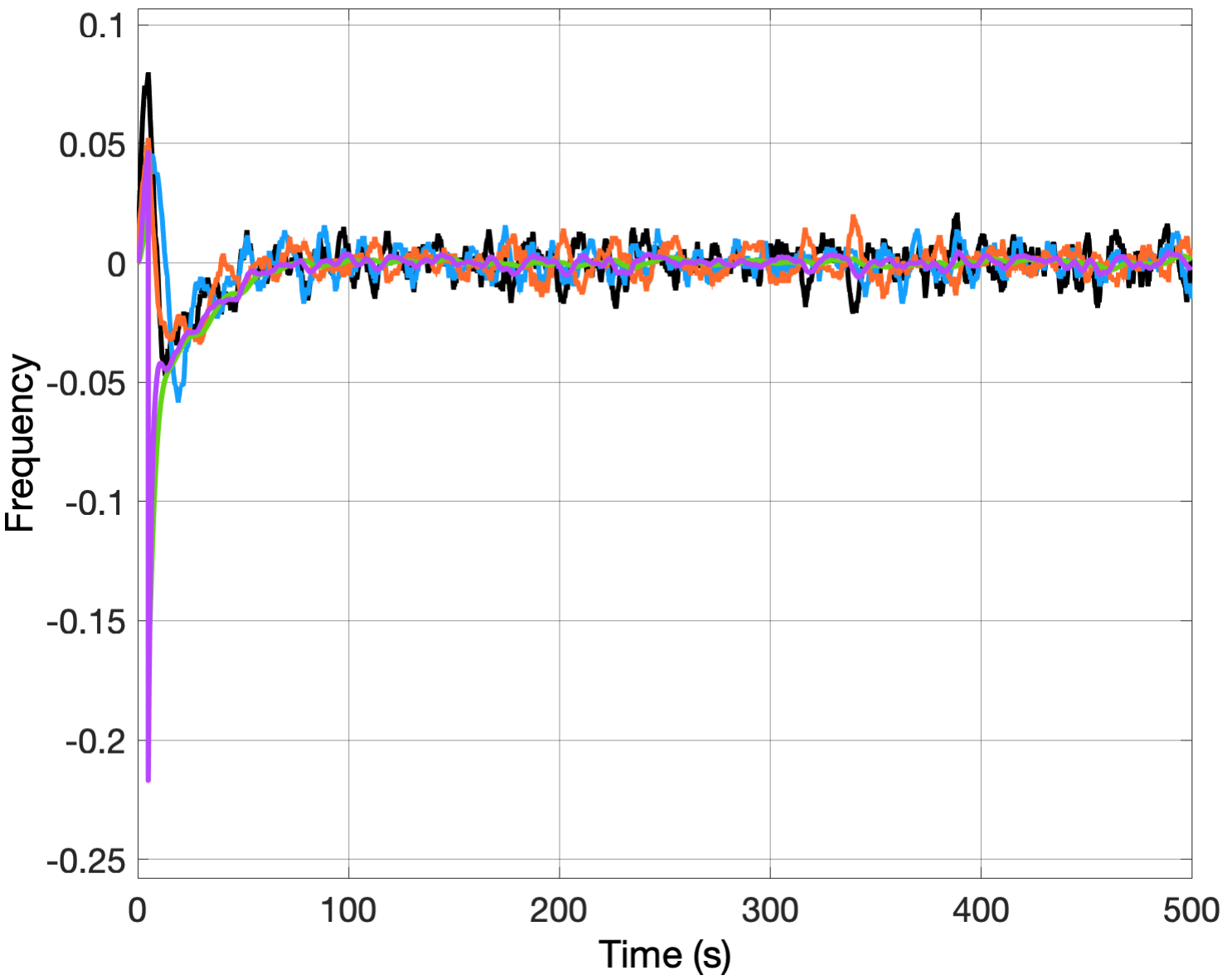}\label{fig:example_disturbance_4}}
	\caption{Frequencies for the primal-dual controller \eqref{eq:controller without delays reformulation primal-dual} under %interconnected
\ic{disturbances in the communicated \il{variables  \mm{described by \eqref{eq:disturbance in scattering variables} and \eqref{eq:disturbance in variables}, respectively, with (right) and without (left)} scattering transformation.}}}
	\label{fig:example_disturbance}
%\todoiny{increase label fonts in the figures on the right}
%\todoiny{make references to equations in the labels}
\end{figure}

In addition, we show that \ic{the control policy with the} scattering transformation \ic{implemented
%is worth implementing as algorithms would have
leads to} better disturbance rejection properties \il{also} %even
when delays are \mm{\ic{small}}.
Consider the same $5$-bus power network and assume that there exist disturbances in the communicating between buses, \ic{with $d_{ij}$ representing} the external disturbance in the channel $i \rightarrow j$.
Then, we have
\begin{align}\label{eq:disturbance in scattering variables}
s_{\overleftarrow{ij}}(t) = E_s s_{\overrightarrow{ji}}(t - T_{ji})+ d_{ji}, ~ s_{\overleftarrow{ji}}(t) = E_s s_{\overrightarrow{ij}}(t - T_{ij}) + d_{ij}
\end{align}
in \eqref{eq:scattering variables under delays primal-dual} for the primal-dual controller incorporated with scattering transformation,
and
\begin{align}\label{eq:disturbance in variables}
\begin{bmatrix}
p_i^c & \zeta_i
\end{bmatrix}^T
\leftarrow
\begin{bmatrix}
p_i^c & \zeta_i
\end{bmatrix}^T
+
d_{ij}
\begin{bmatrix}
1 & 1
\end{bmatrix}^T
\end{align}
in \eqref{eq:controller without delays reformulation primal-dual} for the reformulated primal-dual controller without scattering transformation.
The \il{frequency response at each bus} under square integrable disturbance $d_{ij} = \frac{1}{e_{ij} + t}$, $e_{ij} \in (0,1)$ are shown \ic{with and without the scattering transform} in \cref{fig:example_disturbance_1,fig:example_disturbance_2}, \ic{respectively}, where the former reacts slowly and the latter quickly restores \il{the frequency} to the nominal value.
 The frequency responses when $d_{ij}$ is AWGN with noise power $0.01$ are shown in \cref{fig:example_disturbance_3,fig:example_disturbance_4}.
Note that we only show the disturbance rejection performance of \eqref{eq:controller without delays reformulation primal-dual} instead of the original one \eqref{eq:controller without delays primal-dual} \ml{which is unstable} for even very small disturbance $d_{ij}$ added to $p_j^c$.
Moreover, we have tried various parameters for the algorithm without scattering transformation,  and these do %which does
not give significant improvement on the performance \ic{in}  \cref{fig:example_disturbance_1,fig:example_disturbance_3}.
Overall, the primal-dual controller with scattering transformation \ic{provides} better performance under various \ic{types of} external disturbances.

\subsection*{Example 2}
\begin{figure}[htbp]
	\centering
	\subfigure[Controller \eqref{eq:controller under delays reformulation primal-dual} without scattering transformation under delay $T_{ij} = 0.02s$.]{\includegraphics[width = 0.49\linewidth]{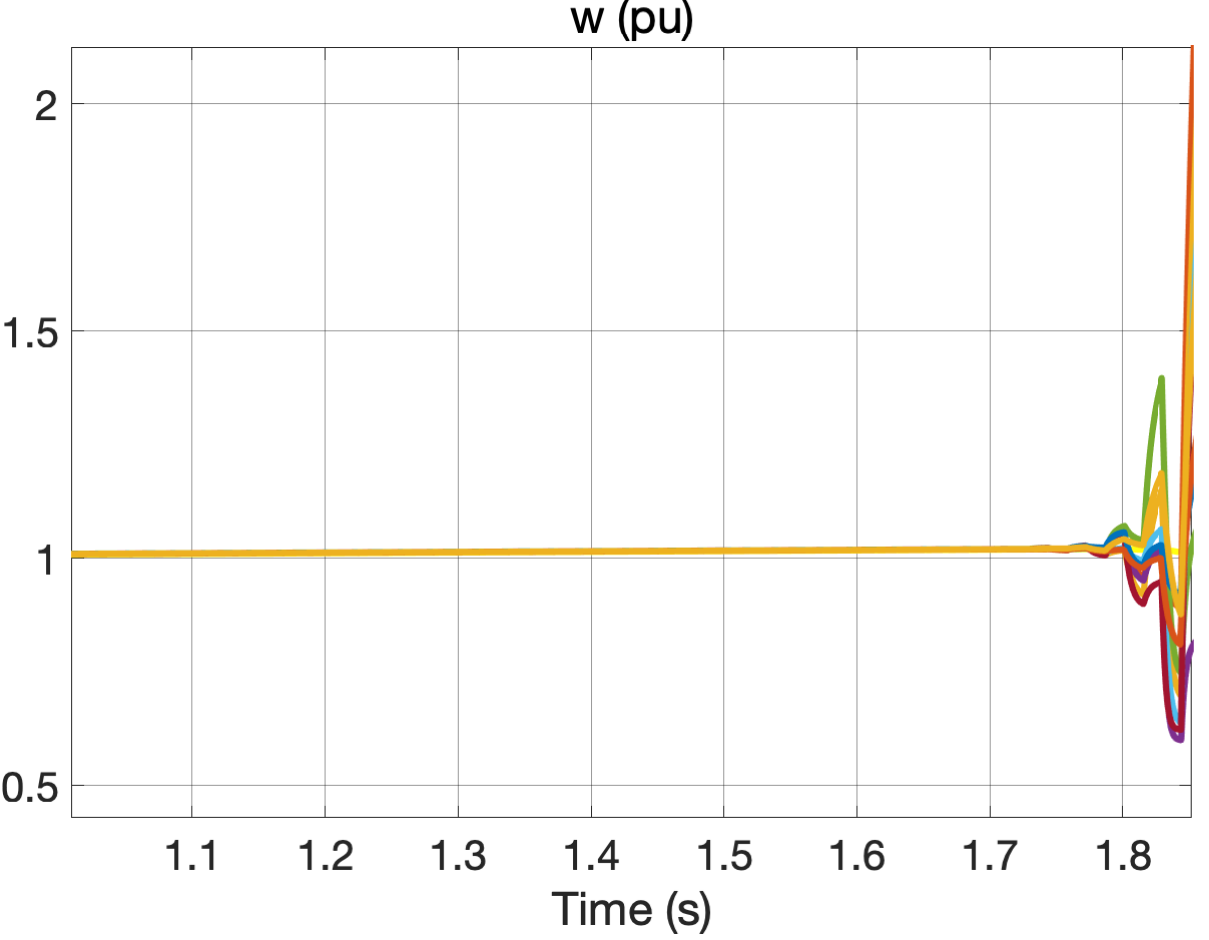}\label{fig:NE39noscattering}}
	\subfigure[Controller \eqref{eq:control with scattering transformation primal-dual} with scattering transformation \eqref{eq:scattering transformation primal-dual}, \eqref{eq:scattering variables under delays primal-dual} under delay $T_{ij} = 0.02s$.]{\includegraphics[width = 0.49\linewidth]{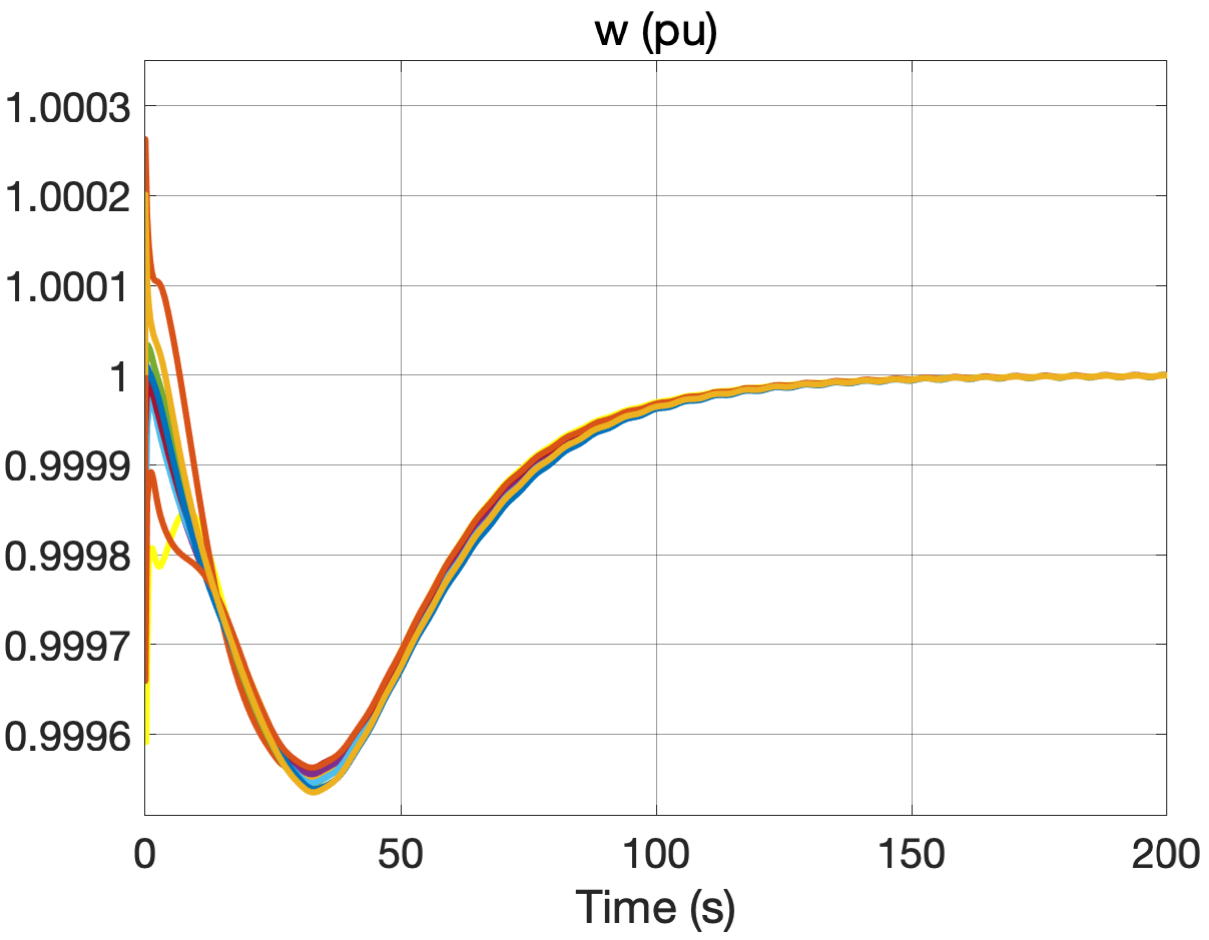}\label{fig:NE39scattering}}
	\caption{Frequencies of the ten machines of the IEEE-39 test system under various secondary controllers and delay $T_{ij} = 0.02s$.}
	\label{fig:ne39example}
\end{figure}

\begin{figure}[htbp]
	\centering
	\subfigure[Controller \eqref{eq:controller without delays tie-line primal-dual} without delays.]{\includegraphics[width = 0.49\linewidth]{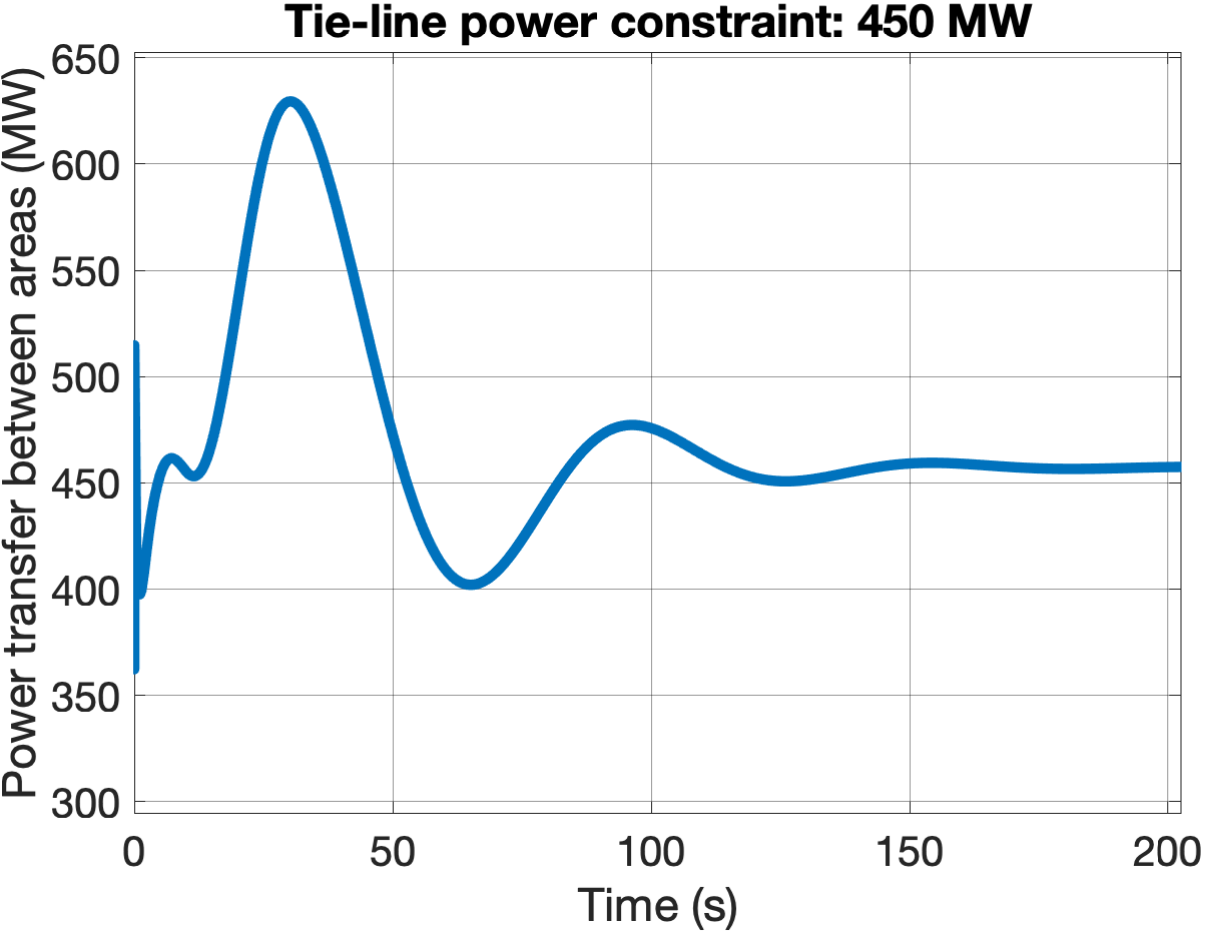}\label{fig:ne39example_tieline_undelay}}
	\subfigure[Extended controller \eqref{eq:control with scattering transformation tie-line primal-dual} with scattering transformation \eqref{eq:scattering transformation tie-line primal-dual}, \eqref{eq:scattering variables under delays tie-line primal-dual} under delays $T_{ij} = 0.02 s$.]{\includegraphics[width = 0.49\linewidth]{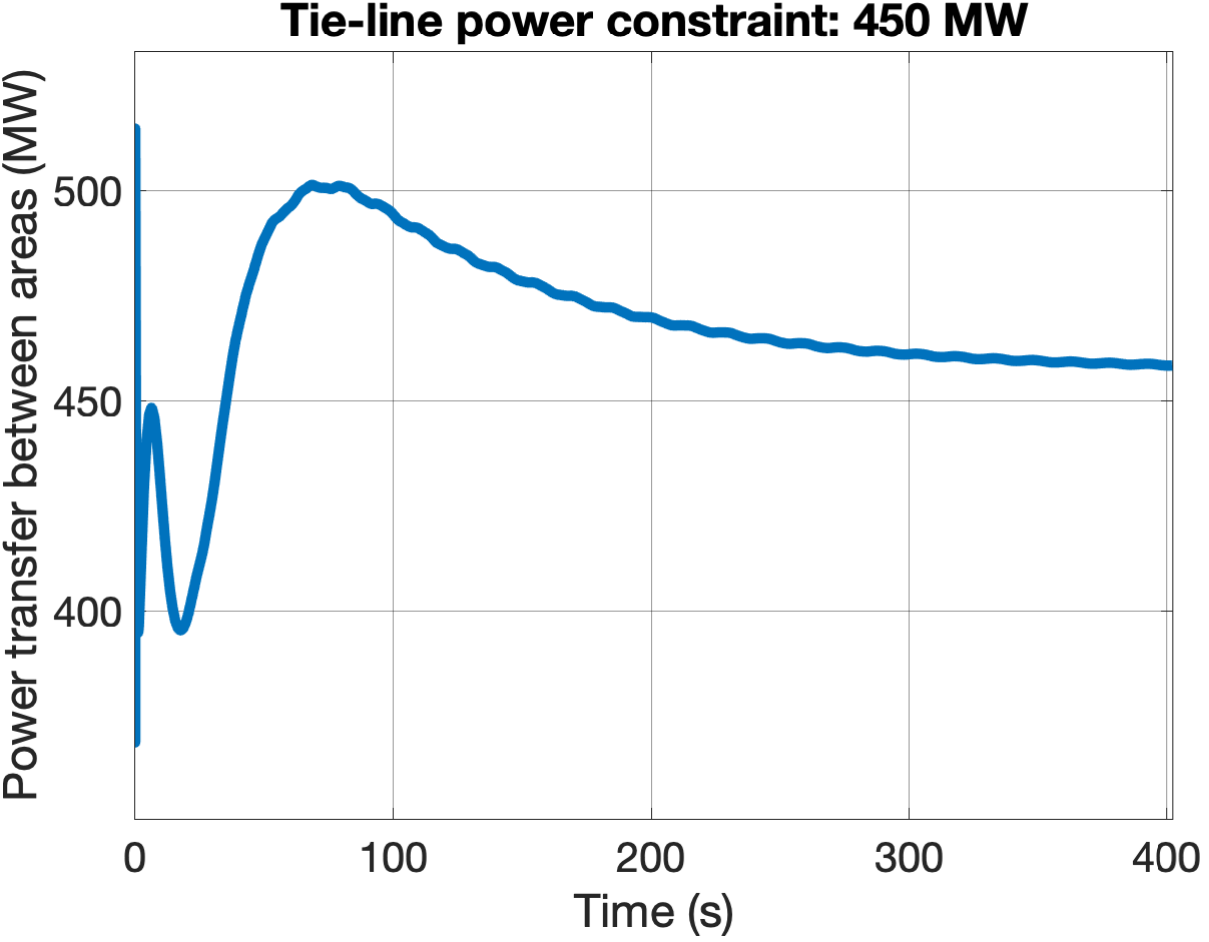}\label{fig:ne39example_tieline_delayed}}
	\caption{Inter-area power flows for the primal-dual controllers in the IEEE 39 bus system.}
	\label{fig:ne39example_tieline}
\end{figure}

We apply the proposed secondary frequency control algorithms on the well-known IEEE New England 39-bus system \cite{newengland39}. The model is highly detailed and includes high-order models of the generators, turbine-governors, exciters, transformers and lines. We compare the controller \eqref{eq:control with scattering transformation primal-dual} with scattering transformation \eqref{eq:scattering transformation primal-dual}, \eqref{eq:scattering variables under delays primal-dual} to the primal-dual scheme \eqref{eq:controller under delays reformulation primal-dual} under a uniform delay of $T_{ij} = 0.02s$. \cref{fig:ne39example} shows that the controller with scattering transformation is able to tolerate this small delay, unlike the original primal-dual scheme without scattering transform.
We also show the ability of the primal-dual controller with scattering transformation to regulate inter-area flows in the presence of delays. To demonstrate \il{this}, we designate buses 21-24, 35 and 36 as Area 2 (with the rest of the network as Area 1) and set the total desired inter-area power flow (over two tie-lines) to 450 MW (without the constraint, the steady-state inter-area power-flow is 507 MW). It can be seen from \cref{fig:ne39example_tieline} that the inter-area power-flow is successfully regulated to (very close to) the desired value and the optimal solution to the \textup{OGR-2} problem \eqref{optimal generation regulation problem 2} is reached. These examples on a realistic and highly detailed model serve to verify the effectiveness of the proposed primal-dual controller with scattering transformation.

\section{Conclusion}\label{sec:Conclusion}
This work has proposed primal-dual controllers with delay independent stability for \icl{distributed} secondary frequency control \ic{in} power systems with unknown and heterogeneous constant communication delays.
An equivalent passive reformulation of the controller has been derived and a novel form of passivity-based scattering transformation has been constructed to robustify %communication channels
\ic{the closed loop dynamics}
against \ic{communication} delays.
Moreover, it has been shown that the proposed \ic{controllers with delay independent stability properties} %controllers
%inherit good extensibility of the conventional
\ic{can adopt various extensions associated with
primal-dual control schemes that allow to incorporate various operational constraints. These include} %such as
extra tie-line power flow and generation boundedness constraints, \ic{and a relaxation of the requirement for demand measurements via an observer}.

\appendices
\section*{Appendix}\label{Appendix}

\subsection{Proof of Lemma~\ref{lem:convergence without delays primal-dual} }\label{appendix proof of lem convergence without delays primal-dual}
\begin{proof}
	Consider the Lyapunov function candidate
	\begin{equation}\label{eq:V_pd}
		V_{pd} = V_F \left(\omega^G \right) + V_P(\eta) + V_C \left(p^c\right) + V_{\psi}\left( \psi \right) + V_D \left( p^M \right)
	\end{equation}
	where the functions on the right hand side are defined by
\begin{align}
	V_F(\omega^{G}) = \frac{1}{2} \sum_{j \in G} M_j \left( \omega_j - \omega_j^* \right)^2 \geq 0, \label{eq:V_F}\\
	V_P(\eta) = \sum_{(i,j) \in E} Y_{ij} \int_{\eta_{ij}^*}^{\eta_{ij}}\left( \sin \theta - \sin \eta_{ij}^* \right) d \theta, \label{eq:V_P}\\
	V_C (p^c) = \frac{1}{2} \sum_{j \in N} \gamma_j \left( p_j^c - p_j^{c,*}\right)^2 \geq 0, \label{eq:V_C}\\
	V_{\psi} (\psi) = \sum_{(i,j) \in \tilde{E}} \frac{1}{2} \gamma_{ij} \left( \psi_{ij} - \psi_{ij}^* \right)^2 \geq 0,\\
	V_D (p^M) = \sum_{j \in G} \frac{\tau_j}{2 k_{g,j} k_{c,j}} \left( p_j^M - p_j^{M,*}\right)^2 \geq 0, \label{eq:V_D}
	\end{align}
respectively, and $V_P \geq 0$ in some neighborhood of $\eta^*$ \ic{as follows from Assumption} \ref{assumption angle}.
The time derivative of $V_{pd}$ along the system trajectories satisfies
\begin{align*}
\begin{array}{rl}
	\hspace{-2mm} \dot{V}_{pd} = \hspace{-3mm} & -\sum_{j \in N} \Lambda_j \left( \omega_j - \omega_j^* \right)^2\\
	& - \sum_{j \in G} \left(Q_j'\left( p_j^M\right) - Q_j'\left( p_j^{M,*}\right)\right) \left( p_j^M - p_j^{M,*}\right) \leq 0.
\end{array}
\end{align*}
\mm{The rest of the proof can be found in \cite{kasis2019stability}.}
\end{proof}

\subsection{Proof of Corollary~\ref{cor:convergence primal-dual}}\label{appendix proof of cor convergence primal-dual}
Consider the Lyapunov function candidate
\begin{align*}
	V_{N} = V_F \left(\omega^G \right) + V_P(\eta) + V_C \left(p^c\right) + V_D \left(p^M \right) + V_{\zeta}(\zeta)
\end{align*}
where $V_F \left(\omega^G \right)$,  $V_P(\eta)$, $V_C \left(p^c\right)$, and $V_D \left(p^M \right)$ are defined by \eqref{eq:V_F}, \eqref{eq:V_P}, \eqref{eq:V_C}, and \eqref{eq:V_D}, respectively, and $V_{\zeta} = \frac{1}{2} \sum_{j \in N} \left( \zeta_j - \zeta_j^*\right)^2$.
The rest of the proof are similar to the proof of \Cref{lem:convergence without delays primal-dual}.

\subsection{Proof of Lemma~\ref{lem:passivity of reformulated controller} }\label{appendix proof of lem passivity of reformulated controller}
\begin{proof}
	Adopt the storage function
\begin{equation}\label{eq:V_B}
V_B = V_F \left(\omega^G \right) + V_P(\eta) + \ml{V_{p^c, \rho}\left(p^c, \rho^{p}\right)} + V_D \left(p^M \right) + V_{\zeta,\rho}(\zeta,\rho^{\zeta})
\end{equation}
where $V_F \left(\omega^G \right)$,  $V_P(\eta)$ and $V_D \left(p^M \right)$ are defined by \eqref{eq:V_F}, \eqref{eq:V_P}, and \eqref{eq:V_D}, respectively, and
\begin{align*}
\begin{array}{rl}
	V_{p^c, \rho} (p^c, \rho^{p})=\displaystyle  \frac{1}{2} (\rho_j^{p})^2 + \frac{1}{2} \sum_{j \in N} \left( p_j^c - \rho_j^{p} - p_j^{c,*} \right)^2, \\
	V_{\zeta,\rho}(\zeta, \rho^{\zeta}) = \displaystyle \frac{1}{2} (\rho_j^{\zeta})^2 + \frac{1}{2} \sum_{j \in N} \left( \zeta_j - \rho_j^{\zeta} - \zeta_j^* \right)^2.
\end{array}
\end{align*}
By Assumption \ref{assumption angle}, there exists an open neighborhood of $\eta^*$ such that $V_B \geq 0$. The time derivative of $V_B$ along the system trajectories \ic{is given by}
\begin{align*}
\begin{array}{rl}
	\dot{V}_B
%	= & (\omega - \omega^*)^T \left( \left(p^M - p^{M,*}\right) + D (p - p^*) \right)\\
%	& + (p - p^*)^T D^T (\omega - \omega^*) - (p^c - p^{c,*})^T  \left( p^M - p^{M,*}\right)\\
%	& + \left( p^M - p^{M,*} \right)^T \left( (p^c - p^{c,*}) - (\omega - \omega^*)\right)\\
%	& -\sum_{j \in N} \left( \Lambda_j \left( \omega_j - \omega_j^* \right)^2 + \ml{(\rho_j^{p})^2 + (\rho_j^{\zeta})^2} \right) \\
%	& - \sum_{j \in G} \left(p_j^M - p_j^{M,*}\right) \left( Q_j'\left( p_j^M \right) - Q_j'\left( p_j^{M,*} \right)\right)\\
%	& - \sum_{j \in N} (p_j^c - p_j^{c,*}) \sum_{i \in \tilde{N}_j} \alpha_{ij} \left( \left(r_{ij}^{\zeta} - r_{ij}^{\zeta,*} \right) - \left(\zeta_j - \zeta_j^*\right) \right) \\
%	& + \sum_{j \in N} ( \zeta_j - \zeta_j^*) \sum_{i \in \tilde{N}_j} \alpha_{ij} \left( \left(r_{ij}^{p} - r_{ij}^{p,*}\right) - \left(p_j^c - p_j^{c,*}\right) \right)\\
	= \hspace{-3mm} & -\sum_{j \in N}  \left(  \Lambda_j \left( \omega_j - \omega_j^* \right)^2 + \ml{(\rho_j^{p})^2 + (\rho_j^{\zeta})^2} \right)\\
	& - \sum_{j \in G} \left(p_j^M - p_j^{M,*}\right) \left( Q_j'\left( p_j^M \right) - Q_j'\left( p_j^{M,*} \right)\right)\\
	& - \sum_{j \in N} (p_j^c - p_j^{c,*}) \sum_{i \in \tilde{N}_j} \alpha_{ij} \left( \left(r_{ij}^{\zeta} - r_{ij}^{\zeta,*} \right)\right)\\
	& + \sum_{j \in N} ( \zeta_j - \zeta_j^*) \sum_{i \in \tilde{N}_j} \alpha_{ij} \left( \left(r_{ij}^{p} - r_{ij}^{p,*}\right) \right)\\
	\leq \hspace{-3mm} & \sum_{j\in N}
	\begin{bmatrix}
		\tilde{\zeta}_j \\ - \tilde{p}_j^c
	\end{bmatrix}^T  \sum_{i \in \tilde{N}_j} \alpha_{ij} \begin{bmatrix}
		\tilde{r}^{p}_{ij} \\ \tilde{r}^{\zeta}_{ij}
	\end{bmatrix}
	=
	\begin{bmatrix}
		\tilde{\zeta} \\ - \tilde{p}^c
	\end{bmatrix}^T
	\begin{bmatrix}
		\tilde{r}^{p} \\ \tilde{r}^{\zeta}
	\end{bmatrix}.
\end{array}
\end{align*}
%In the undelayed case, we have $r_{ij}^{p} = p_i^{c}$, $r_{ij}^{\zeta} = \zeta_i$, and $\dot{V}_B \leq 0$, which leads to the convergence results in \Cref{cor:convergence primal-dual}.
\end{proof}

%\todoiny{\st{I think add a separate proof for Corollary in the same appendix.}}

\subsection{Proof of Lemma \ref{lem:passivity of scattering transformation} }\label{appendix proof of lem passivity of scattering transformation}
\begin{proof}
Define $s_{\overrightarrow{ij}}^* = \frac{1}{\sqrt{2}}\left( - \begin{bmatrix}
		r_{ji}^{p,*} \\ r_{ji}^{\zeta,*}
	\end{bmatrix}  + \begin{bmatrix}
		\zeta_i^* \\ - p_i^{c,*}
	\end{bmatrix} \right)$, $s_{\overrightarrow{ji}}^* = \frac{1}{\sqrt{2}}\left( \begin{bmatrix}
		r_{ji}^{p,*} \\ r_{ji}^{\zeta,*}
	\end{bmatrix}  - \begin{bmatrix}
		\zeta_i^* \\ - p_i^{c,*}
	\end{bmatrix} \right)$ and $s_{\overleftarrow{ij}}^*$, $s_{\overleftarrow{ji}}^*$ are defined similarly. Let the storage functional $V_S^{ij}$ be
\begin{equation}\label{eq:V_S^ij}
\begin{array}{rl}
	\hspace{-4mm} V_S^{ij} \hspace{-1mm} = \hspace{-3mm} & \displaystyle  \frac{1}{2} \int_0^t \left( \left\| E_s s_{\overrightarrow{ij}}(\tau) - E_s s_{\overrightarrow{ij}}^* \right \|^2  \hspace{-0.5mm}  -  \hspace{-0.5mm}  \left\| s_{\overleftarrow{ji}} (\tau) - s_{\overleftarrow{ji}}^* \right \|^2 \right. \\
	&\left. + \left\| E_s s_{\overrightarrow{ji}} (\tau) - E_s s_{\overrightarrow{ji}}^* \right \|^2 - \left\| s_{\overleftarrow{ij}} (\tau) - s_{\overleftarrow{ij}}^* \right \|^2  \right) d \tau
\end{array}
\end{equation}
Recall that $r_{ij}^{p,*} = p_j^{c,*}$, $r_{ij}^{\zeta,*} = \zeta_j^*$, and \mm{\eqref{eq:scattering variables under delays primal-dual}}. Then, we have $E_s s_{\overrightarrow{ij}}^* = s_{\overleftarrow{ji}}^*$, $E_s s_{\overrightarrow{ji}}^* = s_{\overleftarrow{ij}}^*$, and
\begin{align*}
\begin{array}{rl}
	V_S^{ij} = & \hspace{-3mm}  \displaystyle \int_{t - T_{ij}}^t  \frac{1}{2} \left\| E_s s_{\overrightarrow{ij}}(\tau) - E_s s_{\overrightarrow{ij}}^* \right\|^2 d\tau \\
	& \displaystyle + \frac{1}{2} \int_{t - T_{ji}}^t \left\| E_s s_{\overrightarrow{ji}} (\tau) - E_s s_{\overrightarrow{ji}}^* \right\|^2 d \tau \geq 0.
\end{array}
\end{align*}
The time derivative of $V_S^{ij}$ \ic{is given by}
\begin{align*}
\begin{array}{rl}
	\hspace{-2mm} \dot{V}_S^{ij} \hspace{-1mm}
	 =  & \hspace{-3mm}  \frac{1}{2} \left( \left\| E_s \left(s_{\overrightarrow{ij}} - s_{\overrightarrow{ij}}^*\right) \right\|^2 - \left\| s_{\overleftarrow{ji}} - s_{\overleftarrow{ji}}^* \right\|^2 \right. \\
	 & \left. + \left\| E_s \left( s_{\overrightarrow{ji}} - s_{\overrightarrow{ji}}^* \right) \right\|^2 - \left\|s_{\overleftarrow{ij}} - s_{\overleftarrow{ij}}^* \right\|^2\right)\\
= & \hspace{-3mm} \frac{1}{2} \hspace{-1mm} \left( \left\| s_{\overrightarrow{ij}} \hspace{-1mm} - \hspace{-1mm}  s_{\overrightarrow{ij}}^* \right\|^2 \hspace{-2mm}  - \hspace{-1 mm} \left\| s_{\overleftarrow{ji}} \hspace{-1mm}   - \hspace{-1mm}  s_{\overleftarrow{ji}}^* \right\|^2  \hspace{-2mm}
	 + \hspace{-1mm}  \left\| s_{\overrightarrow{ji}} \hspace{-1mm}  - \hspace{-1mm}  s_{\overrightarrow{ji}}^* \right\|^2 \hspace{-2mm} - \hspace{-1mm} \left\|s_{\overleftarrow{ij}} \hspace{-1mm}  - \hspace{-1mm}  s_{\overleftarrow{ij}}^* \right\|^2\right)\\
=  & \hspace{-3mm}  \frac{1}{4}\left( \left\| -  \left[\begin{smallmatrix}
		\tilde{r}_{ji}^{p} \\ \tilde{r}_{ji}^{\zeta}
	\end{smallmatrix}\right]  +  \left[\begin{smallmatrix}
		\tilde{\zeta}_i \\ - \tilde{p}_i^c
	\end{smallmatrix}\right]  \right\|^2 - \left\| - \left[\begin{smallmatrix}
		\tilde{r}_{ji}^{p} \\ \tilde{r}_{ji}^{\zeta}
	\end{smallmatrix}\right]  -  \left[\begin{smallmatrix}
		\tilde{\zeta}_i \\ - \tilde{p}_i^c
	\end{smallmatrix}\right]  \right\|^2 \right.\\
	& + \left. \left\|  \left[\begin{smallmatrix}
		\tilde{r}_{ji}^{p} \\ \tilde{r}_{ji}^{\zeta}
	\end{smallmatrix}\right] - \left[\begin{smallmatrix}
		\tilde{\zeta}_i \\ - \tilde{p}_i^c
	\end{smallmatrix}\right] \right\|^2 - \left\|  \left[\begin{smallmatrix}
		\tilde{r}_{ji}^{p} \\ \tilde{r}_{ji}^{\zeta}
\end{smallmatrix}\right]  +  \left[\begin{smallmatrix}
		\tilde{\zeta}_i \\ - \tilde{p}_i^c
	\end{smallmatrix}\right] \right\|^2 \right)\\
	=  & \hspace{-3mm}  -  \left[\begin{smallmatrix} \tilde{r}_{ij}^{p} &  \tilde{r}_{ij}^{\zeta} \end{smallmatrix}\right]  \left[\begin{smallmatrix} \tilde{\zeta}_j \\ - \tilde{p}_j^c \end{smallmatrix}\right] - \left[\begin{smallmatrix} \tilde{r}_{ji}^{p} & \tilde{r}_{ji}^{\zeta} \end{smallmatrix}\right]  \left[\begin{smallmatrix} \tilde{\zeta}_i \\ - \tilde{p}_i^c \end{smallmatrix}\right]
\end{array}
\end{align*}
\ml{where the second equality holds since $\|E_s\|^2 = 1$.}
It is worth noting that
%the equality holds, i.e., $\dot{V}_{S}^{ij} = -\begin{bmatrix} \tilde{r}_{ij}^{p} \\ \tilde{r}_{ij}^{\zeta} \end{bmatrix}^T \begin{bmatrix} \tilde{\zeta}_j \\ - \tilde{p}_j^c \end{bmatrix} - \begin{bmatrix} \tilde{r}_{ji}^{p} \\ \tilde{r}_{ji}^{\zeta} \end{bmatrix}^T  \begin{bmatrix} \tilde{\zeta}_i \\ - \tilde{p}_i^c \end{bmatrix}$ since
$E_s$ only changes the ordering of vectors without affecting their norms.
%\todoiny{Do we have here an inequality or an equality?}
\end{proof}

\subsection{Proof of Theorem~\ref{thm:convergence under delays primal-dual} }\label{appendix proof of thm convergence under delays primal-dual}

\begin{proof}
%\ilr{The proof is based on an invariance principle applied to the time delayed system under consideration. In particular, }
\ilr{The proof starts by establishing the boundedness of trajectories in an invariant set formulated via a Lyapunov functional, and then uses an invariance principle \revise{applied to the time delayed system under consideration} to deduce convergence to the equilibrium point.}\\
\lmm{\ilr{More precisely,}
let $x = (\omega, \eta, p^M, \rho^{p}, p^c, \rho^{\zeta}, \zeta, r^p, r^{\zeta})$ \ic{with} $x_t \in \mathcal{C}$.}
We adopt the Lyapunov functional candidate
\begin{align}\label{eq:V_all}
\begin{array}{rl}
V_{\text{all}} (x_t, t) = V_B +  \sum_{(i,j) \in \tilde{E}} V_S^{ij},
\end{array}
\end{align}
where $V_B$, and $V_S^{ij}$ are defined in \eqref{eq:V_B} and \eqref{eq:V_S^ij}, respectively, and each edge $(i,j)$ is counted once in the summation.
Following results from \ic{the derivations of} \Cref{lem:passivity of reformulated controller,lem:passivity of scattering transformation}, the time derivative of ${V}_{\text{all}} (x_t, t)$ \ic{is given by}
\begin{align*}
\begin{array}{rl}
	\dot{V}_{\text{all}} = \hspace{-2mm} & -\sum_{j \in N} \left( \Lambda_j \left( \omega_j - \omega_j^* \right)^2 \ml{ + (\rho_j^{p})^2 + (\rho_j^{\zeta})^2} \right) \\
	& - \sum_{j \in G} \left(p_j^M - p_j^{M,*}\right) \left( Q_j'\left( p_j^M \right) - Q_j'\left( p_j^{M,*} \right)\right)\\
&
+ \sum_{j\in N}
	\left( \sum_{ i \in \tilde{N}_j} \alpha_{ij}
	\begin{bmatrix}
		\tilde{r}^{p}_{ij} \\ \tilde{r}^{\zeta}_{ij}
	\end{bmatrix} ^T
	\begin{bmatrix}
		\tilde{\zeta}_j \\ - \tilde{p}_j^c
	\end{bmatrix}
	\right) \\
&
- \sum_{(i,j) \in \tilde{E}} \alpha_{ij}
	\left(
	\begin{bmatrix} \tilde{r}_{ij}^{p} \\ \tilde{r}_{ij}^{\zeta} \end{bmatrix}^T \begin{bmatrix} \tilde{\zeta}_j \\ - \tilde{p}_j^c \end{bmatrix} + \begin{bmatrix} \tilde{r}_{ji}^{p} \\ \tilde{r}_{ji}^{\zeta} \end{bmatrix}^T \begin{bmatrix} \tilde{\zeta}_i \\ - \tilde{p}_i^c \end{bmatrix} \right) \\
	= \hspace{-2mm} & -\sum_{j \in N} \left( \Lambda_j \left( \omega_j - \omega_j^* \right)^2 \ml{ + (\rho_j^{p})^2 + (\rho_j^{\zeta})^2} \right)  \\
	& - \sum_{j \in G} \left(p_j^M - p_j^{M,*}\right) \left( Q_j'\left( p_j^M \right) - Q_j'\left( p_j^{M,*} \right)\right)
\end{array}
\end{align*}
where the last equality holds since the last two terms cancel each other,
\begin{align*}
\begin{array}{rl}
\begin{bmatrix}
		\tilde{\zeta} \\ - \tilde{p}^c
	\end{bmatrix}^T
	\begin{bmatrix}
		\tilde{r}^{p} \\ \tilde{r}^{\zeta}
	\end{bmatrix}
=
\sum_{j\in N}
	\left( \sum_{ i \in \tilde{N}_j} \alpha_{ij}
	\begin{bmatrix}
		\tilde{r}^{p}_{ij} \\ \tilde{r}^{\zeta}_{ij}
	\end{bmatrix} ^T
	\begin{bmatrix}
		\tilde{\zeta}_j \\ - \tilde{p}_j^c
	\end{bmatrix}
	\right) \\
=  \sum_{(i,j) \in \tilde{E}} \alpha_{ij}
	\left(
	\begin{bmatrix} \tilde{r}_{ij}^{p} \\ \tilde{r}_{ij}^{\zeta} \end{bmatrix}^T \begin{bmatrix} \tilde{\zeta}_j \\ - \tilde{p}_j^c \end{bmatrix} + \begin{bmatrix} \tilde{r}_{ji}^{p} \\ \tilde{r}_{ji}^{\zeta} \end{bmatrix}^T \begin{bmatrix} \tilde{\zeta}_i \\ - \tilde{p}_i^c \end{bmatrix} \right)
\end{array}
\end{align*}
where $\tilde{r}^{p}$, $\tilde{r}^{\zeta}$ are defined in \Cref{lem:passivity of reformulated controller}.\\
\ml{{\textbf{\ilrr{Boundedness:}}}
Observe that $V_{\text{all}} (x_t)$ is radially unbounded with respect to $\left( \omega, p^c, \zeta, \rho^{p}, \rho^{\zeta}, p^M \right)$.
Define the level set $\Omega = \left\{x_t: V_{\text{all}} (x_t, t) \leq c \right\}$.
%\ilr{We first show that the trajectories in $\Omega$ are bounded for $c$ sufficiently small.}
The integrand term of $V_{P}$ in \eqref{eq:V_P} is zero at $\eta_{ij}^*$ which implies that $V_P$ \icl{has a strict local minimum} at $\eta_{ij}^*$ \ic{from Assumption \ref{assumption angle}}.
%and so does $V_{\text{all}l}$.
\mm{Then, for sufficiently small $c$, $V_{\text{all}} \geq 0$, and states $(\omega, \eta, p^M, \rho^{p}, p^c, \rho^{\zeta}, \zeta)$ are bounded \ic{within $\Omega$.}}}
%Assumption \ref{assumption angle} ensures that $V_{\text{all}}$ is locally strictly %increasing with respect to $\eta_{ij} - \eta_{ij}^*$ and the fact that $\dot{V}_{\text{all}} %\leq 0$ guarantees that $\eta_{ij}$ remains close to $ \eta_{ij}^*$ at all time.
%Then, $\Omega$ is positively invariant for sufficiently small $c$.
%
%
%\todoiny{\st{Analogous comment as in Lemma 2.}}
%
%\ml{Since $V_{\text{all}} (x_t, t)$ is nonincreasing, for any \mm{equally spaced} increasing time sequence $\{t_{k}\}$, we have
%\begin{align*}
%	V_{\text{all}} (x_{t_{k+1}}, t_{k+1}) - V_{\text{all}} (x_{t_k}, t_k) = \int_{t_k}^{t_{k+1}} \dot{V}_{\text{all}} (x(\tau)) d \tau \leq 0.
%\end{align*}
%As $V_{\text{all}}$ is bounded below by zero, $\lim_{t \to \infty} V_{\text{all}} (x_t, t)$ exists and is finite. Then, $\int_{t_k}^{t_{k+1}} \dot{V}_{\text{all}} (x(\tau)) d \tau \to 0$, as $k \to \infty$. Since $\dot{V}_{\text{all}}$ is continuous, \mm{bounded}, and non-positive, we have $\dot{V}_{\text{all}} \to 0$, as $t \to \infty$.
%}
\lmm{Let $\tilde{x} = (\omega, \eta, p^M, \rho^{p}, p^c, \rho^{\zeta}, \zeta)$. By \Cref{rm:eliminating algebraic loop}, given each initial condition, we can find a vector field $f$ such that \ic{the trajectory generated by} the closed-loop system \eqref{eq:power system model}, \eqref{eq:generation dynamics}, \eqref{eq:generation input primal-dual}, \eqref{eq:control with scattering transformation primal-dual}, \eqref{eq: r_ij eliminating algebraic loop} \ic{satisfies a non-autonomous} system \ic{of the form} $\dot{\tilde{x}} = f (t, \tilde{x})$.
%\todoiny{does the algebraic form affect the arguments in Lasalle's paper? I guess not because the variables that appear in the algebraic relation in \Cref{rm:eliminating algebraic loop} have bounded integrals?\\
%\mm{Mengmou: It turns that we do not need the descriptor system to describe it for individual trajectories.} }
}
\mml{\ic{Within $\Omega$ we have from \eqref{eq:V_S^ij}, \eqref{eq:V_all} that $s_{\overrightarrow{ij}}$, $s_{\overrightarrow{ji}}$ are square integrable over any bounded prescribed time interval with these integrals having a uniform bound.}
As $\zeta$, $p^c$ are bounded, we have that $r_{ij}^{p}$, $r_{ij}^{\zeta}$ are \ic{also} square integrable over any bounded \ic{time} interval \ic{with a uniform bound}.
%It follows that they are also absolutely integrable \lmm{with uniform bound} over any finite interval.
\lmm{\ic{From this and the boundedness of $\tilde x$ in $\Omega$
%In addition, all states are within $\Omega$},
we have the property}} $ \| \int_{\alpha}^{\beta} f (\tau, \tilde{x} (\tau)) d \tau \| \leq \mu (\beta - \alpha)$, for some $\mu > 0$.
\lmm{%This means that states must travel finite distance in any finite time interval.
%Following similar arguments from
\ic{Then from the proof of} \cite[Theorem 1]{lasalle1976stability},\footnote{\ic{It should be noted that proof of \cite[Theorem 1]{lasalle1976stability} involves the analysis of individual \ill{trajectories,
and relies on the integral inequality stated in the previous sentence.}}}}
%and relies on the fact that due to the integral inequality stated in the previous sentence, states travel only bounded distances in finite time. }}}
%The arguments are based on the study of one single trajectory while the same result holds for %arbitrary trajectories generated by different initial conditions.}}
we obtain that $\dot{V}_{\text{all}} \to 0$, as $t \to \infty$, which implies $\rho_j^{\zeta}$, $\rho_j^{p} \to 0$. From \eqref{eq:control with scattering transformation primal-dual}, \ic{for any} increasing sequence $\{t_{k}\}$, \ic{with $t_k\to\infty$ as $k\to\infty$, we have}
%\begin{align*}
%\begin{array}{rl}
$
	 \underset{k \to \infty}{\lim} \left( \zeta_j (t_{k+1}) - \zeta_j (t_k) \right) =  \underset{k \to \infty}{\lim} \int_{t_k}^{t_{k+1}} \dot{\zeta}_j (\tau) d\tau
	  = 2 \underset{k \to \infty}{\lim} \int_{t_k}^{t_{k+1}} \dot{\rho}_j^{\zeta} (\tau) d\tau = 2 \underset{k \to \infty}{\lim} \left( \rho_j^{\zeta} (t_{k+1}) - \rho_j^{\zeta} (t_k) \right)= 0.
$
%\end{array}
%\end{align*}
Similar argument holds for $p_j^c$. Thus, $\zeta_j$, $p_j^c$ also have finite limits.
}
\mml{
%Since $r_{ij}^{p} (t)$, $r_{ij}^{\zeta} (t)$ are continuous by \eqref{eq: r_ij eliminating algebraic loop} and square integrable %by \eqref{eq:V_S^ij}, they are finite over any bounded interval.
\ic{We also} obtain from \eqref{eq: r_ij eliminating algebraic loop} that $r_{ij}^{p} (t) + r_{ij}^{p} (t - T_{ij} - T_{ji})$, $r_{ij}^{\zeta} (t) + r_{ij}^{\zeta} (t - T_{ij} - T_{ji})$ tend to constants \ic{as $t\to\infty$}. In other words, $r_{ij}^{p} (t)$, $r_{ij}^{\zeta} (t)$ tend to periodic functions with period $T = 2T_{ij} + 2T_{ji}$. Thus, we conclude that $r_{ij}^{p} (t)$, $r_{ij}^{\zeta} (t)$ are bounded for all $t$.}\\
%\mm{Now that we have all system \ic{states} %variables
%are bounded \mml{in the set $\Omega$,}
\textbf{\ilrr{Invariance principle:}}
\ic{The invariance principle in the proof of \cite[Theorem 3.1]{hale2013introduction} is applicable.} \lmm{\ic{In particular, since all trajectories in $\Omega$ are bounded and their $\omega$-limit set is an \ic{invariant set  \cite{hale2013introduction},}
%(due to the continuity of the solutions w.r.t. initial conditions),
all trajectories starting in $\Omega$ will converge to the}} \ml{largest invariant set \ic{in} $\{x_t | \dot{V}_{\text{all}} \equiv 0\}$.  In this invariant set, $\omega, p^c, \zeta, \rho^{p}, \rho^{\zeta}, p^M$ are \ic{constant},  \ic{therefore}, $\dot{p}_j^{c}, \dot{\zeta}_j = 0$}. \ic{This implies that for all trajectories in $\Omega$} $r_{ij}^p$, $r_{ij}^{\zeta}$ tend to constants, \ic{as follows from} \eqref{eq:control with scattering transformation primal-dual} and the fact that \ic{for a connected graph the corresponding adjacency matrix is full-rank.}
\mm{The scattering variables in \eqref{eq:scattering transformation primal-dual} converge \ic{to constant values} ${s}_{\overrightarrow{ij}}^*$, ${s}_{\overrightarrow{ji}}^*$, ${s}_{\overleftarrow{ij}}^*$, ${s}_{\overleftarrow{ji}}^*$ \ic{and we also have the corresponding constant vectors} $\begin{bmatrix} {r}_{ij}^{p,*} \\ {r}_{ij}^{\zeta,*} \end{bmatrix}$ and $\begin{bmatrix} {r}_{ji}^{p,*} \\ {r}_{ji}^{\zeta,*} \end{bmatrix}$. We have from \eqref{eq:scattering variables under delays primal-dual} that
$
	{s}_{\overleftarrow{ji}}^*(t) = E_s {s}_{\overrightarrow{ij}}^*(t),~
	 {s}_{\overleftarrow{ij}}^*(t) = E_s {s}_{\overrightarrow{ji}}^*(t).
$
By \eqref{eq:scattering transformation primal-dual}, we have
\begin{align*}
\begin{array}{rl}
	\begin{bmatrix} {r}_{ij}^{p,*} \\ {r}_{ij}^{\zeta,*} \end{bmatrix} = \frac{1}{\sqrt{2}}\left( {s}_{\overrightarrow{ji}}^* + {s}_{\overleftarrow{ji}}^* \right) =\frac{1}{\sqrt{2}} \left( E_s^{-1} {s}_{\overleftarrow{ij}}^* + E_s {s}_{\overrightarrow{ij}}^* \right) = \begin{bmatrix} {p}_i^{c,*} \\ {\zeta}_i^* \end{bmatrix}\\
	\begin{bmatrix} {r}_{ji}^{p,*} \\ {r}_{ji}^{\zeta,*} \end{bmatrix} = \frac{-1}{\sqrt{2}}\left( {s}_{\overrightarrow{ij}}^* + {s}_{\overleftarrow{ij}}^* \right) =\frac{-1}{\sqrt{2}} \left( E_s^{-1} {s}_{\overleftarrow{ji}}^* + E_s {s}_{\overrightarrow{ji}}^* \right) = \begin{bmatrix} {p}_j^{c,*} \\ {\zeta}_j^* \end{bmatrix}
\end{array}
\end{align*}
which recovers the \ic{relation that holds} %interconnection
in the undelayed case \ic{at equilibrium}. In addition, \ic{since also} $\rho_j^{\zeta}, \rho_j^{p} \rightarrow 0$, the optimality of the equilibrium point is guaranteed from Corollary~\ref{cor:convergence primal-dual}.
}

In conclusion, the solutions of \eqref{eq:power system model}, \eqref{eq:generation dynamics}, \eqref{eq:generation input primal-dual}, \eqref{eq:control with scattering transformation primal-dual}, \eqref{eq:scattering transformation primal-dual}, \eqref{eq:scattering variables under delays primal-dual} with initial conditions in the neighborhood of \il{the equilibrium point considered} %$\eta^*$
will converge to \il{a set of equilibrium points that solve} %the optimal solution of
the \textup{OGR} problem \eqref{optimal generation regulation problem} with $\omega^* = \mathbf{0}_{|N|}$.
\il{Convergence to a single equilibrium point follows using arguments analogous to those in \cite[Prop. 4.7]{haddad2011nonlinear}, by noting that each equilibrium point in the set where trajectories converge is also Lyapunov stable.}
\end{proof}

\subsection{Proof of Lemma~\ref{lem:optimality tie-line primal-dual} }\label{appendix proof of lemma optimality}
\begin{proof}
	Problem \eqref{optimal generation regulation problem 2} is solved if the following problem is solved,
\begin{align}\label{optimal generation regulation problem 2 equivalence}
\begin{array}{rl}
\hspace{-4mm} \textbf{OGR-4:} ~ \underset{p^M, p}{\min} &  \displaystyle \sum_{j \in G} Q_j (p_j^M),\\
 \hspace{-10mm} \text{subject to~} &  p_j^M - p_j^L + \displaystyle \hspace{-1mm} \sum_{(i,j) \in E } {D}_{ij} p_{ij} = 0, ~ j \in N, \\
	& \displaystyle \sum_{(i,j) \in {B}_{k}} \hspace{-1mm} \hat{D}_{k,ij} p_{ij} = \hat{P}_{k}, ~ k \in {K},
\end{array}
\end{align}
where $p_j^M : = 0$ for all $j \notin G$, since summing up the first class of constraints for all $j \in N$ gives \eqref{optimal generation regulation problem 2}.
The KKT conditions for \eqref{optimal generation regulation problem 2 equivalence} are
%\todoing{is the summation in the second equation over $E$ or over $\tilde{E}$?}
\begin{subequations}\label{eq:KKT condition tie-line}
	\begin{align}
		Q_j'\left(\bar{p}_j^{M}\right) = \beta_j, ~ j \in G \label{eq:KKT condition tie-line 1}\\
		{D}^T \left( E_{{K}}^T \nu - \beta \right) = \mathbf{0} \label{eq:KKT condition tie-line 2}\\
		\bar{p}^{M} - p^L + {D} \bar{p} = \mathbf{0} \label{eq:KKT condition tie-line 3}\\
		E_{{K}} {D} \bar{p} = \hat{P} \label{eq:KKT condition tie-line 4}
	\end{align}
\end{subequations}
for some constant vector $\beta = (\beta_1, \ldots,\beta_{|N|}) \in \mathbb{R}^{|N|}$ and $ \nu \in \mathbb{R}^{|{K}|}$, where \eqref{eq:KKT condition tie-line 4} follows from \eqref{eq:D=E_K D}.
By similar arguments from Lemma~\ref{lem:optimality without delays primal-dual}, we have that $\omega^* \in \text{Im}(\mathbf{1}_{|N|})$.
The equilibrium of \eqref{eq:equilibrium of generation dynamics}, \eqref{eq:generation input primal-dual} gives $Q_j' \left( {p}_j^{M,*} \right) = {p}_j^{c,*}$, which is equivalent to \eqref{eq:KKT condition tie-line 1} by setting $\beta_j = {p}_j^{c,*}$.
Since it is assumed that each component in the subgraph $\mathcal{G}(N,\tilde{E} / \tilde{B})$ is connected, we have that $E_{{K}} \tilde{L}_{{K}} = \mathbf{0}$, which means that the symmetric matrix $\tilde{L}_{{K}}$ has the null space $\text{Im} \left( E_{{K}}^T \right)$. Then, the equilibrium of \eqref{eq:controller without delays tie-line compact form primal-dual 4} implies that $\pi^* =  E_{{K}}^T{\nu^*}$ for some $\nu^* \in \mathbb{R}^{|{K}|}$.
Let $\nu = \nu^*$ and $\beta_j = {p}_j^{c,*}$, the equilibrium of \eqref{eq:controller without delays tie-line compact form primal-dual 1} implies that ${\pi^*}  - {p}^{c,*} = E_{{K}}^T {\nu} - \beta \in \text{Im}\left( \mathbf{1}_{|N|} \right)$, satisfying \eqref{eq:KKT condition tie-line 2} since the communication graph is connected.
Next, comparing the equilibrium of \eqref{eq:controller without delays tie-line compact form primal-dual 2} with \eqref{eq:equilibrium of swing equations} and \eqref{eq:equilibrium of load balance}, we have that ${D} p^* = - \tilde{L} \zeta^*$. Thus, the equilibrium of \eqref{eq:controller without delays tie-line compact form primal-dual 2} gives \eqref{eq:KKT condition tie-line 3}.
Finally, we can observe that $E_{{K}} J = I_{|{K}|}$, left multiplying the equilibrium of \eqref{eq:controller without delays tie-line compact form primal-dual 3} by $E_{{K}}$, we obtain \eqref{eq:KKT condition tie-line 4}.
\end{proof}

\subsection{Proof of Lemma~\ref{lem:passivity of controller tie-line primal-dual} }\label{appendix proof of lem passivity of controller tie-line primal-dual}
\begin{proof}
	Adopt the storage functional
\begin{align}\label{eq:V_T}
\begin{array}{rl}
	V_T = V_B + V_{\pi, \rho} (\pi, \rho^{\pi}) + V_{\phi,\rho}(\phi, \rho^{\phi})
\end{array}
\end{align}
where $V_B$ is defined in \eqref{eq:V_B}, and
\begin{align*}
\begin{array}{rl}
	V_{\pi, \rho} (\pi, \rho^{\pi}) =\displaystyle \frac{1}{2}  \sum_{j \in N} (\rho_j^{\pi})^2  + \frac{1}{2} \sum_{j \in N} \left( \pi_j - \rho_j^{\pi} - \pi_j^* \right)^2, \\
	V_{\phi,\rho} (\phi, \rho^{\phi})= \displaystyle \frac{1}{2}  \sum_{j \in N} (\rho_j^{\phi})^2 + \frac{1}{2} \sum_{j \in N} \left( \phi_j - \rho_j^{\phi} - \phi_j^*\right)^2.
\end{array}
\end{align*}
Then, $V_T \geq 0$ in some neighborhood of \il{the equilibrium point} %$\eta^*$
\ic{as follows from}
Assumption \ref{assumption angle}. The time derivative of $V_T$ along the system trajectories satisfies
\begin{align*}
\begin{array}{rl}
	\hspace{-2mm} \dot{V}_T \hspace{-1mm} \leq
%	\hspace{-3.5mm}
	& \displaystyle  -  \hspace{-0.5mm} \sum_{j \in N} (p_j^c - p_j^{c,*})  \hspace{-0.5mm} \sum_{i \in \tilde{N}_j}  \hspace{-0.5mm} \alpha_{ij} \hspace{-0.5mm} \left( \left(r_{ij}^{\zeta} \hspace{-0.5mm} - \hspace{-0.5mm} r_{ij}^{\zeta,*} \right) \hspace{-1mm} - \hspace{-1mm} \left(\zeta_j - \zeta_j^* \right) \right) \\
	& \displaystyle +  \hspace{-0.5mm}  \sum_{j \in N} (\zeta_j - \zeta_j^*) \sum_{i \in \tilde{N}_j}  \hspace{-0.5mm} \alpha_{ij} \left( \left(r_{ij}^{p} \hspace{-0.5mm} - \hspace{-0.5mm} r_{ij}^{p,*}\right) \hspace{-1mm}-\hspace{-1mm} \left(p_j^c \hspace{-0.5mm} - \hspace{-0.5mm} p_j^{c,*}\right) \right)\\
	& \displaystyle -  \hspace{-0.5mm}  \sum_{j \in N} (\zeta_j - \zeta_j^*) \hspace{-0.5mm} \sum_{i \in \tilde{N}_j} \hspace{-0.5mm} \alpha_{ij} \hspace{-0.5mm} \left( \left(r_{ij}^{\pi'} \hspace{-0.5mm} - \hspace{-0.5mm} r_{ij}^{\pi,*}\right) \hspace{-1mm} - \hspace{-1mm} \left(\pi_j  \hspace{-0.5mm}-  \hspace{-0.5mm} \pi_j^{*}\right) \right)\\
	& \displaystyle +   \hspace{-0.5mm} \sum_{j \in N} (\pi_j \hspace{-0.5mm} - \hspace{-0.5mm} \pi_j^*)  \hspace{-0.5mm} \sum_{i \in \tilde{N}_j} \alpha_{ij} \hspace{-0.5mm} \left( \left(r_{ij}^{\zeta'} - r_{ij}^{\zeta,*}\right) - \left(\zeta_j \hspace{-0.5mm} -\hspace{-0.5mm}  \zeta_j^{*}\right)\right)\\
	& \displaystyle -  \hspace{-0.5mm}  \sum_{j \in N} (\pi_j \hspace{-0.5mm} - \hspace{-0.5mm} \pi_j^*) \hspace{-0.5mm}  \sum_{i \in \tilde{N}_j \cap {C}_{k}} \hspace{-1mm} \alpha_{ij} \hspace{-0.5mm} \left( \left(r_{ij}^{\phi''} \hspace{-1mm} - \hspace{-0.8mm} r_{ij}^{\phi,*}\right) \hspace{-1mm}  - \hspace{-1mm}  \left(\phi_j   \hspace{-0.5mm} -  \hspace{-0.5mm} \phi_j^{*}\right) \right)\\
	& \displaystyle +  \hspace{-0.5mm} \sum_{j \in N}  (\phi_j \hspace{-0.5mm} - \hspace{-0.5mm} \phi_j^*)  \hspace{-0.5mm}  \sum_{i \in \tilde{N}_j \cap {C}_{k}} \hspace{-1mm} \alpha_{ij} \hspace{-0.5mm} \left( \left(r_{ij}^{\pi''} \hspace{-1.5mm} - \hspace{-0.8mm}  r_{ij}^{\pi,*}\right) \hspace{-1mm}  - \hspace{-1mm}  \left(\pi_j \hspace{-1mm}  - \hspace{-1mm}  \pi_j^{*}\right) \right)\\
\leq
%\hspace{-2mm}
& \begin{bmatrix}
		\tilde{r}^{p} ~ \tilde{r}^{\zeta} ~ \tilde{r}^{\zeta'} ~ \tilde{r}^{\pi'} ~ \tilde{r}^{\pi''} ~ \tilde{r}^{\phi''}
		\end{bmatrix} \begin{bmatrix}
		\tilde{\zeta} ~ -\tilde{p}^c ~ \tilde{\pi} ~ - \tilde{\zeta} ~ \tilde{\phi} ~ - \tilde{\pi}
	\end{bmatrix}^T.
\end{array}
\end{align*}
\end{proof}
%\todoiny{ do not fully understand the last sentence}
\subsection{Proof of Lemma~\ref{lem:passivity of scattering transformation tie-line} }{\label{appendix proof of lem passivity of scattering transformation tie-line}}
\begin{proof}
	\ic{We adopt a storage functional $V_S^{ij}$ which is of the same form as that in \eqref{eq:V_S^ij}, where $s_{\overrightarrow{ij}}^*$, $s_{\overrightarrow{ji}}^*$, $s_{\overleftarrow{ij}}^*$ and $s_{\overleftarrow{ji}}^*$ are values at the equilibrium point of \eqref{eq:scattering transformation tie-line primal-dual}. Then, the Lemma follows in an analogous way to that of Lemma \ref{lem:passivity of scattering transformation}.}
%similar results can be obtained.
\end{proof}

\subsection{Proof of Theorem~\ref{thm:Convergence under delays tie-line primal-dual} }\label{appendix proof of theorem convergence under delays tie-line primal-dual}
\begin{proof}
We adopt the Lyapunov functional candidate $V_{\text{tie}} = V_T + \sum_{(i,j) \in \tilde{E}} V_S^{ij}$, where $V_T$ and $V_S^{ij}$ are defined in \eqref{eq:V_T} and \eqref{eq:V_S^ij} and each edge $(i,j)\in \tilde{E}$ is counted once. Then,
%\todoing{is the first inequality below an equality? \\ Yes, I think so.}
\begin{align*}
\begin{array}{rl}
	\hspace{-2mm}  \dot{V}_{\text{tie}}
	= \hspace{-3mm} & - \hspace{-1mm} \displaystyle \sum_{j \in N} \hspace{-1mm} \left( \Lambda_j \hspace{-0.5mm} \left( \omega_j \hspace{-0.5mm} - \hspace{-0.5mm} \omega_j^* \right)^2 \hspace{-0.5mm} + \hspace{-0.5mm} (\rho_j^{p})^2 \hspace{-0.5mm} + \hspace{-0.5mm}(\rho_j^{\zeta})^2 \hspace{-0.5mm} + \hspace{-0.5mm} (\rho_j^{\pi})^2 \hspace{-0.5mm} + \hspace{-0.5mm} (\rho_j^{\phi})^2 \right) \\
	& \displaystyle - \sum_{j \in G} \left(p_j^M - p_j^{M,*}\right) \left( Q_j'\left( p_j^M \right) - Q_j'\left( p_j^{M,*} \right)\right)\\
& \displaystyle +		\sum_{j \in N} \sum_{i \in \tilde{N}_j} \alpha_{ij}
	\begin{bmatrix}
		\tilde{r}_{ij}^{p} \hspace{-2mm} & \tilde{r}_{ij}^{\zeta} \hspace{-2mm} & \tilde{r}_{ij}^{\zeta'} \hspace{-2mm} & \tilde{r}_{ij}^{\pi'}
	\end{bmatrix}
	\begin{bmatrix}
		\tilde{\zeta}_j \hspace{-2mm} & -\tilde{p}_j^c \hspace{-2mm} & \tilde{\pi}_j \hspace{-2mm} & -\tilde{\zeta}_j
	\end{bmatrix}^T \\
& \displaystyle - \sum_{(i,j) \in \tilde{E}} \alpha_{ij}\left( \begin{bmatrix}
		\tilde{r}^{p}_{ji} \hspace{-2mm} & \tilde{r}^{\zeta}_{ji} \hspace{-2mm} & \tilde{r}^{\zeta'}_{ji} \hspace{-2mm} & \tilde{r}^{\pi'}_{ji}
	\end{bmatrix}
	\begin{bmatrix}
		\tilde{\zeta}_i \hspace{-2mm} & -\tilde{p}_i^c \hspace{-2mm} & \tilde{\pi}_i \hspace{-2mm} & -\tilde{\zeta}_i
	\end{bmatrix}^T
	 \right.\\
&	\left. \qquad \qquad \displaystyle
	+ \begin{bmatrix}
		\tilde{r}^{p}_{ij} \hspace{-2mm} & \tilde{r}^{\zeta}_{ij} \hspace{-2mm} & \tilde{r}^{\zeta'}_{ij} \hspace{-2mm} & \tilde{r}^{\pi'}_{ij}
	\end{bmatrix}
	\begin{bmatrix}
		\tilde{\zeta}_j \hspace{-2mm} & -\tilde{p}_j^c \hspace{-2mm} & \tilde{\pi}_j \hspace{-2mm} & -\tilde{\zeta}_j
	\end{bmatrix}^T \right)\\
&\displaystyle  +	\sum_{j \in N} \sum_{i \in \tilde{N}_j \cap C_{k}} \alpha_{ij} \begin{bmatrix}
		\tilde{r}^{\pi''}_{ji} & \tilde{r}^{\phi''}_{ji} \end{bmatrix}
	\begin{bmatrix}
		\tilde{\phi}_j \\ -\tilde{\pi}_j^c
	\end{bmatrix} \\
& \displaystyle - \hspace{-3mm} \sum_{(i,j) \in \tilde{E} / \tilde{B} } \hspace{-2mm} \alpha_{ij}\left( \begin{bmatrix}
		\tilde{r}^{\pi''}_{ji} & \hspace{-1mm} \tilde{r}^{\phi''}_{ji} \end{bmatrix}
	\begin{bmatrix}
		\tilde{\phi}_i \\ -\tilde{\pi}_i^c
	\end{bmatrix}
	+
	\begin{bmatrix}
		\tilde{r}^{\pi''}_{ij} & \hspace{-1mm} \tilde{r}^{\phi''}_{ij} \end{bmatrix}
	\begin{bmatrix}
		\tilde{\phi}_j \\ -\tilde{\pi}_j^c
	\end{bmatrix} \right)\\
\leq \hspace{-3mm} & \displaystyle - \sum_{j \in N} \hspace{-1mm} \left( \Lambda_j \hspace{-0.5mm} \left( \omega_j \hspace{-0.5mm} - \hspace{-0.5mm} \omega_j^* \right)^2 \hspace{-0.5mm} + \hspace{-0.5mm} (\rho_j^{p})^2 \hspace{-0.5mm} + \hspace{-0.5mm}(\rho_j^{\zeta})^2 \hspace{-0.5mm} + \hspace{-0.5mm} (\rho_j^{\pi})^2 \hspace{-0.5mm} + \hspace{-0.5mm} (\rho_j^{\phi})^2 \right)\\
& \displaystyle - \sum_{j \in G} \left(p_j^M - p_j^{M,*}\right) \left( Q_j'\left( p_j^M \right) - Q_j'\left( p_j^{M,*} \right)\right)
\end{array}
\end{align*}
where the second inequality holds since the sum of the last four terms is zero.
Using arguments
%of the Barbalat's Lemma
\il{analogous to those in the proof of \Cref{thm:convergence under delays primal-dual}, trajectories with initial conditions sufficiently close to the equilibrium point considered will converge to an equilibrium point where $\rho_j^{\zeta}, \rho_j^{\pi}, \rho_j^{\phi} = 0$, for all $j \in N$.}
%can obtain that variables $(\omega, \eta, p^M,\zeta,\pi,\phi)$ converge to constants and $\rho_j^{p},
%\rho_j^{\zeta}, \rho_j^{\pi}, \rho_j^{\phi} \rightarrow 0$, for all $j \in N$.
%As a result, the scattering variables converge to some equilibrium ${s}_{\overrightarrow{ij}}^*$, ${s}_{\overrightarrow{ji}}^*$, ${s}_{\overleftarrow{ij}}^*$, ${s}_{\overleftarrow{ji}}^*$, we have that \begin{align*}
%	&{s}_{\overleftarrow{ji}}^*(t) = E_s {s}_{\overrightarrow{ij}}^*(t - T_{ij}) = E_s {s}_{\overrightarrow{ij}}^*(t)\\
%	 &{s}_{\overleftarrow{ij}}^*(t) = E_s {s}_{\overrightarrow{ji}}^*(t - T_{ji}) = E_s {s}_{\overrightarrow{ji}}^*(t)
%\end{align*}
%and the associated variables satisfy that
%\begin{align*}
%	{r}_{ij}^{p,*} = {p}_j^{c,*}, ~{r}_{ij}^{\pi,*} = {\pi}_j^{*}, ~ {r}_{ij}^{\zeta,*} = {\zeta}_j^{*}, ~ {r}_{ij}^{\phi,*} = {\phi}_j^{*}, ~ \forall j \in N
%\end{align*}
\ml{As a result, %\mm{an}
%equilibrium
\il{an equilibrium point} of \eqref{eq:controller without delays tie-line primal-dual} for the undelayed case \il{is} %is
recovered,
%
%\todoiny{the previous sentence needs to be made more accurate, as there an equilibrium point mentioned in the theorem statement; also you mention "the equilibrium point" - is this unique?}
which guarantees \il{that this solves \textup{OGR-2}} from \Cref{lem:optimality tie-line primal-dual}. }
%Therefore,
%%\il{using also the Lyapunov stability of the equilibrium points as in }
%system states with initial condition in the neighborhood of \il{the equilibrium point considered} %$\eta^*$
%\il{will converge to an equilibrium point that solves} %an equilibrium that solves
%the \textup{OGR-2} problem \eqref{optimal generation regulation problem 2}.
%\todoiny{is the equilibrium point unique? Do we have convergence to an equilibrium point or a set of equilibrium points?}
\end{proof}

\subsection{Proof of Lemma~\ref{lem:optimality generation bounds} }\label{appendix proof of lemma optimality generation bounds}
\begin{proof}
The equilibrium of \eqref{eq:controller without delays primal-dual}, \eqref{eq:generation input generation bounds}, \eqref{eq:local inequality multipliers}, satisfies
\begin{subequations}
	\begin{align}
		& 0 = p_i^{c,*} - p_j^{c,*}, (i,j) \in \tilde{E}, \label{eq:equilibrium of controller generation bounds 1}\\
		& 0 = p_j^L - p_j^{M,*} - \sum_{k: j \rightarrow k} \psi_{jk}^* + \sum_{i: i \rightarrow j} \psi_{ij}^*, ~ j \in N, \label{eq:equilibrium of controller generation bounds 2}\\
		& p_j^{c,*} - \omega_j^* =  Q_j'\left( p_j^{M,*} \right) - {\lambda_j^*}^2 + {\mu_j^*}^2, ~ j \in G, \label{eq:equilibrium of generation input generation bounds}\\
		& 0  \hspace{-0mm} =  \hspace{-0mm} 2 \lambda_j^*  \hspace{-0mm} \left( p_j^{M,\min}  \hspace{-0mm} - \hspace{-0mm} p_j^{M,*} \right), \\
		& 0  \hspace{-0mm} =  \hspace{-0mm} 2 \mu_j^*  \hspace{-0mm} \left( p_j^{M,*}  \hspace{-0mm} -  \hspace{-0mm} p_j^{M, max} \right),
	\end{align}
\end{subequations}
where \eqref{eq:equilibrium of generation input generation bounds} is obtained by combining the equilibrium of \eqref{eq:generation input generation bounds} and \eqref{eq:equilibrium of generation dynamics}.
The KKT conditions for problem \eqref{optimal generation regulation problem 3} are
\begin{subequations}\label{eq:KKT condition generation bounds}
%\todoiny{The KKT conditions seem to be different from the ones previously used for the case without generation constraints.\lmm{Yes, they are slightly different due to different Lagrangians.}}
	\begin{align}
		Q_j'\left(\bar{p}_j^{M}\right) - \bar{\lambda}_j + \bar{\mu}_j = \beta , ~ \forall j \in G, \label{eq:KKT condition generation bounds 1}\\
		 \sum_{j \in G} \bar{p}_j^{M} = \sum_{j \in N} p_j^L,\label{eq:KKT condition generation bounds 2}\\
	 p_j^{M,\min} \leq \bar{p}_j^{M} \leq p_j^{M,\max}, \forall j \in G,\label{eq:KKT condition generation bounds 3}\\
	 \bar{\lambda}_j \left( p_j^{M,\min} - \bar{p}_j^{M} \right) = 0, \\
	 \bar{\mu}_j \left( \bar{p}_j^{M}- p_j^{M, \max} \right) = 0,
	\end{align}
\end{subequations}
for some $\bar{\lambda}_j$,  $\bar{\mu}_j \geq 0$, $\forall j \in G$, and constant $\beta$. The equilibrium equations \eqref{eq:equilibrium of omega} implies $\omega^* \in \text{Im}(\mathbf{1}_{|N|})$. Summing \eqref{eq:equilibrium of controller generation bounds 2} for all $j \in N$, we obtain \eqref{eq:KKT condition generation bounds 2}, and $\omega^* = \mathbf{0}_{|N|}$. The equations \eqref{eq:equilibrium of controller generation bounds 1} implies that $p^c \in \text{Im} (\mathbf{1}_{|N|})$.
\ml{Moreover, by solving the ordinary differential equations in \eqref{eq:local inequality multipliers}, we can obtain that $\lambda_j$, $\mu_j \geq 0$ given positive initial conditions, and $\lambda_j^* = 0$,  $\mu_j^* = 0$ only if the inequalities for generation bounds are satisfied \ill{at equilibrium}.}
Then, \ill{letting} \ml{$\beta = p_j^{c,*}$, $\bar{\lambda}_{j} = {\lambda_j^*}^2$, $\bar{\mu}_j = {\mu_j^*}^2$},  the rest of the KKT conditions are satisfied.
\end{proof}

\subsection{Proof of Theorem~\ref{thm:convergence under delays generation bounds} }\label{appendix proof of theorem convergence under delays generation bounds}
\begin{proof}
Consider the Lyapunov functional candidate $V_{\text{gb}} = V_{\text{all}} + V_{G}(\lambda, \mu) $, where $V_{\text{all}}$ is defined in \eqref{eq:V_all}, and
\begin{align}\label{eq:V_G}
\begin{array}{rl}
\hspace{-5mm} V_G (\lambda, \mu) \hspace{-1mm} = \hspace{-1mm} \sum_{j \in G} \frac{1}{4} \left( \lambda_j^2 - {\lambda_j^{*}}^2 \right) - \frac{1}{2} {\lambda_j^*}^2 \left( \ln{\lambda_j} - \ln{\lambda_j^*}  \right)\\
+ \frac{1}{4} \left(\mu_j^2 - {\mu_j^{*}}^2 \right) - \frac{1}{2} {\mu_j^*}^2 \left( \ln{\mu_j} - \ln{\mu_j^*}  \right)
\end{array}
\end{align}
where $x \ln{x} := 0$ if $x = 0$.
Since $\ln{x} \leq \frac{x}{y} + \ln{y} -1$ for any $x$, $y > 0$, we have
%\begin{align*}
%\begin{array}{rl}
$
V_G \geq \sum_{j \in G}\frac{1}{4} (\lambda_j - \lambda_j^*)^2 + \frac{1}{4} (\mu_j - \mu_j^*)^2 \geq 0.
$
%\end{array}
%\end{align*}
Following results from \Cref{thm:convergence under delays primal-dual}, the time derivative of ${V}_{\text{gb}}$ along the system trajectories gives
\begin{align*}
\begin{array}{rl}
\dot{V}_{\text{gb}} = \hspace{-3mm}
& -\sum_{j \in N} \left( \Lambda_j \left( \omega_j - \omega_j^* \right)^2 + (\rho_j^{p})^2 + (\rho_j^{\zeta})^2 \right)  \\
	& -  \sum_{j \in G} \left(p_j^M - p_j^{M,*}\right) \left( Q_j'\left( p_j^M \right) - Q_j'\left( p_j^{M,*} \right)\right) \\
	& + \sum_{j \in G}  \left( \lambda_j^2 - {\lambda_j^*}^2  \right)  \left(p_j^{M,\min} - p_j^{M,*}\right) \\
	& +  \sum_{j \in G}  \left(\nu_j^2 - {\nu_j^*}^2 \right) \left(p_j^{M,*} - p_j^{M,\max}\right)\\
	\leq  \hspace{-3mm}
	& -\sum_{j \in N} \left( \Lambda_j \left( \omega_j - \omega_j^* \right)^2 + (\rho_j^{p})^2 + (\rho_j^{\zeta})^2 \right)  \\
	& -   \sum_{j \in G} \left(p_j^M - p_j^{M,*}\right) \left( Q_j'\left( p_j^M \right) - Q_j'\left( p_j^{M,*} \right)\right)
\end{array}
\end{align*}
\ml{where the inequality follows from the saddle point property $\mathcal{L} (p^{M,*}, \lambda, \nu)  - \mathcal{L} (p^{M,*}, \lambda^*, \nu^*)  \leq 0$ of the Lagrangian in \eqref{eq:generalized Lagrangian}.}
The rest of the \il{proof is analogous to that}  %arguments are similar to the proof
of \Cref{thm:convergence under delays primal-dual}.
%\todoiny{needs checking as the proof of \Cref{thm:convergence under delays primal-dual} has changed.\\
%\lmm{Yes, the result is unchanged.}}
\end{proof}

\subsection{Proof of Theorem~\ref{thm:convergence observer-based primal-dual} }\label{appendix proof of theorem convergence observer-based primal-dual}
\begin{proof}
Adopt the Lyapunov functional candidate
	\begin{align}
	\begin{array}{rl}
		V_O = V_B + \sum_{(i,j) \in \tilde{E}}V_S^{ij} + V_E(b, \chi, \omega)
	\end{array}
	\end{align}
where $V_B$ and $V_S^{ij}$ are defined in \eqref{eq:V_B}, \eqref{eq:V_S^ij}, and $V_E(b, \chi, \omega)$ is given by
\begin{align}
\begin{array}{rl}
	\hspace{-3mm} V_E(b, \chi, \omega) = & \hspace{-3mm} \frac{1}{2} \sum_{j \in G} \left( M_j \left( ( b_j - b_j^*) - (\omega_j - \omega_j^*) \right)^2 \right. \\
	& \left. \qquad \qquad + \tau_{\chi,j} \left( \chi_j - \chi_j^* \right)^2 \right)
\end{array}
\end{align}
where $\chi_j^*$, $b_j^*$ are the equilibrium of $\chi_j$, $b_j$, respectively. The time derivative of $V_O$ along system \eqref{eq:power system model}, \eqref{eq:generation dynamics}, \eqref{eq:generation input primal-dual}, \eqref{eq:controller observer-based primal-dual}, \eqref{eq:scattering transformation primal-dual}, \eqref{eq:scattering variables under delays primal-dual} gives
\begin{align*}
\begin{array}{rl}
	\hspace{-2mm} \dot{V}_O \leq & \hspace{-3mm} - \sum_{j \in N}  \left( \Lambda_j  \left( \omega_j  - \omega_j^* \right)^2 \hspace{-1mm}  - \hspace{-1mm}  \left(\chi_j \hspace{-0.5mm}  - \hspace{-0.5mm}  \chi_j^*\right)^2 - (\rho_j^{p})^2 - (\rho_j^{\zeta})^2 \right) \\
	& \hspace{-3mm} - \sum_{j \in G} \left(p_j^M - p_j^{M,*}\right) \left( Q_j'\left( p_j^M \right) - Q_j'\left( p_j^{M,*} \right)\right).
\end{array}
\end{align*}
%Using arguments
%%of Barbalat's Lemma
%similarly to the proof of \Cref{thm:convergence under delays primal-dual}, we can prove the convergence of the trajectories to the optimal solution of the \textup{OGR} problem \eqref{optimal generation regulation problem}.
%\todoiny{as above, needs checking as the proof of \Cref{thm:convergence under delays primal-dual} has changed.\\
%\lmm{Yes, the result remains unchanged.}}
\il{The rest of the \il{proof is analogous to that}  %arguments are similar to the proof
of \Cref{thm:convergence under delays primal-dual}.}
\end{proof}

\ill{
\section*{Acknowledgement}
\revise{The authors would like to thank %the \icr{Editor} and
\icr{the Reviewers} for \ilr{their valuable comments.}}
%, which significantly improved the paper.}

This work was supported by ERC starting grant 679774. 
The first author is partly supported by National Natural Science Foundation of China under grant number 62173155.
For the purpose of open access, the authors have applied a Creative Commons Attribution (CC BY) licence  to any Author Accepted Manuscript version arising.}

%\section*{Acknowledgment}
%The authors would like to thank the Editors and the Reviewers
%for their time.

\ifCLASSOPTIONcaptionsoff
  \newpage
\fi



%\bibliographystyle{IEEEtran}

\bibliographystyle{IEEEtran}
\bibliography{References}

%\begin{biography}[{\includegraphics[width=1in,height=1.25in,clip,keepaspectratio]{mshell}}]{Michael Shell}

\begin{IEEEbiography}[{\includegraphics[width=1in,height=1.25in,clip,keepaspectratio]{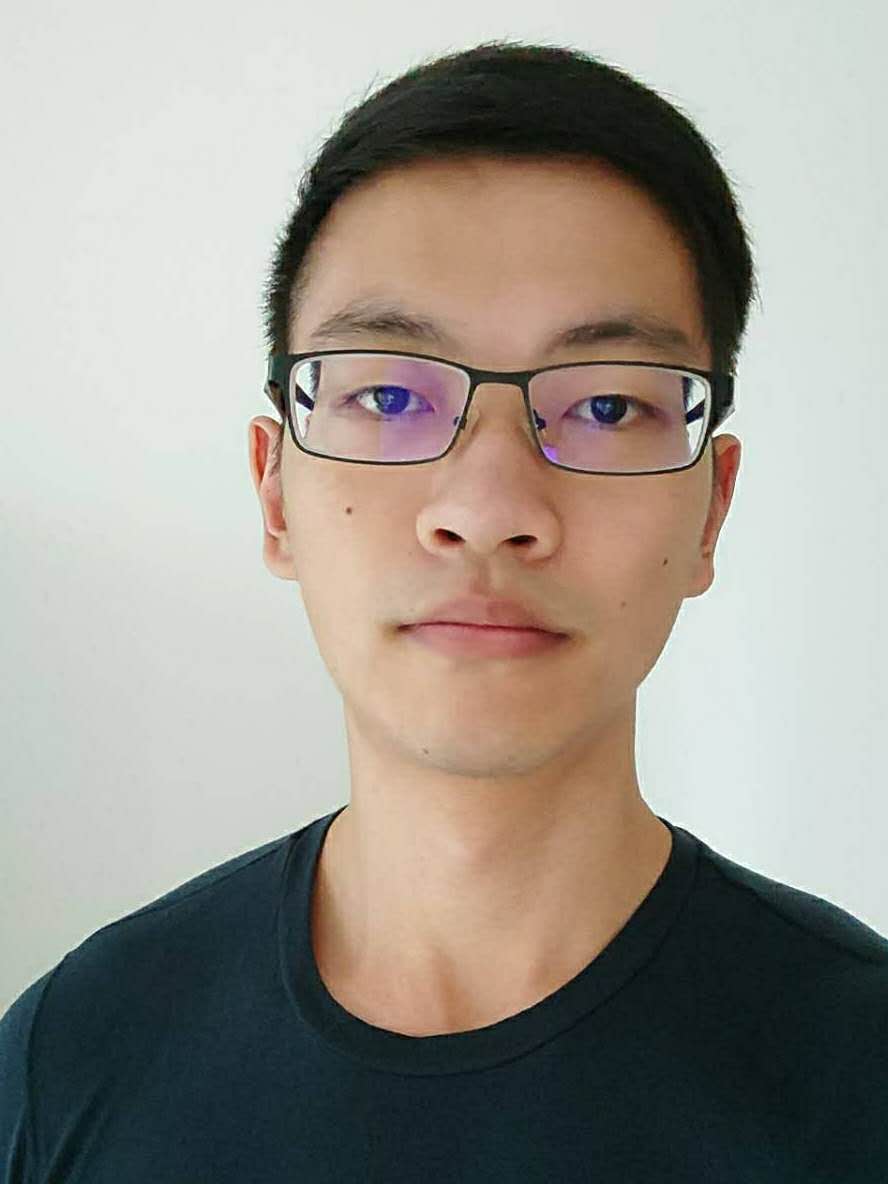}}]{Mengmou Li} is currently a Specially Appointed Assistant Professor with the Department of Systems and Control Engineering at Tokyo Institute of Technology, Japan. He received his B.S. degree in Physics from Zhejiang University, China, in 2016 and his Ph.D. degree in Electrical and Electronic Engineering from the University of Hong Kong in 2020. From February 2021 to July 2022, he served as a research associate with the Control Group at the University of Cambridge. His research interests include power systems, optimization, and robust control.
\end{IEEEbiography}

\begin{IEEEbiography}[{\includegraphics[width=1in,height=1.25in,clip,keepaspectratio]{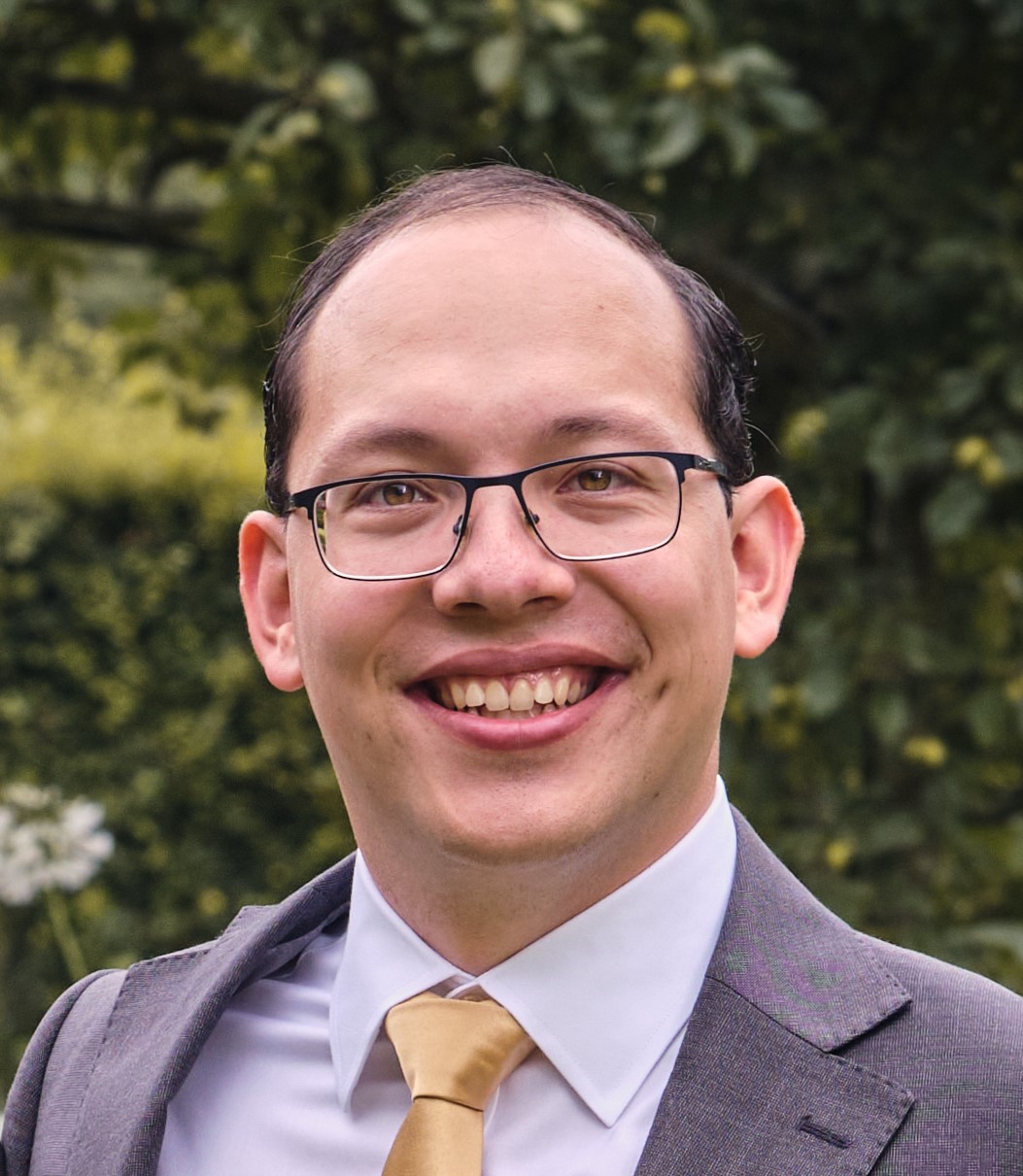}}]{Jeremy~D.~Watson} received the B.E. degree (First Class Hons.) in electrical engineering from the University of Canterbury, Christchurch, New Zealand, in 2015, and the Ph.D. degree in engineering from the University of Cambridge, Cambridge, United Kingdom, in 2021, where he was a Research Associate with the Control Group, Department of Engineering. He is currently a lecturer at the University of Canterbury, Christchurch, New Zealand. His research interests include control, analysis, and optimization of power networks, focusing especially on hybrid AC/DC networks and microgrids.
\end{IEEEbiography}

\begin{IEEEbiography}[{\includegraphics[width=1in,height=1.25in,clip,keepaspectratio]{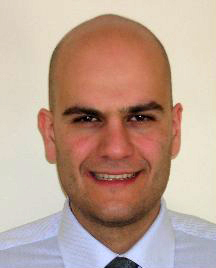}}]{Ioannis~Lestas} received the B.A. (Starred First) and M.Eng. (Distinction) degrees in Electrical and In formation Sciences and the Ph.D. in control theory from the University of Cambridge (Trinity College) in 2002 and 2007, respectively. His doctoral work was performed as a Gates Scholar. He has been a Junior Research Fellow of Clare College, University of Cambridge and he was awarded a five year Royal Academy of Engineering research fellowship. He is also the recipient of a five year ERC starting grant. He is currently a Professor at the University of Cambridge, Department of Engineering. His research interests include analysis and control of large scale networks with applications in power systems and smart grids.
\end{IEEEbiography}

%\begin{IEEEbiography}{Ioannis Lestas}
%Biography.
%\end{IEEEbiography}

%\newpage


%\vfill


\end{document}